\documentclass[10pt,reqno,fleqn]{amsart}

\usepackage[utf8]{inputenc}		
\usepackage[english]{babel}		

\usepackage{amsmath}				
\usepackage{mathtools}			
\usepackage{amsthm}				
\usepackage{amssymb}				
\usepackage[nointlimits]{esint}	
\usepackage{mathrsfs} 			

\usepackage{graphicx}			
\usepackage{epsfig}				
\usepackage{tikz}				
\usepackage{xxcolor}				
\usepackage{wrapfig}
\usepackage{caption}

\usepackage{enumitem}			
\usepackage{ifdraft} 			

\usepackage[]{todonotes}			

\usepackage[linktocpage=true,colorlinks=true,linkcolor=blue,citecolor=green]{hyperref}
\usepackage[
	backend=biber,style=alphabetic, natbib=false, mcite=false, casechanger=auto,
	sorting=nyt, sortcites=false, pluralothers=true, maxnames=4, minnames=1,
	backref=false, arxiv=abs, 
	date=year,
	isbn=false, url=false, doi=true, eprint=true, related=false,
	giveninits=true
]{biblatex}
\renewbibmacro{in:}{}	

\usepackage[hmarginratio=1:1,top=32mm,columnsep=20pt]{geometry} 

\numberwithin{equation}{section}				

\theoremstyle{plain}							
\newtheorem{theorem}{Theorem}[section]		
\newtheorem*{theorem*}{Theorem}				
\newtheorem*{maintheorem*}{Main Result}
\newtheorem{lemma}[theorem]{Lemma}			
\newtheorem{proposition}[theorem]{Proposition}
\newtheorem{corollary}[theorem]{Corollary}

\theoremstyle{definition}					
\newtheorem{definition}[theorem]{Definition}

\theoremstyle{remark}						
\newtheorem{remark}[theorem]{Remark}

\theoremstyle{definition}					
\providecommand{\customtheoremname}{}			
\newtheorem{placeholderthhm}{\customtheoremname}
\newcommand{\newcustomtheorem}[2]{%
  \newenvironment{#1}[1]						
  {											
   \renewcommand\customtheroremname{#2}
   \renewcommand\theplaceholderthhm{##1}	
   \placeholderthhm							
  }											
  {\endplaceholderthhm}						
}
\newcustomtheorem{customass}{Assumption}		
\newcustomtheorem{customdef}{Definition}

\DeclarePairedDelimiter{\Vpair}{\|}{\|} 			
\DeclarePairedDelimiter{\vpair}{|}{|} 			
\DeclarePairedDelimiter{\bpair}{[}{]} 			
\DeclarePairedDelimiter{\ppair}{(}{)} 			
\DeclarePairedDelimiter{\Bpair}{\{}{\}} 			
\DeclarePairedDelimiter{\ceil}{\lceil}{\rceil} 	
\DeclarePairedDelimiterX{\cond}[2]{[}{]}{#1  \,\delimsize\vert\,  \mathopen{}#2}

\newcommand\mathvshift[2]{\raisebox{#1}{$\displaystyle{#2}$}}
\newcommand\strutb{\vphantom{\bigm|}}

\newcommand\norminner[2][0pt]{\mathchoice{\mathvshift{#1}{\strutb\mathopen{}#2}}{\smash{#2}}{#2}{#2}}
\newcommand\bracketinner[2][0pt]{\mathchoice{\mathvshift{#1}{\strutb\mathopen{}#2}}{#2}{#2}{#2}}
\ifoptiondraft{
	\newcommand\norm[2][0pt]{\Vpair*{#2}} 	
	\newcommand\abs[2][0pt]{\vpair*{#2}} 		
	\newcommand\pp[2][0pt]{\ppair*{#2}} 		
	\newcommand\bp[2][0pt]{\bpair*{#2}} 		
	\newcommand\Bp[2][0pt]{\Bpair*{#2}} 		
}{
	\newcommand\norm[2][0pt]{\Vpair*{\norminner[#1]{#2}}}	 	
	\newcommand\abs[2][0pt]{\vpair*{\norminner[#1]{#2}}} 		
	\newcommand\pp[2][0pt]{\ppair*{\bracketinner[#1]{#2}}} 	
	\newcommand\bp[2][0pt]{\bpair*{\bracketinner[#1]{#2}}} 	
	\newcommand\Bp[2][0pt]{\Bpair*{\bracketinner[#1]{#2}}}		
}
\newcommand\delimstyle[2]{\setlength{\delimitershortfall}{#1}\delimiterfactor=#2}
\delimstyle{5pt}{850}

\newcommand{\ov}{\overline}
\newcommand{\wt}{\widetilde}

\newcommand{\vphi}{\varphi}

\newcommand{\data}{\mathrm{data}}
\newcommand{\vrectangle}{{\setlength{\fboxsep}{0pt}\fbox{\phantom{t}}}}

\newcommand{\N}{\mathbb{N}}
\newcommand{\Nz}{\mathbb{N}_0}
\newcommand{\Z}{\mathbb{Z}}

\newcommand{\R}{\mathbb{R}}
\newcommand{\Zt}{\Z^2}

\newcommand{\Rt}{{\R^{2}}}
\newcommand{\Rtt}{{\R^{2\times 2}}}
\newcommand{\Rttsym}{{\R^{2\times 2}_{\mathrm{sym}}}}
\newcommand{\Sone}{{\mathbb{S}^1}}

\newcommand{\pat}{\partial_t}

\newcommand{\diam}{\operatorname{diam}}
\newcommand{\dist}{\operatorname{dist}}

\newcommand{\charfun}[1][_]{{\if_#1 \chi \else \chi_{#1} \fi}}			
\newcommand{\comp}[1]{{{#1}^c}}										

\newcommand{\Idt}{\operatorname{Id}_{2\times 2}}
\newcommand{\tr}{\operatorname{tr}}		
\newcommand{\supp}{\operatorname{supp}}

\newcommand{\Pm}{{\mathbb{P}}}											
\newcommand{\PM}[2][_]{ {\if_#1 \Pm\bp{#2} \else \Pm\cond*{#2}{#1} \fi} }	
\newcommand{\F}{\mathcal{F}}												
\newcommand{\G}{\mathcal{G}}												

\newcommand{\vmin}{v_{\infty}}
\newcommand{\vmineff}{v_{\mathrm{min}}}
\newcommand{\vmax}[1][_]{{\if_#1 v_{\mathrm{max}} \else v_{\mathrm{max},#1} \fi}}
\newcommand{\cw}{c_w}

\newcommand{\ratebox}{{\alpha_{\scalebox{0.5}{$\uparrow$}}}}
\newcommand{\Cbox}{{C_{\scalebox{0.5}{$\uparrow$}}}}

\newcommand{\Cgrid}{{C_{\scalebox{0.5}{$\square$}}}}

\newcommand{\restcomp}{Q^{\oslash}}		

\newcommand{\Eup}[2][_]{{E^{\if_#1 \scalebox{0.6}{$\uparrow$} \else \scalebox{0.6}{$\uparrow$},#1 \fi}_{#2}}}
\newcommand{\Er}[2][_]{{E^{\if_#1 \scalebox{0.6}{$\rightarrow$} \else \scalebox{0.6}{$\rightarrow$},#1 \fi}_{#2}}}
\newcommand{\Ed}[2][_]{{E^{\if_#1 \scalebox{0.6}{$\downarrow$} \else \scalebox{0.6}{$\downarrow$},#1 \fi}_{#2}}}
\newcommand{\El}[2][_]{{E^{\if_#1 \scalebox{0.6}{$\leftarrow$} \else \scalebox{0.6}{$\leftarrow$},#1 \fi}_{#2}}}
\newcommand{\Qup}[1]{{Q^{\scalebox{0.5}{$\uparrow$}}_{#1}}}
\newcommand{\Qr}[1]{{Q^{\scalebox{0.6}{$\rightarrow$}}_{#1}}}
\newcommand{\Qd}[1]{{Q^{\scalebox{0.5}{$\downarrow$}}_{#1}}}
\newcommand{\Ql}[1]{{Q^{\scalebox{0.6}{$\leftarrow$}}_{#1}}}
\newcommand{\Esquare}[2][_]{{E^{\if_#1 \scalebox{0.5}{$\square$} \else \scalebox{0.5}{$\square$},#1 \fi}_{#2}}}
\newcommand{\Egrid}[2][_]{{E^{\if_#1 \scalebox{0.65}{$\boxplus$} \else \scalebox{0.65}{$\boxplus$},#1 \fi}_{#2}}}

\newcommand{\Espeed}[2][_]{{E^{\if_#1 \scalebox{0.65}{$\circledcirc$} \else \scalebox{0.65}{$\circledcirc$},#1 \fi}_{#2}}}
\newcommand{\ratespeed}{{\alpha_{\scalebox{0.5}{$\circledcirc$}}}}
\newcommand{\Cspeed}{{C_{\scalebox{0.5}{$\circledcirc$}}}}
\newcommand{\cspeed}{{c_{\scalebox{0.5}{$\circledcirc$}}}}

\newcommand{\subJ}[1][]{\mathcal{J}^{2,-}_{#1}}			

\newcommand{\reach}[3][]{{{\mathscr{R}}_{#2}^{#1}\pp{#3}}}
\newcommand{\uls}[2][_]{\if_#1 u_{#2} \else u_{#2, #1} \fi}	
\newcommand{\ur}[2][_]{\if_#1 u_{#2}^{\mathscr{R}} \else u_{#2, #1}^{\mathscr{R}} \fi}		

\newcommand{\kvert}{k_{\scalebox{0.5}{$\parallel$}}}
\newcommand{\khor}{k_{\scalebox{0.6}{$=$}}}
\newcommand{\rup}{r_{\scalebox{0.5}{$\parallel$}}}
\newcommand{\rbase}{r_{\scalebox{0.6}{$=$}}}
\newcommand{\xup}[2]{{x^{\scalebox{0.5}{$\uparrow$}}_{#1,#2}}}

\newcommand{\xr}[2]{{x^{\scalebox{0.6}{$\rightarrow$}}_{#1,#2}}}
\newcommand{\xl}[2]{{x^{\scalebox{0.6}{$\leftarrow$}}_{#1,#2}}}
\newcommand{\xrr}[2]{{x_{\scalebox{0.6}{$\rightarrow$}}}}
\newcommand{\xll}[2]{{x_{\scalebox{0.6}{$\leftarrow$}}}}
\newcommand{\jup}[1]{{j_{#1\scalebox{0.5}{$\uparrow$}}}}
\newcommand{\jr}[1]{{j_{#1\scalebox{0.6}{$\rightarrow$}}}}
\newcommand{\jl}[1]{{j_{#1\scalebox{0.6}{$\leftarrow$}}}}

\newcommand{\Cv}{C_v}

\newcommand{\hsubs}[1][_]{\if_#1 \subset_h \else \subset_{#1} \fi}

\newcommand{\inner}[2][]{{\operatorname{int}_{#1}\pp{#2}}}
\newcommand{\Cmeas}{c_{\F}}
\newcommand{\cstable}{c_s}

\newcommand{\yup}{{y^{\scalebox{0.5}{$\uparrow$}}}}
\newcommand{\yd}{{y^{\scalebox{0.5}{$\downarrow$}}}}
\newcommand{\yr}{{y^{\scalebox{0.6}{$\rightarrow$}}}}
\newcommand{\yl}{{y^{\scalebox{0.6}{$\leftarrow$}}}}

\newcommand{\CG}[1][]{C_{\G}}

\newcommand{\vhom}[1][]{\ov{v}_{#1}}

\graphicspath{ {./MCF_box_images/} }
\counterwithin{figure}{section}

\begin{document}

\title[A criterion for ballistic behavior of geometric motions]{A sufficient criterion for ballistic behavior of forced geometric motions in random environments}

\author{Julian Fischer}
\address{Institute of Science and Technology Austria (ISTA), Am~Campus~1, 
3400~Klosterneuburg, Austria}
\author{Jonas Ingmanns}
\address{Universit\'e Libre de Bruxelles (ULB), Boulevard~du~Triomphe, 
1050~Bruxelles, Belgium}

\begin{abstract}
We introduce a sufficient criterion for ballistic motion of interface evolutions in random environments. 
More precisely, we consider interfaces moving by forced mean curvature flow in dimension $d=2$ with a spatially fluctuating (possibly negative) driving force. Such models arise e.g. from the motions of dislocation lines or magnetic domain walls in a material with impurities. We show that the evolution admits an `effective minimum speed on large scales' if for a sample of the random environment in some long and thin box there is a sufficiently large probability that an interface initialized in the bottom of the box evolves to reach its top.
The concept of an `effective minimum speed on large scales' was introduced and is required for the quantitative stochastic homogenization result in the companion paper \cite{FI25}. 
We show that our criterion and hence the requirements for the homogenization result in \cite{FI25} are satisfied for compactly supported obstacle shapes distributed according to a Poisson point process with sufficiently low intensity.
\end{abstract}

\thanks{
This project has received funding from the European Research Council (ERC) under the European Union's Horizon 2020 research and innovation programme (grant agreement No 948819)
\smash{
\begin{tabular}{@{}c@{}}\includegraphics[width=6ex]{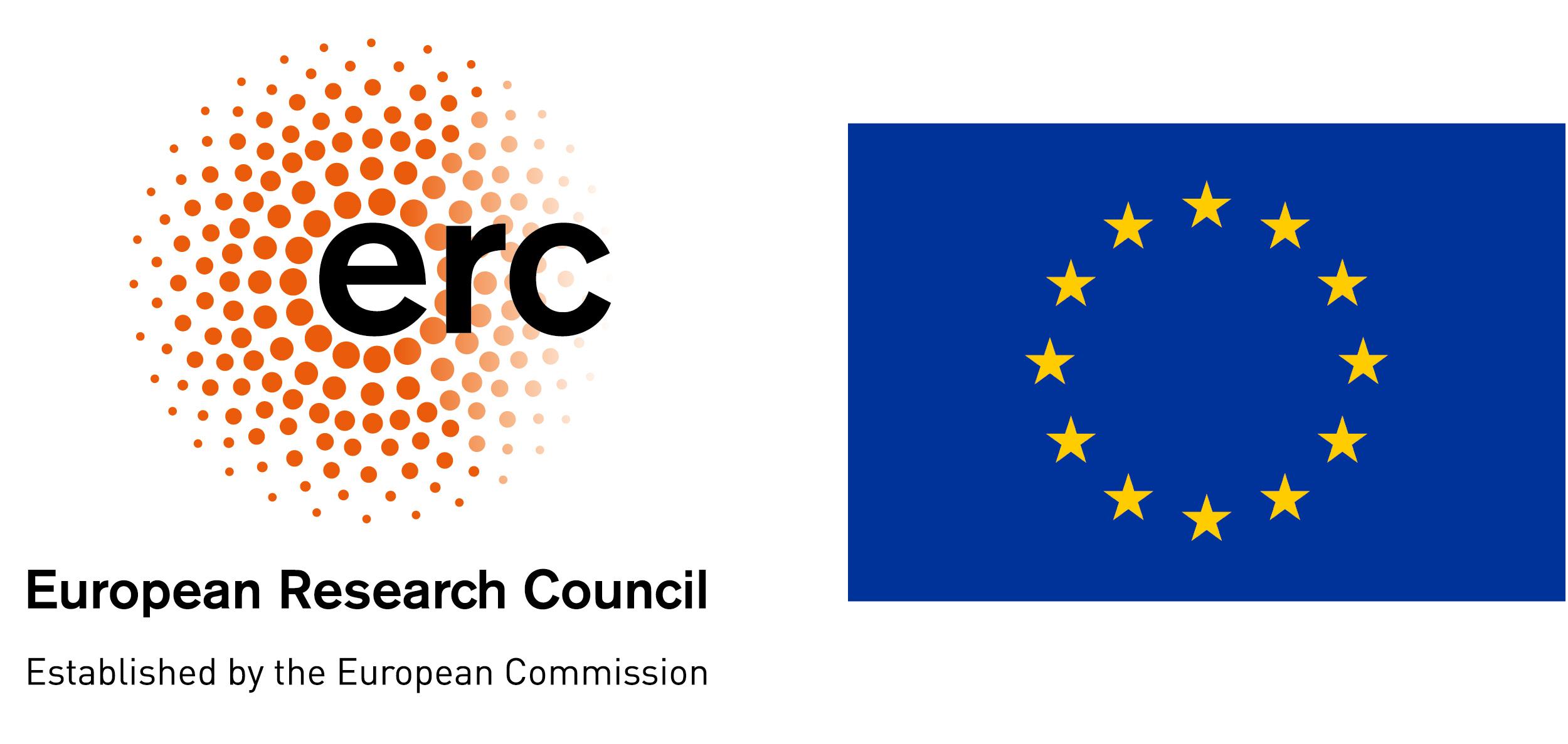}\end{tabular}
}
}

\maketitle
{\small{\bfseries \href{https://mathscinet.ams.org/mathscinet/msc/msc2020.html}{MSC2020 subject classification}:} 60K35, 53E10, 35R60, 35B27, 35F21}


\setcounter{tocdepth}{1}		
\tableofcontents
\setcounter{tocdepth}{3}		


\section{Introduction}\label{s_intro}
We consider the question of when ballistic motion occurs for interface evolutions in random environments, that is, whether on large time scales the interface propagates with positive velocity.
In the present work, we introduce a sufficient criterion for ballistic motion in dimension $d=2$ for arbitrary random environments with a finite range of dependence.
While we expect that our criterion and methods can be extended to general interface evolutions with a comparison principle, we focus on interfaces evolving by forced mean curvature flow through a random environment:
The normal velocity $V$ of an interface at $x\in\Rt$ is given in terms of the normal $\nu\in\Sone$ and the (mean) curvature $H$ by
\begin{equation}\label{intro_eq_normal}
	V= a\pp{x, \nu}H + F\pp{x, \nu},
\end{equation}
with the curvature coefficient field $a=\sigma^2\geq 0$ and the forcing field
\begin{equation*}
	F(x,\nu)=F_{uni}(\nu)-F_{obst}(x,\nu).
\end{equation*}
Here, $F_{uni}$ can be thought of as a deterministic, external uniform forcing while $F_{obst}$ encodes local obstacles.
We assume that $\sigma,F\in C^{0,1}(\Rt\times\Sone)$ are sampled from a stationary probability measure $\Pm$ with finite range of dependence, see Section~\ref{s_assump} for a detailed description of the assumptions.

For a given random environment in the form of  $(\sigma,F_{obst})$, qualitatively very different behavior might emerge depending on the strength of the uniform forcing $F_{uni}$.
For \eqref{intro_eq_normal} and related interface motions such as the quenched Edwards-Wilkinson model, the occurrence of pinning and depinning phenomena has been established both for periodic obstacle fields \cite{DiYi06} and for certain stationary random obstacle fields with a finite range of dependence \cite{CDL10, DoSch12, BoTe14, DoSch17, DJS21, DJS22}.
In short, for prototypical obstacle fields, there are critical forcings $0<F_* \leq F^*$ such that  
\begin{itemize}
	\item for $F_{uni}< F_*$, any interface gets pinned, meaning that its evolution is bounded by a stationary surface (pinning regime),
	\item for $F_{uni}>F^*$, for increasingly large initial interfaces the probability that their evolution is ballistic converges to one, meaning that the main front covers a strictly positive average distance per time unit
	(ballistic regime).
\end{itemize}
Both $F_*=F^*$ \cite{BoTe14} or $F_*<F^*$ \cite{DoSch17} are possible, depending on the specific model. 
In the case of $F_*<F^*$, there is an intermediate sub-ballistic regime in which the interface propagates indefinitely but the limit velocity is zero.

Pinning and depinning can be observed for a range of physical phenomena in material science such as dislocation lines, magnetic domain walls, or grain boundaries \cite{IoVi87, Kar98, BraNa04}.
The motions of these lines or interfaces is driven by the competition between three effects:
A preference for the minimization of the interface area, an external, uniform forcing pushing the interface forwards, and material impurities which may act as local obstacles to the interface motion. 
Forced mean curvature flow as in \eqref{intro_eq_normal} is used to model these phenomena with the strength of the effects encoded in $a$, $F_{uni}$, and $F_{obst}$.

Most results for random environments have been established for explicitly modeled obstacle fields, usually obstacles distributed according to a Poisson point process or a discrete analogue. 
We introduce a criterion applicable to arbitrary random environments which implies ballistic behavior of the interface motion.
In order to show that this criterion is satisfied, it is sufficient to check statistical properties of a sample from the random fields in a box not necessarily much larger than the range of dependence.
For a discrete version of the quenched Edwards-Wilkinson model, \cite{BoTe14} also provides a criterion for ballistic behavior given arbitrary random environments while allowing for slowly decaying long-range correlations rather than the finite range of dependence which we require.
However, their criterion already requires control on all scales for the probability that there is a pinned solution.
Moreover, our criterion implies not just ballistic behavior but the existence of an almost guaranteed `effective minimum speed of propagation', a notion that has been introduced and is required for the homogenization result established in the companion paper \cite{FI25}.

As an application, we show that our criterion is satisfied for Poisson distributed obstacles with low enough intensity compared to the uniform forcing.
In particular, having thus established the existence of an almost guaranteed `effective minimum speed of propagation', we obtain the homogenization result \cite[Theorem~1.1]{FI25} announced in the companion paper.
Ballistic motion in this setting has already been established in \cite{CDO22}: 
Starting with a flat interface, the average height of the main front above the original line grows with positive speed.
The notion of an almost guaranteed `effective minimum speed of propagation' is stronger: 
On large scales, we show that with high probability any sufficiently regular interface expands at least according to some minimum speed in any direction -- up to relatively small enclosures left behind the main front -- and that this holds not just for the initial interface but at any stage of the evolution even when the initial regularity assumption might not be satisfied anymore.

As a special case, \eqref{intro_eq_normal} includes percolation models:
For vanishing curvature coefficient $a\equiv 0$ and the forcing given by $F(x,e)=e\cdot g^{-1}(x)e$ on $M=\Bp{x\in\Rt\,;\,F(x)>0}$ with $g(x)\in\Rttsym$ positive definite for all $x\in M$, the evolution corresponds to Riemannian first passage percolation \cite{LaGW10} on $M$ equipped with the Riemannian metric $g$: 
The time by which the set $\Bp{x}$ evolving by \eqref{intro_eq_normal} -- see Definition \ref{not_def_reach} for the precise interpretation -- corresponds to the Riemannian distance between $x$ and $y$, that is the infimum over the Riemannian length of all paths in $M$ connecting $x$ and $y$. 
Our work can hence be viewed as a generalization or extension of results from percolation theory which show that the chemical distance in infinite clusters is comparable to the Euclidean distance with high probability, see \cite{AnPi96}. 
However, in contrast to those works, our notions and arguments are build for general set evolutions with a comparison principle. 

\begin{figure}[h]
\centering
\includegraphics[width=0.5\textwidth]{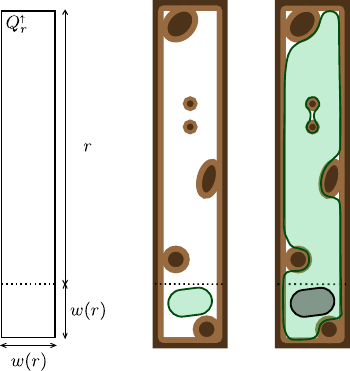}
\caption[An illustration of the event {${\Eup[v]{r}}$} introduced for our box criterion]{	An illustration of the event $\Eup[v]{r}$ introduced for our box criterion.
		Left: The dimensions of the box $\Qup{r}$ and its `starting area'.
		Center: A sample of the obstacle field (depicted in brown) inside the box with the artificial obstacle surrounding the boundary (also depicted in brown) since the coefficient field is restricted to the box. 
		For the event $\Eup[v]{r}$ to hold, there has to be a stable set $S$ in the starting area of the box (depicted in green).
		Right: After time $t=v^{-1}r$, the evolution $\reach{t}{S}$ (depicted in green) of the set $S$ is required to have reached the top of the box.
}
\label{intro_fig_box}
\end{figure}

We now state an informal version of our \emph{box criterion}, which roughly speaking requires that it is sufficiently likely for an interface starting in the bottom of a long but thin box sample of the environment to reach the top of the box.
The threshold for this probability depends on the dimensions of the box and can be arbitrarily small, see Assumption~\ref{sett_aP_box} in Section~\ref{s_assump} for the rigorous formulation of the criterion.
In the special case of first passage percolation, this corresponds to a lower bound on the crossing probabilities in those boxes and covers the whole supercritical phase for suitable parameters.
Using this criterion, we will later obtain improved probability bounds for the same event in larger boxes with the width $w(r)$ of the boxes scaling sublinearly compared to the length $r$.
Since this sublinear scaling is central for obtaining a positive speed in the limit, we incorporate it into the criterion to streamline the notation.
\begin{enumerate}[label=`(P$\,\vrectangle$)', left=0pt]
	\item	\label{sett_aP_pseudoBox}
			\textbf{Box criterion.}
			We assume that there exist a velocity $v_0>0$, a width scaling exponent $q>1$, and a box length $r_0>C(q)$ such that within the thin boxes $\Qup{r}=[0,w(r)]\times [-w(r),r]$ with width $w(r)\sim r^{\frac{1}{q}}$, as in Figure~\ref{intro_fig_box}, it holds that
			\begin{align}
				&\qquad\Pm\Big[\big\lbrace
					\text{There is a `stable' set in the bottom $[0,w(r)]\times[-w(r),0]$ of $\Qup{r}$,} \notag\\
				&\qquad\qquad\qquad	
					\text{which evolving `within $\Qup{r}$' by \eqref{intro_eq_normal}}\notag\\
				&\qquad\qquad\qquad
					\text{reaches the top $[0,w(r)]\times\Bp{r}$  by the time $t=v_0^{-1}r$}
				\big\rbrace\Big]\notag\\
				&\qquad\quad	\geq		1-p_0,\label{sett_eq_pseudoBox}
			\end{align}
			for all $r\in [r_0,r_1)$ with the acceptable probability of failure $p_0=p_0(r_0,w(\cdot))\geq 0$ depending on the dimensions of the boxes.
			The upper end $r_1$ of the range for the box lengths is chosen such that a rotation of $\Qup{r_0}$ fits into $\Qup{r_1}$, see Figure~\ref{sett_fig_r1}, meaning that $w(r_1)\approx r_0$. 
			Here, `stable' means that the set is non-shrinking and `within $\Qup{r}$' means that the coefficient fields are modified to restrict the evolution to the box.
			In addition, we assume that the same is true for the accordingly rotated boxes $\Qr{r},\Qd{r}, \Ql{r}$.	
\end{enumerate}
The following informal version of our main result is a precise statement up to replacing \ref{sett_aP_pseudoBox} with the rigorous formulation stated in Assumption~\ref{sett_aP_box} and the informal version \ref{sett_aP_pseudoStar} of the guarantee for an `effective minimum speed of propagation' from \cite{FI25}, which is stated below the result, with its rigorous formulation \ref{veff_aP_star} given in Section~\ref{subs_sett_veff}.

\begin{theorem}[Informal main result]\label{intro_thm_informal}
Let the coefficient fields $a=\sigma^2$ and $F$ be sampled from a stationary probability measure $\Pm$ with range of dependence one and let $\sigma,F\in C^{0,1}(\Rt\times\Sone)$ be uniformly bounded and Lipschitz. 
Assume that $\Pm$ satisfies the box criterion \ref{sett_aP_pseudoBox}. 
Then the evolution \eqref{intro_eq_normal} with the coefficient field $(a,F)$ admits an almost guaranteed `effective minimum speed' $\vmineff=\vmineff(r_0,w(\cdot),v_0)>0$ in the sense of \ref{sett_aP_pseudoStar} from \cite{FI25}, see below.
In particular, this evolution is in the ballistic regime.
\end{theorem}

We now provide the aforementioned informal version of the almost guarantee for an `effective minimum speed'. 
We will state the rigorous formulation and the required definitions of the respective notions from \cite{FI25} later in Section~\ref{subs_sett_veff}.
\begin{enumerate}[label=`(P$\circledcirc$)', left=0pt, resume]
	\item	\label{sett_aP_pseudoStar}\textbf{Almost guaranteed effective minimum speed}.
			There exist $\vmineff,C>0$ and a rate $\ratespeed>0$ such that for all $h\geq h_0$ and $M\subset\Rt$
			\begin{align*}
				&\qquad\Pm\Big[\big\lbrace
					\text{In $M$, on the scale $h$, every `fat' set has a `stable approximation',} \\
				&\qquad\qquad\qquad\quad			
					\text{which propagates with an `effective minimum speed' of $\vmineff$}
				\big\rbrace\Big]\\
				&\qquad\quad	\geq		1-C\pp{1+\diam\pp{M}}^2\exp\pp{-C^{-1}h^\ratespeed},
			\end{align*}
			where `fat' sets $S\subset\Rt$ are unions of sufficiently large balls and `stable approximations' are non-shrinking interior approximations $S_h\subset S$.
			An `effective minimum speed' on the scale $h$ means that for any $t\geq 0$ the evolved set $\reach{t}{S_h}$ after another time step of $\vmineff^{-1}h$ has advanced by at least $h$ in all directions while allowing for enclosures of width smaller than $2h$ to be left behind the main front.
\end{enumerate}

Such an almost guaranteed effective minimum speed allows us to apply the homogenization results from the companion paper under the additional assumptions that the probability measure is isotropic and that the curvature coefficient field is constant.

\begin{corollary}[{Application of \cite[Theorem~4.1 and Theorem~4.2]{FI25}}] 
Let the coefficient field $F$ be sampled from a stationary probability measure $\Pm$ with range of dependence one and $F\in C^{0,1}(\Rt\times\Sone)$ uniformly bounded and Lipschitz. 
Let $\Pm$ be isotropic, that is $\operatorname{law} F(\cdot,\cdot) = \operatorname{law} F(R\cdot,R\cdot)$ for all rotations $R\in \operatorname{SO}(2)$.
Assume that $\Pm$ satisfies the box criterion \ref{sett_aP_pseudoBox} for \eqref{intro_eq_normal} with the coefficient field $(a,F_0+F)$ for some constants $a\geq 0$ and $F_0\in\R$. 
Then for Lebesgue-almost every $F_{uni}>F_0$ there exists $\vhom=\vhom(a,\Pm, F_{uni})>0$ such that on large scales the evolution \eqref{intro_eq_normal} with $(a,F_{uni}+F(x))$ can be approximated by a constant-speed motion
\begin{equation*}
	V= \vhom
\end{equation*}
in terms of arrival times in target areas which are small with respect to the large scale but large compared to the range of dependence, see \cite{FI25}.
\end{corollary}

As an application of the main result, we show that \ref{sett_aP_pseudoBox} and hence \ref{sett_aP_pseudoStar} is satisfied for a compactly supported obstacle shape distributed according to a Poisson point process with sufficiently low intensity compared to the uniform forcing.
We obtain the same scaling relationship (``Taylor scaling'') between the intensity and $F^*$, the uniform forcing marking the transition to the ballistic regime, as in \cite{CDO22}.
Note however that \ref{sett_aP_pseudoBox} is stronger than the ballistic property shown in \cite{CDO22}.

\begin{theorem}[Almost guaranteed effective minimum speed for Poisson obstacles]\label{intro_thm_Poisson}
Consider Poisson-distributed obstacles in dimension $d=2$ represented by 
\begin{equation*}
	F_{obst}\colon\R^2\rightarrow \R,\quad F_{obst}(x) \coloneqq \max_{i\in\N}\vphi(\abs{x-X_i})
\end{equation*}
for
\begin{itemize}
\item	a Poisson point process $\Bp{X_i}_{i\in\N}\subset\Rt$ with intensity $\rho>0$, and
\item	a Lipschitz continuous obstacle shape $\vphi\in C^{0,1}(\R)$ with support $\operatorname{\supp}\vphi\subset\ov{B}_1$.
\end{itemize}
There exist $\rho_0>0$ and $C_0>0$ independent of the obstacle shape $\vphi$ such that -- given any intensity $\rho\in(0, \rho_0]$ -- for  every $F_{uni}\geq C_0\sqrt{\rho}$ the interface motion given by
\begin{equation*}
	\mathrm{Normal~velocity} = \mathrm{curvature} + F_{uni}-F_{obst}(x) 
	\qquad\text{at $x\in \R^2$}
\end{equation*}
satisfies the box criterion \ref{sett_aP_box} for some $v_0=v_0(F_{uni})>0$, $r_0= r_0(F_{uni})>1$, and $w(r)= \pp{\frac{r}{\cw}}^{\frac{1}{q}}$ with $q>1$ and $\cw=\cw(F_{uni})>0$.
\end{theorem}

Combining this result with Theorem~\ref{intro_thm_informal} allows the application of \cite[Theorem~4.1 and Theorem~4.2]{FI25} to sufficiently sparse Poisson obstacles and hence yields the quantitative homogenization result announced in \cite[Theorem~1.1]{FI25}.

\subsection{Structure of the paper}
In Section~\ref{s_lit}, we discuss the related literature, in particular with respect to results on pinning-depinning phenomena.
In Section~\ref{s_assump}, we state our assumptions including the box criterion and give precise definitions of the notions used throughout the present work such as set evolutions, restricted coefficient fields, stable sets,  stable approximations and an effective minimum speed of propagation.
In Section~\ref{s_str}, we provide the formal version of our main result and describe the strategy of the proof. 
In Section~\ref{s_pre}, we gather preliminary results on the basic properties of the set evolutions which will be used in the main part of the proof.
In Section~\ref{s_vmin} and Section~\ref{s_veff}, the main part of the proof is presented, see Section~\ref{subs_str_str} for a more detailed description.
In Section~\ref{s_poisson}, we will apply our box criterion for obstacles distributed according to a Poisson point process, establishing Theorem~\ref{intro_thm_Poisson}.

\subsection{Statement on AI Use}
Large Language Models were only used as a first additional proofreading pass.
All issues found during this process were fixed without the use of AI.
Outside of this single proofreading pass, no generative AI or AI tools were used in the research or writing process.

\section{Related literature}\label{s_lit}
The present work is a companion paper to \cite{FI25}, in which the authors established a stochastic homogenization result for sets evolving by forced mean curvature flow through random obstacle fields assuming an `effective minimum speed', see Section~\ref{subs_sett_veff}.
One of the major motivations of the present work is to provide an easy-to-check, tangible criterion in dimension $d=2$ for the quite general assumption of an `almost guaranteed effective minimum speed' introduced in \cite{FI25}.
For a more detailed review of related literature on first passage percolation \cite{Kes93, Ale97, BKS03, LaGW10, Zha10, BSS26}, which corresponds to the situation without surface tension ($a\equiv 0$), and homogenization in a periodic setting \cite{LiSo05a, DKY08, CaMo14, GaKi19} as well as a random setting \cite{ArCa18, Fel19} we refer to \cite[Section~2]{FI25}.
Here, we focus on results for pinning and depinning phenomena and some more related results from percolation theory.

\subsection{Physical motivation}
The competition between some form of surface tension, a uniform driving force and local obstacles or impurities as represented by \eqref{intro_eq_normal} governs a range of physical phenomena such as the random field Ising model \cite{KoLe85, Kar98}, the evolution of domain boundaries in magnetic materials \cite{IoVi87, Kar98, BraNa04} or the evolution of dislocation lines in a slip plane \cite{IoVi87, BraNa04}.
In addition, \eqref{intro_eq_normal} is used as a simplification of more complex models, for example immiscible-fluid displacement in porous media \cite{KoLe85, Kar98}.
A key phenomenon observed in these settings is a pinning-depinning transition when increasing the uniform driving force:
The interface remains pinned even for a small positive driving force. 
Only when increasing the forcing past a critical threshold will the interface propagate at a positive effective speed \cite{Kar98, BraNa04}.
Such behavior is also investigated with respect to charge density waves \cite{FuLe78, IoVi87, Kar98, BraNa04}, contact lines \cite{Ge85, Kar98}, Wigner crystals \cite{NaSch00, BraNa04} and vertex (or flux) lines in type-II superconductors \cite{Kar98, NaSch00, BraNa04}.

\subsection{Pinning and depinning}\label{subs_lit_pinn}
Pinning and depinning effects have been extensively studied for the \emph{quenched Edwards-Wilkinson (QEW) model}, which heuristically approximates the evolution by forced mean curvature flow as in \eqref{intro_eq_normal} if the interface can be represented as a graph with small slope.
For a constant, isotropic curvature coefficient ($a\equiv 1$ in \eqref{intro_eq_normal}) and a forcing $F(x)=F_{uni}-F_{obst}(x)$, the height $h$ of the approximating interface in the QEW model is given by 
\begin{align}\label{lit_QEW}
	\partial_t h(y,t) &= \Delta_y h(y,t) + F_{uni} - F(y, h(y))
	&&\text{on }\R^{d-1}.
\end{align}
However, there are some fundamental differences compared to forced mean curvature flow: 
This model enforces the graph structure while an interface evolving by forced mean curvature flow can ``wrap around'' obstacles, leaving enclosures behind the main front.

The main objects of study in the related mathematical literature on pinning and depinning phenomena are the critical thresholds $F_*,F^*$ for the uniform forcing $F_{uni}$ introduced at the start of Section~\ref{s_intro}, below and above which the behavior of the interface evolution changes qualitatively. 
For the QEW model, these can be defined as follows.
\begin{itemize}
\item		$F_*$ is the supremum of all $F\geq 0$ such that the height $h$ of the interface remains pointwise bounded for all $F_{uni}\leq F$, that is, there exists a stationary supersolution to \eqref{lit_QEW} for all $F_{uni}\leq F$ (pinning regime).
\item		$F^*$ is the infimum of all $F\geq 0$ such that $\liminf_{t\rightarrow\infty} \frac{h(y,t)}{t}>0$ for all $y\in\R^{d-1}$ and for all $F_{uni}>F$  (ballistic regime).
\end{itemize}
One might assume that there is a dichotomy between these two behaviors ($F_*=F^*$), but in general there might be a sub-ballistic regime between them, in which no pinning occurs but in which the effective speed vanishes, $\limsup_{t\rightarrow\infty} \frac{h(y,t)}{t}=0$.
Note that for \eqref{intro_eq_normal} without a curvature component and for compactly supported, separated obstacle shapes we would have $F_*=F^*=0$ since the interface would propagates with velocity $F_{uni}$ along the obstacle-free areas. 
However, with a curvature component -- corresponding to surface tension --, the interface might remain pinned in this setting even for $F_{uni}>0$ with the obstacles acting as `anchors', see the results cited below.

\subsubsection{Pinning and depinning for periodic media}
For the QEW model in arbitrary dimension $d\in\N$ and a periodic obstacle field $F_{obst}$, it was shown in \cite{DiYi06} that $0<F_*<+\infty$ and, under a non-degeneracy assumption on the obstacle field, that $F^*=F_*$ with the height advancing at an effective speed $\vhom\sim\pp{F_{uni}-F^*}^{\frac{1}{2}}$ for $0<F_{uni}-F^*\ll 1$. 
For the QEW model in dimension $d=2$ and periodically spaced, soft, and compactly supported obstacles, a scaling relation between the strength of the obstacles and the critical forcings $F_*$, $F^*$ was established in \cite{CBD20}.

\subsubsection{Pinning and depinning for random media}
For the QEW model in dimension $d=2$ and periodically spaced, compactly supported obstacles with random strength independently sampled from an exponential distribution, it was shown in \cite{CDL10} that $F_*<+\infty$.
In the same setting, $F^*<+\infty$ was established in \cite{DoSch12} and extended to arbitrary dimension $d\in\N$ for a spatial semi-discretization in \cite{DoSch17}.
For a full time-space-discretization of this model, $F_*=F^*$ was obtained in \cite{BoTe14}.
However, in \cite{DoSch17} it was shown that $F_*=F^*$ does not hold in general, using parallel planar obstacles with random sparse holes as counterexample with a sub-ballistic regime.

For general discrete interface motions in the graph setting, in \cite{BoTe14} it was shown that the evolution is ballistic assuming that the probability for the evolution getting blocked in a box decays polynomially as the box size increases, using a general mixing condition rather than the strict independence required for most of the other results in the literature. 

For compactly supported obstacles distributed according to a Poisson point process in arbitrary dimension $d\in\N$, $0<F_*<+\infty$ was established in \cite{DDS11} for both the QEW model and forced mean curvature flow using a percolation result from \cite{DDGHS10}.
For the QEW model in $d=2$, the same has been shown while allowing for obstacles with different signs in \cite{DJS21}.
If in addition the strength of the individual obstacles is sampled from a probability distribution with a sufficiently thick tail, then for the QEW model it holds that $F_*=+\infty$, that is, pinning occurs no matter how strong the uniform forcing is, see \cite{DJS22}.
In contrast, an interface evolving by mean curvature flow with a sufficiently strong uniform forcing would just ``wrap around'' the obstacles, no matter how strong they are.
In fact, given impenetrable obstacles distributed according to a Poisson point process with intensity $\rho$ in dimension $d=2$, it was shown in \cite{CDO22} that $c\sqrt{\rho}\leq F_* \leq F^* \leq C\sqrt{\rho}$ holds with suitable constants $c,C>0$ for forced mean curvature flow and related geometric motions.
An alternative proof for $F^* \leq C\sqrt{\rho}$ will be given in Section~\ref{s_poisson} by applying our main result to this setting.

\subsection{Percolation theory}\label{subs_lit_percolation}
Without its curvature component ($a\equiv 0$) and forcing and obstacles independent of the direction ($F(x,e)=F(x)$), the evolution \eqref{intro_eq_normal} would correspond to Riemannian first passage percolation as introduced in \cite{LaGW10} with the speed at $x\in\Rt$ given by $F(x)$.
The domain is hence split into an `open' area $\Bp{x\,:\,F(x)>0}$ and a `closed' area $\Bp{x\,:\,F(x)\leq 0}$.
Rather than having to consider the global evolution of a set in order to determine if or when a point is reached, it is sufficient to check for paths along open areas connecting the point of the initial set and compute the passage time.
In classical (first passage) percolation theory \cite{Gri99, ADH17}, typically a discrete grid is considered instead of continuum space.
In this setting, our box criterion breaks down into two questions. 
First, in a given box, is there path in the open area connecting opposite sides of the box?
Second, is there such a path along which the travel time is faster than $v_0^{-1}r$ for some fixed speed $v_0>0$?

In percolation theory \cite{Gri99}, i.e., the study of the connectivity properties, the first question plays a prominent role. 
The probability that the answer is yes is called the crossing probability. 
For example, RSW theory allows the transfer of probability bounds for crossing probabilities between rectangles with different aspect ratios, e.g., lower bounds for crossing the short side of a box imply lower bounds for crossing parallel to the long side of a thin box. 
This theory originated with the seminal works of Russo \cite{Rus78} and Seymour and Welsh \cite{SeWe78} on Bernoulli percolation.

Regarding the behavior on large scales, in the supercritical phase (where an infinite open cluster exists), the crossing probability converges to one exponentially fast as the box size increases, see \cite{DeuPi96,DCT20}.
This is similar to our Corollary~\ref{vmin_cor_assBoxScaled}. 
We however only obtain stretched exponential bounds, but note that the aspect ratio of our boxes gets worse as the size increases. 
Our argument in Section~\ref{s_vmin} corresponds to a renormalization or coarse graining scheme which in spirit is similar to such approaches in percolation theory, e.g., \cite{AnPi96, DCT20}.

In first passage percolation on infinite percolation clusters \cite{GaMa04, GaMa07, GaMa09, Mou11, CeThe16, Dem24},  i.e. the study of passage times or the so-called `chemical distance', a key tool has been an almost guaranteed upper bound for passage times on large scales established in \cite{AnPi96}.
For the whole supercritical regime, it is shown that there is a minimum speed such that, assuming two vertices $x,y$ are part of the same cluster, the probability of the passage time exceeding $\vmineff|x-y|$ decays exponentially with rate proportional to $|x-y|$, which is sharp.
Together with control of the size of holes in the infinite cluster for dimension $d=2$, this yields an `almost guaranteed effective minimum speed' as in \ref{veff_aP_star} for first passage percolation with exponential rates. 

Our work can be viewed as a generalization or extension to arbitrary set evolutions in $\Rt$ satisfying a comparison principle.
In particular, applying our main theorem to this setting yields a comparable statement.
Our box criterion is sharp in the sense that it is satisfied for the whole supercritical regime with suitably chosen parameters.
Instead of the exact passage time conditioned on the target point being reachable, we consider the time at which the set evolution reaches the point directly or fully encloses it.
Starting from some point $x\in\Rt$ in the infinite component, our result implies that the probability of this time exceeding $\vmineff|x-y|$ for a target point $y\in\Rt$ decays exponentially with rate proportional to $|x-y|^\ratespeed$ for any $\ratespeed\in(0,1)$ but with $\vmineff=\vmineff(\ratespeed)$.
For higher dimensions, note that an `almost guaranteed effective minimum speed' can not hold for the whole supercritical regime since the complement of the infinite cluster might contain an infinite component too, preventing the existence of a continuously advancing main front up to left-behind enclosures.

\section{Setting and Assumptions}\label{s_assump}
For convenience, we denote dependencies on the parameters introduced with the assumptions below via
\begin{equation*}
	\data\coloneqq 	\pp{C_{1A}, C_{1F}, q, \cw, v_0, p_0, r_0}.
\end{equation*}

\subsection{Admissible coefficient fields}\label{subs_sett_AF}
We say that $(A,F)$ is an admissible coefficient field if the following conditions are satisfied.
Here, $F$ takes the role of the forcing in \eqref{intro_eq_normal} while $A$ corresponds to the (degenerate) diffusion matrix for the level-set formulation -- see Section~\ref{subs_sett_reach} below -- derived from a curvature coefficient field $a=\sigma^2$ as in \eqref{intro_eq_normal} by combining it with a projection onto the tangent space.
We use this diffusion matrix rather than the curvature coefficient field in order to keep the notation aligned with \cite{FI25}.

We assume that $A$ is uniformly Lipschitz continuous and bounded, given by 
\begin{multline}\label{sett_eq_assA}
	A(x,e)=\sigma^2(x,e) \pp{\Idt-e\otimes e}\\
	\,\,\text{ with }
	\norm{\sigma}_{C^{0,1}(\Rt\times\Sone)},\norm{A}_{C^{0,1}(\Rt\times\Sone;\,\Rttsym)} \leq C_{1A}
\end{multline}
for some $\sigma \in C^{0,1}(\Rt\times\Sone)$ and a uniform constant $C_{1A}$.
We further assume that $F$ is also uniformly bounded and Lipschitz, that is there exists $C_{1F}>0$ such that
\begin{equation}\label{sett_eq_assF}
	F\in C^{0,1}(\Rt\times\Sone;\,\R) 
	\text{ with }
	\norm{F}_{C^{0,1}\pp{\Rt\times\Sone;\,\R}} \leq C_{1F}.
\end{equation}
To simplify the notation we use the $0$-homogeneous extension of $A$ and $F$ to $\pp{\Rt\times\Rt\setminus\Bp{0}}$, for $A$ given by $A(x,\xi)\coloneqq A\pp{x,\abs{\xi}^{-1}\xi}$.
For these admissible coefficient fields, we will now introduce the level-set formulation of \eqref{intro_eq_normal}, which is well-posed in the viscosity solution framework.

\subsection{Set evolutions and the level-set formulation}\label{subs_sett_reach}
In general, interface evolutions such as \eqref{intro_eq_normal} can lead to topological changes of the evolving interfaces, allowing for different interpretations of \eqref{intro_eq_normal}. 
In order to give rigorous meaning to the evolution, a variety of frameworks have been introduced, which for some cases are equivalent but differ in others. 
As in \cite{FI25}, we work with viscosity solutions of the level-set formulation.
For a collection of results from the theory of viscosity solutions relevant for our setting, which have been gathered from \cite{CIL92, Giga06, GGIS91,Juu01}, see \cite[Appendix A]{FI25}.

Rather than as an interface evolution, we view 
\begin{equation}
	V= a\pp{x, \nu}H + F\pp{x, \nu}\tag{\ref{intro_eq_normal}}
\end{equation}
as a set evolution with $V$ denoting the outwards normal speed of the boundary, $\nu$ its normal direction, and $H$ the curvature of the boundary. 
In our chosen framework, these set evolutions are given as the sublevel sets of a function $u\colon\Rt\times[0,\infty)\rightarrow\R$ satisfying the level-set equation, a viscous Hamilton--Jacobi equation derived from \eqref{intro_eq_normal} with degenerate diffusion matrix $A(x,e)= a(x,e)\pp{\Idt-e\otimes e}$ for $x\in\Rt$, $e\in\Sone$. 
While it is also possible to view the level sets $t\mapsto\Bp{x \,:\, u(x,t)=0}$ as an interface evolution, `fattening' might occur \cite{BSS93}, meaning that the level set develops non-empty interior.

\begin{definition}[The set evolution in the heterogeneous environment]\label{not_def_reach}
Let $(A,F)$ be an admissible coefficient field satisfying \eqref{sett_eq_assA} and \eqref{sett_eq_assF} respectively. 
Given a closed set $S\subset\Rt$, we define the evolution $\pp{\reach{t}{S}}_{t\geq 0}$ of $S$ with respect to $(A,F)$ as
\begin{equation*}
	\reach{t}{S}=\reach[(A,F)]{t}{S}\coloneqq \Bp{x\in\Rt \,:\, \uls{S}(x,t)\leq 0},
\end{equation*}
where $\uls{S}$ with continuous initial data $\uls{S}(\cdot,0)$ satisfying $S=\Bp{x\in\Rt \,:\, \uls{S}(x,0)\leq 0}$ is the (viscosity) solution of 
\begin{align}
	&\pat u + G_{(A,F)}(x,\nabla u, D^2 u)=0 
	\quad\,\text{in }\Rt\times(0,\infty)\label{intro_eq_uls}\\
	&	\qquad\qquad\quad \text{with }\quad
	G_{(A,F)}(x,\xi,X) \coloneqq - \tr\pp{A\pp{x, \frac{\xi}{\abs{\xi}}}X} + F\pp{x, \frac{\xi}{\abs{\xi}}}\abs{\xi}\label{intro_eq_GAF}
\end{align}
for $(x,\xi,X)\in \Rt\times\pp{\Rt\!\setminus\!\Bp{0}}\times \Rttsym$ and the singularity at $\xi=0$ treated via semi-continuous envelopes in the framework of viscosity solutions.
Such a solution always exists and $\reach{t}{S}$ is unique and independent of the specific initial data $\uls{S}(\cdot,0)$.
We can assume that $\uls{S}$ is uniformly continuous on $\Rt\times(0,T)$ for any $T>0$, see \cite[Theorem~A.14, Remark~A.15, Theorem~A.17]{FI25}.
\end{definition}

%

\subsection{Assumptions on the probability distribution}\label{subs_aP}
In short, we assume that the coefficient field $(A,F)$ in \eqref{intro_eq_uls} is a random variable sampled from a stationary probability measure with finite range of dependence, which satisfies the box criterion as informally described in Assumption~\ref{sett_aP_pseudoBox}.

More explicitly, the probability space is defined as the set of admissible coefficient fields
\begin{equation}\label{sett_eq_omega}
	\Omega\coloneqq \Bp{(A,F)\,:\,A\text{ satisfies }\eqref{sett_eq_assA},\,F\text{ satisfies }\eqref{sett_eq_assF}}.
\end{equation}
On $\Omega$ we define the $\sigma$-algebras $\Bp{\F(U)}$ for all Borel sets $U\subset\Rt$,
\begin{multline}\label{sett_eq_sigmaF}
  	\F(U)\coloneqq 	\text{the $\sigma$-algebra generated by the family of maps}\\
  					\Omega\rightarrow \Rtt\times\R,\,
  					(A,F)\mapsto \pp{A(x,e), F(x,e)}\quad
  					\text{with }x\in U,\,e\in\Sone.
\end{multline}  
The $\sigma$-algebras $\F(U)$ can be seen as containing all of the information derived from observing the coefficient field in $U\subset\Rt$.

We assume that the randomness of the coefficient field is described by a probability measure $\Pm$ on $(\Omega,\F(\Rt))$, which satisfies the following conditions.
The first two assumptions are standard in the theory of stochastic homogenization:
\begin{enumerate}[label=(P\arabic*), left=8pt]
	\item	\label{sett_aP_stat} \textbf{Stationarity.}
			 $\Pm$ is invariant under translations of the coefficient field: 
			\begin{equation*}
				\Pm=\Pm\circ \tau_x, \quad\text{for any }x\in\Rt,
			\end{equation*}
			 where 
				$\tau_x(\wt{\Omega})=\lbrace\pp{A(\cdot+x,\cdot),F(\cdot+x,\cdot)}\,:\,(A,F)\in\wt{\Omega}\rbrace$  
			for $\wt{\Omega}\in\F(\Rt)$.
	\item	\label{sett_aP_fin} \textbf{Finite range of dependence.}
			For all Borel sets $U,V\subset\Rt$ with $\dist(U,V)\geq 1$,
			\begin{equation*}
				\F(U)\text{ and }\F(V)\text{ are $\Pm$-independent.}
			\end{equation*}
\end{enumerate}
The final assumption is our box criterion, which is illustrated in Figure~\ref{intro_fig_box} and the informal version of which was stated in the introduction as \ref{sett_aP_pseudoBox}.
For sampling the coefficient field we use boxes scaling sub-linearly in width, that is, given a length scale $r>0$ we consider the box
\begin{align}\label{box_eq_Qup}
	\Qup{r}	&\coloneqq [0,w(r)]\times[-w(r),r]
	&\text{with }\quad	
	w(r)		&\coloneqq \pp{\frac{r}{\cw}}^{\frac{1}{q}}.
\end{align}
for some fixed $q>1$ and $\cw>0$.
For such a box $x+\Qup{r}$ of size $r>0$ anchored at $x\in\Rt$, for a prescribed average speed $v>0$ we consider the event $\Eup[v]{r,x}$ that there is a non-shrinking set in the bottom of the box which within the box evolves to reach the top of the box with an average speed of at least $v$.
Formally, we define this event as
\begin{align}
	\label{box_E}\Eup[v]{r,x}=
	\Big\lbrace
		(A,F)\in\Omega \,:\,\,		&	\text{there is a closed set }S\subset x+[0,w(r)]\times[-w(r),0]\\
	\label{box_E_stable}			&	\text{such that } S \text{ is stable with respect to } (A,F)_{x+\Qup{r}}\\
	\label{box_E_top}			&	\text{and } \reach[(A,F)_{x+\Qup{r}}]{v^{-1}r}{S} \cap \pp{x+[0,w(r)]\times\Bp{r}} \neq \emptyset
	\Big\rbrace.
\end{align}
where \emph{`stable'} sets and the \emph{`restricted coefficient field'} $(A,F)_{x+\Qup{r}}$ are defined below in Section~\ref{subs_box}.
The evolution of stable sets is non-shrinking, ensuring that the evolved set stays connected to the \emph{`starting area'} $x+[0,w(r)]\times[-w(r),0]$, see Lemma~\ref{sti_lem_charEup}.
The restriction of the coefficient field confines the evolution to a small neighborhood of the box and ensures the measurability of $\Eup[v]{r,x}$ with respect to $\F(x+\Qup{r}+B_1)$, see Lemma~\ref{meas_lem_Eup}. 

\begin{figure}[h]
\centering
\includegraphics[width=0.95\textwidth]{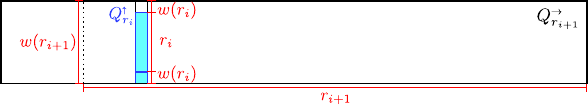}
\caption[An illustration of how $r_0$ and $r_1$ relate]{	An illustration of how $r_0$ and $r_1$ relate.
		The length scales are chosen such that translations of $\Qup{r_0}$ and $\Qd{r_0}$ fit into $\Qr{r_1}$ without their starting areas overlapping.
		}
\label{sett_fig_r1}
\end{figure}

We write 
\begin{equation*}
	\Eup[v]{r}\coloneqq \Eup[v]{r,0} = \tau_{x}\pp{\Eup[v]{r,x}},
\end{equation*}
where $\tau_{x}$ is a translation as defined in Assumption~\ref{sett_aP_stat}.
Due to the stationarity from Assumption~\ref{sett_aP_stat}, we have $\PM{\Eup[v]{r}}=\PM{\Eup[v]{r,x}}$ for all $x\in\Rt$.
We require that there is an initial length scale $[r_0,r_1)$ and a speed $v_0>0$ such that this event and its counterparts $\Er[v_0]{r}, \Ed[v_0]{r}, \El[v_0]{r}$ for the other directions are sufficiently likely, with the threshold depending on the length scale and the width-scaling:

\begin{enumerate}[label={(P$\,\vrectangle$)}, left=0pt]
\item	\label{sett_aP_box} \textbf{Box criterion for $d=2$.}
		There exist a velocity $v_0>0$, a width scaling exponent $q>1$ and a width scaling factor $\cw>0$ such that
		\begin{equation}\label{sett_eq_assBox}
			\min\Bp{\PM{\Eup[v_0]{r}},\PM{\Er[v_0]{r}},\PM{\Ed[v_0]{r}},\PM{\El[v_0]{r}}}\geq 1-p_0,
			\quad\text{for all }	r\in[r_0,r_1),
		\end{equation}	
		with $r_1\coloneqq w^{-1}(r_0+2w(r_0))$, see Figure~\ref{sett_fig_r1}, for some initial box length $r_0>0$ and an acceptable probability of failure $p_0\geq 0$ satisfying
		\begin{align}
				r_0 &\geq \max\Bp{3^q\cw,~ \pp{24q}^{\frac{q}{q-1}}\cw^{-\frac{1}{q-1}}},\label{sett_eq_assBox_r}\\					
				p_0	&< \min\Bp{	\exp\pp{-\frac{\log(6)+q\log(k_0)}{k_0-1}}, 
								\exp\pp{-\frac{2q}{k_0}}},\notag\\
				&\hspace{160pt}\text{where }\quad	k_0	\coloneqq	\left\lfloor \frac{1}{12q}\cw^{\frac{1}{q}}r_0^{\frac{q-1}{q}} \right\rfloor\geq 2.\label{sett_eq_assBox_p}
		\end{align}
\end{enumerate}
Here, the events $\Er[v]{r,x},\Ed[v]{r,x},\El[v]{r,x}$ above are defined analogously to $\Eup[v]{r,x}$, replacing the vertical box $\Qup{r}$ with 
\begin{align*}
	\Qr{r}		&\coloneqq [-w(r),r]\times[-w(r),0]	& \text{with starting area }&[-w(r),0]^2,\\
	\Qd{r}		&\coloneqq [-w(r),0]\times[-r,w(r)]	& \text{with starting area }&[-w(r),0]\times[0,w(r)],\\
	\Ql{r}		&\coloneqq [-r,w(r)]\times[0,w(r)]	& \text{with starting area }&[0,w(r)]^2.
\end{align*}

Note that the threshold for $p_0$ in \eqref{sett_eq_assBox_p} can be arbitrarily close to one if the initial scale is sufficiently large (and the corresponding boxes sufficiently thin).

\begin{remark}\label{sett_rem_single-box-check}
Note that for \ref{sett_aP_box} it can be sufficient to check a single box size rather than every size in a given length scale, albeit with harsher requirements on the probability bound.
Probability bounds for a continuous range of sizes and thus Assumption~\ref{sett_aP_box} can then be obtained using a construction similar to the one in Section~\ref{s_vmin}.
\end{remark}

\subsection{Stable sets and restricted coefficient fields}\label{subs_box}
We now make precise the notions of `stable' sets and `restricted coefficient fields' used to define the event $\Eup[v]{r,x}$ in \eqref{box_E} for the box criterion \ref{sett_aP_box}.

\subsubsection{Introduction of Stable Sets}
Throughout this paper, sets which do not shrink under the interface evolution will be central to many arguments.
We call these sets or the corresponding interfaces \emph{stable}.
We will discuss properties of these sets in Section~\ref{subs_stable}.

\begin{definition}[Stable Sets]\label{box_def_stable}
	A closed set $S\subset\Rt$ is called stable (with respect to $(A,F)\in\Omega$) if $u\coloneqq 1-\charfun[S]$ is a supersolution of 
	\begin{equation}\label{box_eq_stable}
		- \tr\pp{A\pp{x, \frac{\nabla u}{\abs{\nabla u}}}D^2 u} 
		+ F\pp{x, \frac{\nabla u}{\abs{\nabla u}}}\abs{\nabla u}
		= 0
		\quad\text{in }\Rt.
	\end{equation}
\end{definition}

\subsubsection{Introduction of restricted coefficient fields}
For the box criterion, we want to modify the forcing to ensure that the initial set evolves and stays within a small neighborhood of the given box.
This should also allow us to apply the finite range of dependence from Assumption~\ref{sett_aP_fin} for boxes which are sufficiently far apart.
For this purpose, we introduce suitably restricted coefficient fields.

\begin{definition}[Restricted coefficient fields]\label{restr_def}
Let $Q\subset\Rt$ be a Borel Set. 
Given $(A,F)\in\Omega$ we define $(A,F)_Q\coloneqq \pp{A_Q,F_Q}$, the coefficient fields restricted to $Q$, as 
\begin{equation}\label{box_eq_defAQ}
	A_Q(x,e)\coloneqq
	\left\lbrace\begin{aligned}
		&A(x,e)							&&\text{for }\dist(x,Q)\leq \frac{2}{3},\\
		&9\pp{1-\dist(x,Q)}^2 A(x,e)		&&\text{for }\dist(x,Q)\in\bp{\frac{2}{3},1},\\
		&0								&&\text{for }\dist(x,Q)\geq 1,\\
	\end{aligned}\right.
\end{equation}
\begin{equation}\label{box_eq_defFQ}
	F_Q(x,e)\coloneqq
	\left\lbrace\begin{aligned}	
		&F(x,e)			\qquad\qquad\qquad	\qquad\qquad\qquad\,		\text{for }\dist(x,Q)\leq 0,\\					
		&3\pp{\frac{1}{3}-\dist(x,Q)} F(x,e) - 3\dist(x,Q)\max\Bp{C_{1F},\,6C_{1A}}	\\			
		&\qquad\quad		\qquad\qquad\qquad	\qquad\qquad\qquad\,		\text{for }\dist(x,Q)\in\bp{0,\frac{1}{3}},\\				
		&-\max\Bp{C_{1F},\,6C_{1A}}	\,		\qquad\qquad\qquad\,		\text{for }\dist(x,Q)\geq \frac{1}{3}.
	\end{aligned}\right.
\end{equation}
\end{definition}

By definition, the evolution of a set with respect to $(A,F)_Q$ only depends on $(A,F)$ restricted to $Q+B_1$ and hence should be $\F(Q+B_1)$ measurable, see e.g. Lemma~\ref{meas_lem_Eup}.
As a consequence of this definition, any $S\subset Q$ stays stuck in $Q+B_{\frac{2}{3}}$ when evolving by \eqref{intro_eq_uls} with coefficients $(A,F)_Q$.
In addition, in this case the evolution of a set with respect to $(A,F)_Q$ provides a lower bound to the evolution of the same set with respect to the unmodified coefficient field $(A,F)$.
For the rigorous statement see Lemma~\ref{restr_lem_PropRestr} in Section~\ref{subs_restr}.

\begin{remark}[Well-posedness of evolution for restricted coefficient fields]
Existence, uniqueness and uniform continuity for solutions of the level-set equation as well as invariance for the evolution of the level sets as stated in Definition~\ref{not_def_reach} for $(A,F)\in\Omega$ also hold if $(A,F)=(\wt{A},\wt{F})_Q$ for some $(\wt{A},\wt{F})\in\Omega$ and $Q\subset\Rt$.
\end{remark}

\begin{remark}[Choice of \eqref{box_eq_defAQ} and \eqref{box_eq_defFQ}]
Instead of transitioning from $(A,F)(\cdot)$ on $Q$ to $(0,\ov{F})$ with $\ov{F}=-\max\Bp{C_{1F},\,6C_{1A}}$ via interpolation on a scale of $1$, we could choose an arbitrarily small scale $s>0$ -- as long as we adapt $\ov{F}$ accordingly (see \cite[Proof of Lemma~5.10]{FI25}, adaptation only necessary for non-convex $Q$).
Then the evolution of a set with respect to $(A,F)_Q$ would be $\F(Q+B_s)$ measurable.
This and being more precise in the proof of Proposition~\ref{vmin_prop_assBoxScaled} would for example allow us to weaken the assumption on $r_0$ in \eqref{sett_eq_assBox_r} and \eqref{sett_eq_assBox_p} from \ref{sett_aP_box}.
The arguments in Section~\ref{s_vmin} stay the same, except that we place boxes with a distance of $1+2s$ instead of $3$ for independence.
Therefore, we just choose a scale of $1$ for a more stringent notation.
\end{remark}

\subsection{The effective minimum speed assumption from \cite{FI25}}\label{subs_sett_veff}
We now precisely state the effective minimum speed assumption \ref{veff_aP_star} including the notions of an `effective minimum speed' of propagation, `approximate stability' and `fat' sets as introduced in \cite{FI25}.
For detailed descriptions and illustrations of these notions see \cite[Section~3.3]{FI25}.

Intuitively, a set evolves with an `effective minimum speed' on a scale $h>0$ if -- at any step of the evolution -- the set is guaranteed to expand by at least $h$ in every direction within a time proportional to $h$, up to ignoring enclosures of size up to $h$ left behind the main front. 
To formalize ignoring such enclosures, we introduce the notion of $h$-envelopes.
A set is `approximately stable' on a scale $h>0$ if there is a good interior approximation which is stable, i.e., non-shrinking, and $h$-envelops the original set after a time proportional to $h$.
A set is $h$-`fat' if it can be written as the union of balls with radius $h$.

\begin{definition}\label{veff_def_fatstab}
Let $(A,F)\in\Omega$.
For a set $S\subset\Rt$, let $t\mapsto \reach[(A,F)]{t}{S}$ denote the set evolution from Definition~\ref{not_def_reach}.
Let $h>0$, $c>0$.
\begin{enumerate}[label=(\roman*)]
\item	We say that a set $S_1\subset\Rt$ \emph{$h$-envelops} another set $S_2\subset\Rt$, denoted by $S_2\hsubs S_1$,
		if for each connected component $H$ of $\comp{(S_1)}$ either $H\cap S_2 = \emptyset$ or $H\subset S_1+\ov{B}_h$.	
\item	We say that $(A,F)$ admits $v>0$ as an \emph{effective minimum speed of propagation} on the scale $h$ in $M\subset\Rt$ for a set $S\subset\Rt$ if for any $t\geq 0$ it holds that
		\begin{equation}\label{veff_eq_effspeed}
			\quad	\pp{\reach[(A,F)]{t}{S}\cap M} + \ov{B}_h  \hsubs \reach[(A,F)]{t+v^{-1}h}{S}.
		\end{equation} 
\item	We say that a set $S\subset\Rt$ is \emph{$(h,c)$-approximately stable} with respect to $(A,F)$
		if there exists a set $S_h\subset S$ such that
		\begin{align*}
				\quad	S_h \text{ is stable},
			&&		S	\hsubs \reach[(A,F)]{c h}{S_h},
			&&		S_h	\supset \Rt\setminus\pp{\comp{S}+c\vmax B_h},
		\end{align*}		
		where $\vmax=\vmax(\data)$ corresponds to the maximum speed with which the interface propagates on larger scales, see Lemma~\ref{not_lem_vmax}.
		In this case, $S_h$ is called a \emph{stable $(h,c)$-approximation} of $S$.
\item	We say that a set $S\subset \Rt$ is \emph{$h$-fat} if 
		\begin{equation*}
			\quad	S=\bigcup\Bp{\ov{B}_h(x)	\,:\,	x\in\Rt \text{ with }	\ov{B}_h(x)\subset S}.
		\end{equation*}
\end{enumerate}
\end{definition}

Using these notions, we now state the rigorous version of the effective minimum speed assumption \ref{veff_aP_star} from \cite{FI25}, with the informal version \ref{sett_aP_pseudoStar} stated in the introduction.
In short, for coefficient fields satisfying \ref{veff_aP_star}, there is a velocity $\vmineff>0$ such that -- as the scale increases -- it becomes increasingly likely that all sufficiently fat sets in an area have stable approximations for which $\vmineff$ is admitted as an effective minimum speed of propagation.
In order to ensure measurability with respect to a small neighborhood of the area under consideration, we require a family of events guaranteeing the effective minimum speed rather than defining these events directly.

In the present work, our goal is to establish \ref{veff_aP_star} for suitable parameters based on Assumptions \ref{sett_aP_stat}, \ref{sett_aP_fin}, and \ref{sett_aP_box} from Section~\ref{subs_aP}.

\begin{enumerate}[label=(P$\circledcirc$), left=0pt]
\item	\label{veff_aP_star}	\textbf{Almost guaranteed effective minimum speed}.
		There is $\vmineff>0$, $h_0\geq 1$ and a rate $\ratespeed>0$ with constants $\cspeed, \Cspeed, \Cmeas, \cstable>0$, such that for each $h\geq h_0$ there is a family of events $\pp{\Espeed{h}(M)}_{M\subset\Rt}$ with 
		\begin{equation}\label{veff_eq_probBound}
			\quad	\Pm[\Espeed{h}(M)]\geq 1 - \Cspeed\pp{1+\frac{\diam(M)}{h}}^2\exp(-\cspeed h^{\ratespeed}),
		\end{equation}			
		which satisfies
		\begin{enumerate}[label=($*$\roman*), left=5pt]
			\item	\label{veff_Pmeas}
					$\Espeed{h}(M)$ is measurable with respect to $\F(M+\Cmeas \ov{B}_h)$,
			\item	\label{veff_Porder}
					for all $M_1\subset M_2$ it holds that $\Espeed{h}(M_1)	\supset	\Espeed{h}(M_2)$,
			\item	\label{veff_Pspeed}
					for all $(A,F)\in \Espeed{h}(M)$, every $h$-fat set $S$ with $\partial S\subset M$ has a stable $(h,\cstable)$-approximation, 
					for which $(A,F)$ admits $\vmineff>0$ as an effective minimum speed of propagation on the scale $h$ in $M\subset\Rt$, see Definition~\ref{veff_def_fatstab} above for the precise introduction of these notions.
		\end{enumerate}	
\end{enumerate} 

In particular, \ref{veff_aP_star} implies ballistic motion on large scales. 
However, ballistic motion does not necessarily imply the existence of an (almost guaranteed) effective minimum speed as introduced here, at least not in general dimension. 
For ballistic motion it is sufficient to compare the main front to the regular initial set.
But for an effective minimum speed, the main front has to expand compared to the state of the set evolution at any time, which is potentially much less regular than the initial set. 

For an illustration, consider the counterexample from \cite{CaMo14} to periodic homogenization of forced mean curvature flow despite a uniformly positive combined forcing in dimension $d=3$:
Thin, infinitely long, parallel cylinders are placed periodically in space with the forcing inside the cylinders much stronger than outside with a transition layer in-between. 
Due to the mean curvature, sets initialized inside the cylinders are not able to pass the transition layer.
The cylinders can hence be thought of as pipes along which sets expand much faster but to which these sets remain constrained.
Choosing the half space orthogonal to the pipes as the initial set, the main front will expand in the direction parallel to the pipes at different speeds: 
With a fast head speed inside the pipes and with a much slower but still positive tail speed in the remaining space, lagging increasingly far behind the evolution inside the pipes.
We hence observe a ballistic motion, but the evolution does not admit an effective minimum speed.
For the points reached early within a pipe, the set will not expand into the close-by area between the pipes: 
This area remains connected to the complement of the set evolution until it will be covered much later by the part of the front which advances with the slower tail speed.

\section{Main Result and the Strategy for the Proof}\label{s_str}

\subsection{The main result}\label{subs_str_main}
Our main result establishes an effective minimum speed and thus ballistic behavior for sets evolving by \eqref{intro_eq_normal} for random coefficient fields sampled from a stationary coefficient field with a finite range of dependence and which satisfies our box criterion.

\begin{theorem}[Main result]\label{main_thm}
Let $(A,F)\in\Omega$ be a random coefficient field sampled from a probability measure $\Pm$ on $(\Omega,\F(\Rt))$ which
\begin{itemize}
\item	is stationary, that is invariant with respect to translations (Assumption~\ref{sett_aP_stat}),
\item	has finite range of dependence (Assumption~\ref{sett_aP_fin}), and
\item	satisfies our box criterion (Assumption~\ref{sett_aP_box} with width-scaling $q>1$).
\end{itemize}
Then the coefficient field admits an almost guaranteed effective minimum speed:
the coefficient field satisfies \ref{veff_aP_star} for every $\ratespeed\in \pp{0, \pp{1-\frac{1}{q}}\pp{1+\frac{1}{q^2}}}$ with coefficients $\vmineff=\vmineff(\ratespeed,\data)>0$, $h_0=r_0$, $\cspeed=\cspeed(\ratespeed,\data)>0$, $\Cspeed>0$, $\Cmeas>0$, and $\cstable=\cstable(\data)>0$. 
\end{theorem}

As a consequence, we obtain Theorem~\ref{intro_thm_Poisson} by showing that the assumptions of Theorem~\ref{main_thm} are satisfied for obstacles distributed according to a Poisson point process and a sufficiently large uniform forcing, see Section~\ref{s_poisson}. 
Theorem~\ref{main_thm} then yields an almost guaranteed effective minimum speed, which allows the application of the homogenization results from \cite{FI25} to Poisson distributed obstacles and yields the announced \cite[Theorem~1.1]{FI25}.

\subsection{Strategy of the proof}\label{subs_str_str}
\begin{figure}[h]
\centering
\includegraphics[width=0.4\textwidth]{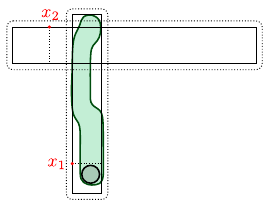}
\hspace{10pt}
\includegraphics[width=0.4\textwidth]{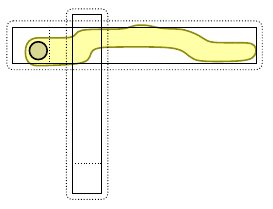}\\\vspace{15pt}
\includegraphics[width=0.4\textwidth]{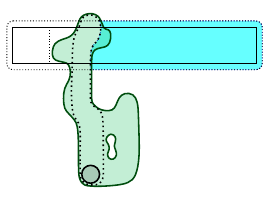}
\hspace{10pt}
\includegraphics[width=0.4\textwidth]{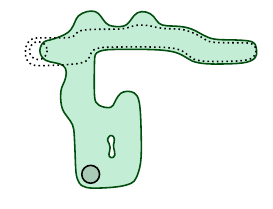}
\caption[An illustration of our `stitching lemma' applied to chain the propagation along two consecutive boxes]{	An illustration of our `stitching lemma' applied to chain the propagation along two consecutive boxes.
			Assume that for two intersecting boxes the respective events $\Eup[v]{r_1,x_1}$ and $\Er[v]{r_2,x_2}$ hold for the coefficient field $(A,F)$.
			Top left: Due to $\Eup[v]{r_1,x_1}$, there is a stable set $S_1$ (depicted in grey) which reaches the top of the first box after evolving for $t_1=v^{-1}r_1$ with respect to the coefficient field restricted to the first box (depicted in green).
			Top right: Due to $\Er[v]{r_2,x_2}$, there is a stable set $S_2$ (depicted in grey) which reaches the right end of the second box after evolving for $t_2=v^{-1}r_2$ with respect to the coefficient field restricted to the second box (depicted in yellow).
			Bottom left: By comparison, the evolution of $S_1$ with respect to the unrestricted coefficient field at $t_1$ (depicted in green) separates $S_2$ from the right side of the second box (depicted in blue).
			Bottom right: Since $S_1$ and thus its evolution are stable, the stitching lemma allows us to compare the evolution of $S_1$ and $S_2$ in that area and we obtain that the evolution of $S_1$ reaches the right end of the second box at $t_1+t_2$ (depicted in green)}
\label{str_fig_box-stitching}
\end{figure}

In Section~\ref{s_pre}, we establish some preliminary results about the properties of stable sets and set evolutions with respect to restricted coefficient fields.
The proof of Theorem~\ref{main_thm} then consists of two major steps:
\begin{enumerate}
	\item	The box criterion guarantees lower bounds for the probability that the interface propagates along boxes of length $r\in[r_0,r_1)$ with average speed $v_0$.
			The key step of our argument in Section~\ref{s_vmin} is to iteratively extend these probability bounds to larger box lengths $r\in[r_i,r_{i+1})$ with average speed $v_i$ for $i\in\N$, while at the same time improving the probability bounds with each step and sufficiently controlling the average speed such that $\liminf_{i\rightarrow\infty} v_i>0$.
	\item	In Section~\ref{s_veff}, we show that a suitable network of such boxes along which the interface propagates guarantees an effective minimum speed in the sense of \ref{veff_aP_star}. 
			Due to the results from Section~\ref{s_vmin}, the existence of such a network is sufficiently likely for the probability bounds in \ref{veff_aP_star} to be satisfied.
\end{enumerate}
A key tool in the present work is Lemma~\ref{restr_lem_stitch}, our `stitching lemma'. 
Combining our notion of stable sets from Definition~\ref{box_def_stable} with the properties of the restricted coefficient fields from Definition~\ref{restr_def} allows us to chain the guaranteed propagation for consecutive boxes crossing each other, see Figure~\ref{str_fig_box-stitching}. 	
The main method for our proofs is to construct such overlapping box configurations and derive the targeted property on the larger scale.
The finite range of dependence from Assumption~\ref{sett_aP_fin} is applied for subgroups of suitably spaced boxes to control or improve the probability.

\subsubsection{Lifting the box criterion to larger scales while improving the probability bound}
The box criterion \ref{sett_aP_box} guarantees that $\PM{\Eup[v_0]{r}}\geq 1-p_0$ for $r\in[r_0,r_1)$ and the same for the other directions.
The first goal is to iteratively lift these to larger scales while improving the probability bound. 
That is, for all $i\in\N$ we want to show that 
\begin{align}\label{str_eq_iteration}
	\min\Bp{\PM{\Eup[v_i]{r}},\PM{\Er[v_i]{r}},\PM{\Ed[v_i]{r}},\PM{\El[v_i]{r}}}&\geq 1-p_i
	&&	\text{for all }r\in[r_i,r_{i+1})
\end{align}
for $v_i>0$ and $p_i\in[0,1)$ while ensuring that $v_i\geq \vmin$ for some $\vmin>0$ and $p_i\rightarrow 0$ for $i\rightarrow\infty$ as fast as possible.
As $r_1$ in \ref{sett_aP_box}, the sequence $(r_i)_i$ is given by $w(r_{i+1}) = r_i+2 w(r_i)$, meaning that the box-width of the next scale corresponds to the box-length (and an additional starting area) of the current scale, see Figure~\ref{sett_fig_r1}.

Given a large box $\Qup{r}$, the idea is to place layers of smaller boxes within and show that if the interface propagates along suitable subsets of the smaller boxes then the interface propagates from the starting area of the large box to its top.
To be more precise, starting from the top of the starting area, we place multiple layers of smaller boxes in $\Qup{r}$, see Figure~\ref{vmin_fig_box-config}. 
Each layer consists of three subgroups of size $\kvert, \khor\in\N$ to be fixed later: 
\begin{enumerate}
\item	$\kvert$ parallel translations of $\Qup{\rup}$ with $\rup=\rup(r,\kvert)\ll r$,
\item	$\khor$ parallel translations of $\Qr{\rbase}$ with $\rbase + 2 w(\rbase)=w(r)$ all cutting through the top of the first subgroup from left to right, potentially enabling the application of the stitching lemma,
\item	$\khor$ parallel translations of $\Ql{\rbase}$ all cutting through the top of the first subgroup from right to left.
\end{enumerate}
Within each subgroup, the boxes are spaced with distance $3$ such that the respective restricted coefficient fields are independent.
The next layer is then placed such that the horizontal boxes of the previous layer cut through the bottom of the new vertical boxes, until the top of  $\Qup{r}$ is reached.
For the number of layers we write $N=N(r,\kvert,\khor, w(\cdot))$.
Note that the boxes are placed in such a way that, if for each subgroup there is at least one box with a stable set in its starting area propagating to the opposite end, then the stable set from the first subgroup is in the starting area of $\Qup{r}$ and -- due to the stitching lemma -- propagates layer by layer until it reaches the top of $\Qup{r}$, see Figure~\ref{vmin_fig_stitch}.

Suppose that 
\begin{align*}
	\min\Bp{\PM{\Eup[v]{\rup}},\PM{\Er[v]{\rbase}},\PM{\El[v]{\rbase}}}&\geq 1-p
\end{align*}
for some $v>0$ and $p\in[0,1)$.
Then due to the stationarity from Assumption~\ref{sett_aP_stat} and due to the finite range of dependence from Assumption~\ref{sett_aP_fin}, with the argument from above we obtain
\begin{align}
	\PM{\Eup[\wt{v}]{r}} &\geq 1 - N\pp{1-\PM{\Eup[v]{\rup}}}^{\kvert} - N\pp{1-\PM{\Er[v]{\rbase}}}^{\khor} - N\pp{1-\PM{\El[v]{\rbase}}}^{\khor}\notag\\
	&\geq 1 - N p^{\kvert} - 2N p^{\khor} \label{str_eq_Pbound}
\end{align}
for some average speed $\wt{v}=\wt{v}(v,r,\kvert,\khor, w(\cdot))$.

While the average speed along the smaller boxes is $v$, this average speed $\wt{v}$ for reaching the top of the large box is slower due to the extra distance covered. 
For each of the $N$ layers, the horizontal distance $\rbase$ is traversed at average speed $v$. 
In addition, between each of the layers there is an overlap of the vertical subgroups corresponding to the height $\khor \pp{w(\rbase)+3}$ of the horizontal subgroups.
Finally, the last layer might have to be pushed down by up to $\rup$ to stay within the larger box $\Qup{r}$.
In summary, we obtain 
\begin{equation}
	\wt{v}\geq v\frac{r}{r+N\rbase+N\khor\pp{w(\rbase)+3}+\rup}. \label{str_eq_vbound}
\end{equation}
Iteratively repeating this argument for all directions yields \eqref{str_eq_iteration} for all $i\in\N$ with \eqref{str_eq_Pbound} governing the behavior of $(p_i)_{i\in\N}$ and with \eqref{str_eq_vbound} governing the behavior  $(v_i)_{i\in\N}$.
Due to the sublinear scaling of the box-width, an increasingly smaller percentage of the minimum average speed is shaved off in \eqref{str_eq_vbound} while at the same time we are still able to choose increasingly large $\khor,\kvert$ as $r\rightarrow\infty$.
This leads to exponentially improving probability bounds provided that relative to each other $r_0$ is large enough and $p_0$ is small enough, corresponding to the bounds in Assumption~\ref{sett_aP_box}.

As a consequence, we do indeed obtain $\vmin>0$ such that $v_i\geq \vmin$ for all $i\in\N$ and a positive rate $\ratebox>0$ such that 
\begin{align}\label{str_eq_pbound_exp}
	\min\Bp{\PM{\Eup[\vmin]{r}},\PM{\Er[\vmin]{r}},\PM{\Ed[\vmin]{r}},\PM{\El[\vmin]{r}}}&\geq 1-\exp\pp{-c r^\ratebox}
\end{align}
for all $r\geq r_0$ and some $c>0$, which is the main result from Section~\ref{s_vmin}, see Corollary~\ref{vmin_cor_assBoxScaled}.

\subsubsection{Effective minimum speed}
In order to show that there is an almost guaranteed effective minimum speed $\vmineff>0$ in the sense of \ref{veff_aP_star}, for a given set $M\subset\Rt$ and a scale $h\gg 1$ we need to define an event $\Espeed{h}(M)$ which implies that all $h$-fat sets with boundary in $M$ have a stable approximation which within $M$ propagates with effective minimum speed $\vmineff$.

In Section~\ref{s_veff}, we choose $\Espeed{h}(M)$ to be the event that there is a net of boxes along which the interface propagates with mesh width $\sim h$.
More precisely, we set 
\begin{align*}
	\Espeed{2h}(M)	&\coloneqq \bigcap_{x\in h\Zt\,:\,\pp{x+[-2h,3h]^2}\cap M\neq\emptyset} \Esquare[\vmin]{h,x},
\end{align*}
where $\Esquare[\vmin]{h,x}$ is the event that there are boxes $\yup+\Qup{2h}$, $\yr+\Qr{2h}$, $\yd+\Qd{2h}$, and $\yl+\Ql{2h}$ which first cut through the grid cell $x+[0,h]^2$ and then the respective neighboring cell and for which the events $\Eup[\vmin]{2h,\yup}$, $\Er[\vmin]{2h,\yr}$, $\Ed[\vmin]{2h,\yd}$, and $\El[\vmin]{2h,\yl}$ hold, see Figure~\ref{veff_fig_Egrid}.
Probability bounds as required in \ref{veff_aP_star} follow from \eqref{str_eq_pbound_exp} in a straightforward fashion.

For $(A,F)\in\Espeed{2h}(M)$, we obtain the stable approximation of a $2h$-fat set $S$ with $\partial S\subset M$ by combining stable sets from grid cells included in $S$ and close to the boundary.
More precisely, for each boundary point $y\in \partial S$ there is $y'\in\Rt$ such that $y\in\partial B_{2h}(y')$ and $\ov{B}_{2h}(y')\subset S$. Hence, the grid cell $x+[0,h]^2$ with $x\in h\Zt$ which includes $y'$ intersects $M+[-2h,2h]^2$ and its neighborhood $x+[0,h]^2+\ov{B}_{\frac{h}{2}}$ is included in $\ov{B}_{2h}(y')\subset S$.
Using that the evolution of the stable sets along the boxes associated to these grid cells preserves the stability, see Lemma~\ref{stable_lem_preserv}, and that the union of stable sets is again stable, see Lemma~\ref{stable_lem_union}, we obtain a stable set which roughly tracks the boundary of $S$ from the inside, see Figure~\ref{veff_fig_int-groups}. 
Adding the interior of $S$ enclosed by this set removes the inwards facing part of the interface -- therefore preserves the stability --  and yields our candidate for a stable approximation of $S$.
Iteratively applying the stitching lemma guarantees that $S_h$ propagates along the boxes associated to the surrounding grid cells with average speed at least $\vmin$.
Since these boxes form a net with mesh width smaller than $2h$, we obtain $S\hsubs \reach[(A,F)]{Ch}{S_h}$ for $C$ sufficiently large, see Figure~\ref{veff_fig_int-expand} and Figure~\ref{veff_fig_int-expand-small}, and thus have shown that $S_h$ is a stable approximation of $S$.

It remains to show that $(A,F)$ admits an effective minimum speed $\vmineff>0$ in $M$ for the stable approximation $S_h$. 
Suppose that $y\in\reach[(A,F)]{t}{S_h}\cap M$ for some $t>0$. 
We need to show that $\ov{B}_h(y)\hsubs \reach[(A,F)]{t+\vmineff^{-1}h}{S_h}$.
Consider the grid cell $x+[0,h]^2$ with $x\in h\Zt$ which includes $y$ and the ring of eight grid cells surrounding it.
Without loss of generality let $\dist(y,S_h)\geq 3h$, meaning that the initial set $S_h$ does not intersect these cells.
Since $S_h$ is stable, $\reach[(A,F)]{t}{S_h}$ has to cut through the ring of surrounding grid cells, otherwise $y$ could not have been reached.
In particular, $\reach[(A,F)]{t}{S_h}$ has to cut side to side through a single cell or two adjacent cells and hence also through the respective boxes associated to the cells, see Figure~\ref{veff_fig_reach-paths}.
This allows the iterative application of the stitching lemma, guaranteeing the propagation along the net of boxes and hence that the set evolution $h$-envelops $\ov{B}_h(x)$ after a fixed amount of time steps $\Delta t = 2\vmin^{-1}h$, see Figure~\ref{veff_fig_reach-expand1} and Figure~\ref{veff_fig_reach-expand2}.
This yields the desired result if $\vmineff$ is chosen accordingly.

Note that for this last step it is critical that we work in dimension $d=2$.
In principle, we expect that it is possible to translate all of the previous arguments to higher dimensions, for example by replacing the boxes with thin sheets in dimension $d=3$ and requiring that the evolution of the stable set does not just reach the top but also separates the front and back of the sheet.
However, a stable set for which we want to apply the stitching lemma would have to cut through the whole sheet, splitting it in half.
The connection between any reached point and the original set is sufficient to split a box in dimension $d=2$ but of course not a sheet in dimension $d=3$.

\section{Preliminaries}\label{s_pre}
Throughout this Section, if not specified further, we always assume $(A,F)\in\Omega$ or $(A,F)=(\wt{A},\wt{F})_Q$, that is that $(A,F)$ is the restriction of some admissible coefficient field $(\wt{A},\wt{F})\in\Omega$ to a closed set $Q\subset\Rt$ as in Definition~\ref{restr_def}.
We first state some basic properties of the set evolutions from Definition~\ref{not_def_reach}.
Then we will present some more specific properties related to stable sets and restricted coefficient fields.
While some of these have already been established in \cite{FI25}, we also introduce one of the key tools for the present work in form of Lemma~\ref{restr_lem_stitch}, our `stitching lemma', and apply it to stitch together set evolutions along overlapping boxes.
Finally, we will recall a tool to check measurability with respect to $\F(U)$ from \cite{FI25} and use it to show the measurability of  the event $\Eup[v_0]{r}$ defined in \eqref{box_E}.

\begin{remark}[Semigroup property of $\reach{t}{\cdot}$]\label{not_rem_SGprop}
Letting closed sets $S\subset\Rt$ evolve with respect to a fixed coefficient field $(A,F)$ for different times can be seen as a semigroup, that is
\begin{equation*}
	\reach[(A,F)]{t+s}{S}=\reach[(A,F)]{t}{\reach[(A,F)]{s}{S}}
	\qquad\text{for any $t,s\geq 0$,}
\end{equation*}
see e.g. \cite[Theorem~4.5.1]{Giga06}.
We will often implicitly use this property.
\end{remark}

When working with the level-set formulation based on properties of a set evolution, we will often use the following lemma to obtain a supersolution.

\begin{lemma}[{\cite[Lemma~5.2]{FI25}}]\label{not_lem_ur}
Let $S\subset\Rt$ be a closed set.
The function $\ur{S}=\ur[(A,F)]{S}$ given by
\begin{equation*}
	\ur[(A,F)]{S}(x,t)
	\coloneqq 1-\charfun[{\reach[(A,F)]{t}{S}}](x)
\end{equation*}
is a supersolution of \eqref{intro_eq_uls} with the coefficient field $(A,F)$ in $\Rt\times(0,\infty)$.
\end{lemma}

The comparison principle for the level-set formulation naturally translates to a comparison principle for the set evolutions.

\begin{lemma}[The comparison principle for set evolutions {\cite[Lemma~5.5(a)]{FI25}}] \label{not_lem_comp}
		Let $S_1\subset S_2$. 
		Then $\reach{t}{S_1}\subset \reach{t}{S_2}$ holds for all $t\geq 0$.
\end{lemma}

For long distances, there is a maximum speed at which the set evolutions can propagate. 
For shorter distances, the speed might be arbitrarily large due to curvature effects.

\begin{lemma}[Maximum speed of front propagation {\cite[Lemma~5.3]{FI25}}]\label{not_lem_vmax}
Let $S\subset\Rt$ be a closed set.
Then for any $\delta>0$ it holds that 
	$\reach[(A,F)]{T}{S}\subset S+\pp{\delta+T\vmax[\delta]}\ov{B}_1$ 
for all $T>0$, where
\begin{equation*}
 	\vmax[\delta]\coloneqq \max\Bp{C_{1F},\,6C_{1A}}+ C_{1A}\delta^{-1}.
\end{equation*} 
We write $\vmax\coloneqq \vmax[1]$.
\end{lemma}

\subsection{Properties of stable sets}\label{subs_stable}
We now state some properties for stable sets and their evolutions with respect to \eqref{intro_eq_normal}.

\begin{lemma}[Unions of stable sets]\label{stable_lem_union}
Any union of stable sets which is closed is again stable.
\end{lemma}

\begin{proof}
The result follows immediately from the fact that supersolutions are closed under the infimum: The infimum over a set of subsolutions still is a subsolution, see \cite[Lemma~2.4.1]{Giga06}.
\end{proof}

A key property of stable sets is that they do not shrink while evolving according to \eqref{intro_eq_normal} and that the stability is preserved.

\begin{lemma}[Preservation of Stability {\cite[Lemma~5.4]{FI25}}]\label{stable_lem_preserv}
Let $S\subset\Rt$ be a closed set.
If $S$ is stable with respect to $(A,F)$,
then also $\reach[(A,F)]{t}{S}$ is stable and $\reach[(A,F)]{s}{S}\subset \reach[(A,F)]{t}{S}$ for any $0\leq s \leq t$.
\end{lemma}

The evolution of a stable set is continuous in time with respect to the Hausdorff distance due to the maximum speed of propagation from Lemma~\ref{not_lem_vmax}.
For sets which are not stable, this is not the case since a connected component might fully disappear in finite time.

\begin{corollary}\label{stable_cor_HDcont}
Let $S\subset\Rt$ be a closed set.
If $S$ is stable, then $t\mapsto\reach[(A,F)]{t}{S}$ is continuous with respect to the Hausdorff distance.
\end{corollary}

\begin{proof}[Proof of Corollary~\ref{stable_cor_HDcont}]
Using the monotonicity from Lemma~\ref{stable_lem_preserv}, the semigroup property introduced in Remark~\ref{not_rem_SGprop}, and $\vmax[\delta]$ as in Lemma~\ref{not_lem_vmax}, we have for $t>s$
\begin{equation*}
	\reach{s}{S}\subset \reach{t}{S} = \reach{t-s}{\reach{s}{S}}
	\subset \bigcap_{\delta>0}\pp{\reach{s}{S} + (\delta+(t-s)\vmax[\delta])\ov{B}_1}.
\end{equation*}
\end{proof}

\subsection{Properties of restricted solutions}\label{subs_restr}
The reasoning behind restricting coefficient fields to a set $Q$ as in Definition~\ref{restr_def} is to ensure that the set evolution remains in a neighborhood of $Q$ and provides a `lower bound' for the evolution with respect to the unrestricted coefficient field.

\begin{lemma}[Properties of restricted coefficient fields {\cite[Lemma~5.10]{FI25}}]\label{restr_lem_PropRestr}
Let $u$ be a supersolution to \eqref{intro_eq_uls} with coefficients $(A,F)_{Q}$ for $(A,F)\in\Omega$ and some closed set $Q\subset\Rt$, which satisfies 
\begin{align*}
	u(\cdot,0)&\leq C	\quad\text{in }\Rt,	&
	\Bp{x\,:\, u(\cdot,0)<C} \subset \comp{\pp{\restcomp}},
\end{align*}
where
\begin{equation}\label{restr_eq_SQ}
	\restcomp\coloneqq \bigcup\Bp{B_{\frac{1}{3}}(x)\,:\,x\in\Rt\text{ with }\dist(x,Q)\geq \frac{2}{3}}.
\end{equation}
Then for any $t\geq 0$ it holds that
\begin{align*}
	u(\cdot,t)&\leq C	\quad\text{in }\Rt,	&
	\Bp{x\,:\, u(\cdot,t)<C} \subset \comp{\pp{\restcomp}}.
\end{align*}
In addition, $u$ is a supersolution also to \eqref{intro_eq_uls} with coefficients $(A,F)_{\wt{Q}}$ for any $\wt{Q}\supset Q$. 

In particular, for any $S\subset \comp{\pp{\restcomp}}$ and $\wt{Q}\supset Q$ it holds that
\begin{enumerate}[label=(\alph*), left=0pt]
\item	\label{restr_prop1_restr}
		$\reach[(A,F)_Q]{t}{S}\subset \comp{\pp{\restcomp}}\subset Q+ \ov{B}_{\frac{2}{3}}$ for any $t\geq 0$.
\item	\label{restr_prop2_stable}
		If $S$ is stable with respect to $(A,F)_Q$, then the set is stable with respect to $(A,F)_{\wt{Q}}$.
\item	\label{restr_prop3_super}
		$\ur[(A,F)_Q]{S}=1-\charfun[{\reach[(A,F)_Q]{\cdot}{S}}]$ as in Lemma~\ref{not_lem_ur} is a supersolution to \eqref{intro_eq_uls} with coefficients $(A,F)_{\wt{Q}}$.
\end{enumerate}
\end{lemma}

\begin{figure}[h]
\centering
\includegraphics[width=0.45\textwidth]{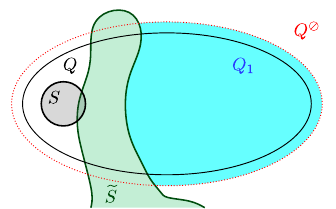}
\hspace{10pt}
\includegraphics[width=0.45\textwidth]{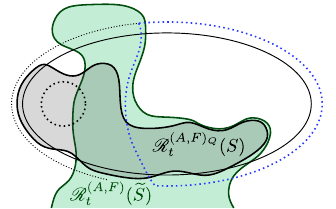}

\caption[An illustration of Lemma \ref{restr_lem_stitch}, our `stitching lemma']{	An illustration of Lemma~\ref{restr_lem_stitch}, our `stitching lemma'.
		Left: The stable set $\wt{S}$ (depicted in green) cuts through $\comp{\pp{\restcomp}}$ (depicted with red dotted line) separating the open area $Q_1$ (depicted in blue) from  $S$ (depicted in grey).
		Right: The set evolution $\reach[(A,F)_{\wt{Q}}]{t}{\wt{S}}$ (depicted in green) with respect to the less restricted coefficient field stays ahead of $\reach[(A,F)_Q]{t}{S}$ (depicted in grey) in $Q_1$ (depicted with blue dotted line) }
\label{restr_fig_stitch}
\end{figure}

Combining the properties of stable sets and our restricted coefficient fields yields a `stitching lemma', one of the key tools in the present work.
It allows us to guarantee the expansion of stable set by comparison to a set evolving with respect to a restricted coefficient field, even if there is no overlap between the sets. 
It is sufficient if the stable set cuts through the area to which the coefficient field is restricted, separating the second set from a part of that area. 
The stitching lemma allows us to compare the set evolutions in that part of the area, see Figure~\ref{restr_fig_stitch}.

\begin{lemma}[Stitching lemma]\label{restr_lem_stitch}
Let $(A,F)\in\Omega$ and $Q$ be a closed and bounded set, let $\wt{Q}\supset Q$ be a closed set.
Let $\wt{S}$ be a stable set with respect to $(A,F)_{\wt{Q}}$.
Let $Q_1\subset \comp{\pp{\restcomp}}$ be such that $\partial Q_1\subset \wt{S}\cup\partial \restcomp$, where $\restcomp$ has been defined in \eqref{restr_eq_SQ}.
Then, for any closed set $S\subset \wt{S}\cup\comp{\pp{\restcomp}}\setminus Q_1$ and $t\geq 0$ it holds that
\begin{equation*}
	\reach[(A,F)_Q]{t}{S}\cap Q_1 \subset \reach[(A,F)_{\wt{Q}}]{t}{\wt{S}}.
\end{equation*}
\end{lemma}

\begin{proof}
Let $v_{\wt{S}}\coloneqq \min\Bp{1,\,u_{\wt{S}}}$, where $u_{\wt{S}}$ as in Definition~\ref{not_def_reach} is a solution of \eqref{intro_eq_uls} for the coefficient fields $(A,F)_{\wt{Q}}$ with $\wt{S}=\Bp{x\,:\,u_{\wt{S}}(x,0)\leq 0}$.
In particular, by \cite[Theorem~A.16]{FI25}, $v_{\wt{S}}$ also is a solution.
By Lemma~\ref{restr_lem_PropRestr}, $\ur[(A,F)_Q]{S}=1-\charfun[{\reach[(A,F)_Q]{\cdot}{S}}]$ is a supersolution to \eqref{intro_eq_uls} for the coefficient fields $(A,F)_{\wt{Q}}$ and $\ur[(A,F)_Q]{S}=1$ on $\partial \restcomp$.

We apply \cite[Theorem~A.8]{FI25}, the comparison principle, to $\ur[(A,F)_Q]{S}$ and $v_{\wt{S}}$ in $Q_1$.
For this, we first need to compare the two functions on the parabolic boundary. 
Because $S\subset\wt{S}\cup\comp{\pp{Q_1}}$, we have $\ur[(A,F)_Q]{S}(\cdot,0)\geq 1-\charfun[\wt{S}] \geq  v_{\wt{S}}(\cdot,0)$ in $Q_1$.
Since $\wt{S}$ is stable, Lemma~\ref{stable_lem_preserv} yields
	$v_{\wt{S}}(\cdot,t)\leq 0\leq \ur[(A,F)_Q]{S}(\cdot,t)$ on $\wt{S}$ for $t\geq 0$. 
As mentioned above, on $\partial \restcomp\times[0,\infty)$ we have $\ur[(A,F)_Q]{S}=1\geq v_{\wt{S}}$. 
Hence, $v_{\wt{S}}\leq \ur[(A,F)_Q]{S}$ on $\partial Q_1\times[0,\infty)$.
Now, by comparison we obtain $v_{\wt{S}}\leq \ur[(A,F)_Q]{S}$ in $Q_1\times[0,\infty)$ and therefore in particular
\begin{multline*}
	\reach[(A,F)_Q]{t}{S}\cap Q_1 =\Bp{x\in Q_1\,:\, \ur[(A,F)_Q]{S}(x,t)=0}\\
	\subset \Bp{x\,:\, v_{\wt{S}}(x,t)\leq 0} = \reach[(A,F)_{\wt{Q}}]{t}{\wt{S}}.
\end{multline*}
\end{proof}

\subsection{Stitching together set evolutions along overlapping boxes}

Recall that the events $(\Eup[v]{r,x})_{r,x}$, $(\Er[v]{r,x})_{r,x}$, $(\Ed[v]{r,x})_{r,x}$, $(\El[v]{r,x})_{r,x}$ as introduced in Section \ref{subs_aP} guarantee the existence of a stable set in the starting area of the respective associated box which will propagate along that box at some average speed at least $v>0$.
We want to apply the stitching lemma to guarantee the propagation of such a stable set not just within the initial box but along networks of suitably placed overlapping boxes for which the respective events hold.

For this purpose, we introduce the notion of a set cutting vertically or horizontally through a box, which means that the set separates the respective opposite sides of that box. 
If such a set is stable, this separation will allow the application of the stitching lemma in the respective halves of the box.

\begin{definition}\label{sti_def_cut}
Let $S\subset\Rt$ be a closed set and $x\in\Rt$ with $R_1,R_2>0$.
We say that $S$ \emph{cuts vertically (horizontally) through the box} $Q=x+[0,R_1]\times[0,R_2]$ if there is no path in $Q\setminus S$ connecting the opposing sides $x+\Bp{0}\times[0,R_2]$ and $x+\Bp{R_1}\times[0,R_2]$ ($x+[0,R_1]\times\Bp{0}$ and $x+[0,R_1]\times\Bp{R_2}$).
\end{definition}

In order for a stable set in the starting area of a box associated to an event $\Eup[v]{r,x}$ to reach the opposite side, its evolution will have to cut along the main body of that box. 
Due to the maximum speed of front propagation from Lemma~\ref{not_lem_vmax} we expect that this implies that the set evolution will contain a path connecting the starting area to the top.
However, we did not easily find a proof and hence chose to fall back to the weaker notion of cuts as introduced here, which is sufficient for our purposes.

\begin{lemma}[Characterization of ${\Eup[v]{r,x}}$]\label{sti_lem_charEup}
For $x\in\Rt$, $v>0$ assume that $(A,F)\in\Eup[v]{r,x}$ as defined in \eqref{box_E}.
By definition, there exists a stable set $S\subset x+[0,w(r)]\times[-w(r),0]$ in the starting area of the box $x+\Qup{r}$ associated to the event $\Eup[v]{r,x}$ with 
\begin{equation*}
	\reach[(A,F)_{x+\Qup{r}}]{v^{-1}r}{S} \cap \pp{x+[0,w(r)]\times\Bp{r}}\neq \emptyset.
\end{equation*}
Then $\reach[(A,F)_{x+\Qup{r}}]{v^{-1}r}{S}$ cuts vertically through the main body $Q=x+[-1,w(r)+1]\times[0,r]$ of the box $x+\Qup{r}$.
\end{lemma}

\begin{proof}
Without loss of generality let $x=0$.
We prove the result by contrapositive, showing that if the set evolution does not cut through the box, then due to the comparison principle the top of the box can not have been reached.
Suppose that $\reach[(A,F)_{\Qup{r}}]{v^{-1}r}{S}$ does not cut vertically through the main body $Q$. 
Then there exists a path 
\begin{multline*}
	\gamma\in C\pp{[0,1],Q\setminus \reach[(A,F)_{\Qup{r}}]{v^{-1}r}{S}}\\
	\qquad\text{with }	\gamma(0)\in \Bp{-1}\times[0,r],~
						\gamma(1)\in \Bp{w(r)+1}\times[0,r].
\end{multline*}
Without loss of generality let $[0,w(r)]\times\Bp{r}\nsubseteq\gamma([0,1])$, otherwise we are done.
Consider the part of the thickened box $\Qup{r}+B_1$ above $\gamma$, that is, the connected component $U$ of the top of the box $[0,w(r)]\times\Bp{r}$ in $\pp{\Qup{r}+B_1}\setminus \gamma([0,1])$.
Note that $U\cap S=\emptyset$ because this would imply that the top and bottom side of the main body $Q$ are connected in $Q\setminus \gamma$ which is not possible due to the Jordan curve theorem.
In addition, $\partial U\subset \gamma([0,1])\cup \partial\pp{\Qup{r}+B_1}$ and hence 
\begin{equation*}
	\partial U\cap \reach[(A,F)_{\Qup{r}}]{t}{S}=\emptyset
	\qquad\text{for all $t\in[0,v^{-1}r]$}
\end{equation*}
because $\partial\pp{\Qup{r}+B_1}$ cannot be reached with the coefficient field restricted to $x+\Qup{r}$ by Lemma~\ref{restr_lem_PropRestr} and if $\gamma([0,1])$ had been reached at any time $t< v^{-1}r$, then those parts would still be covered at time $t=v^{-1}r$ by Lemma~\ref{stable_lem_preserv} due to the stability of $S$, contradicting the definition of $\gamma$.

We can thus compare 
\begin{align*}
	u(y,t)&=1-\charfun[{\reach[(A,F)_{\Qup{r}}]{t}{S}}](y)
	&&\text{and}
	&	v(y,t)&= 1
\end{align*}
on $U\times(0,v^{-1}r)$ with $u$ a supersolution by Lemma~\ref{not_lem_ur} and $v$ a trivial solution:
On the parabolic boundary of $U\times(0,v^{-1}r)$ we have $u=1=v$.
The comparison principle for bounded domains \cite[Theorem~A.8]{FI25} and the continuity of the set evolution from Lemma~\ref{stable_cor_HDcont} yield $u\geq 1$ on $\ov{U}\times[0,v^{-1}r]$ and therefore in particular the contrapositive
\begin{equation*}
 	\reach[(A,F)_{\Qup{r}}]{v^{-1}r}{S} \cap \pp{[0,w(r)]\times\Bp{r}}= \emptyset,
\end{equation*} 
concluding the proof.
\end{proof}

\begin{figure}[h]
\centering
\includegraphics[width=0.95\textwidth]{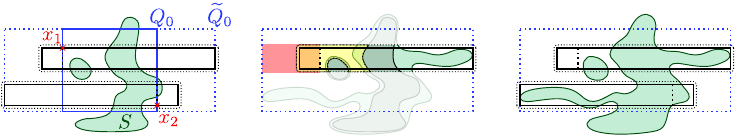}
\caption[An illustration of the box configuration in Lemma~\ref{sti_lem_base-form} and its proof]{	An illustration of the box configuration in Lemma~\ref{sti_lem_base-form} and its proof.
			Left: The set $S$ (depicted in green) cuts vertically through the box $Q_0$ (depicted in blue), which is crossed by the main parts of $x_1+\Qr{r}$ and $x_2+\Ql{r}$ (depicted in black). Together with $Q_0$, they span the box $\wt{Q}_0$ (depicted with blue dotted lines).
			Center: The layer $\wt{Q}_1\subset\wt{Q}_0$ corresponding to $x_1+\Qr{r}$ is cut horizontally by the union of the extension $Q'_1$ (depicted in red), the connected component $U_1\subset\Rt\setminus S$ of the starting area (depicted in yellow), and the set evolution $\reach[(A,F)]{v^{-1}r}{S}$ (depicted in green).
			Right:  The set evolution $\reach[(A,F)]{v^{-1}r}{S}$ (depicted in green) cuts horizontally through $\wt{Q}_0$}
\label{sti_fig_base-form}
\end{figure}

\begin{figure}[h]
\centering
\includegraphics[width=0.67\textwidth]{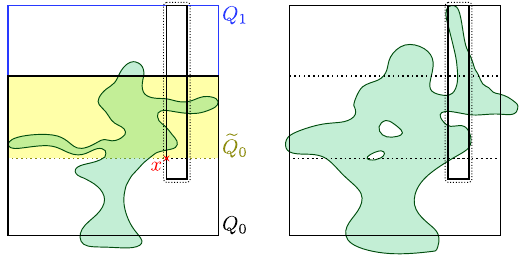}
\caption[An illustration of the box configuration in Lemma~\ref{sti_lem_extend}]{	An illustration of the box configuration in Lemma~\ref{sti_lem_extend}.
			Left: The set $S$ (depicted in green) cuts vertically through the box $Q_0$ (depicted in black) and horizontally through its top horizontal layer $\wt{Q}_0$ (depicted in yellow), which is crossed by the main part of $x+\Qup{r}$ (depicted in black). The box $Q_1$ (depicted in blue) extends $Q_0$ to the top of $x+\Qup{r}$.
			Right:  The set evolution $\reach[(A,F)]{v^{-1}r}{S}$ (depicted in green) cuts vertically through $Q_1\supset Q_0$}
\label{sti_fig_extend}
\end{figure}

In order to iteratively apply the stitching lemma along a network of overlapping boxes, we alternate horizontal and vertical boxes.
If, for example, $\Eup[v]{r}$ holds, we obtain a stable set which evolves to cut vertically along the associated box as in Lemma~\ref{sti_lem_charEup}.
For horizontal boxes crossed by this vertical box, we can hence apply the stitching lemma, see Figure~\ref{str_fig_box-stitching}.
Lemma~\ref{sti_lem_base-form} below states that, given two such horizontal boxes oriented in opposing directions along which the set evolution propagates, the original stable set will evolve to cut horizontally across the area spanned by the these two boxes, see Figure~\ref{sti_fig_base-form}.
This cut can be used as a base for an extension of the original vertical cut: 
Lemma~\ref{sti_lem_extend} states that, given an additional vertical box crossed by the horizontal cut and along which the set evolution propagates, the original stable set will further evolve to cut vertically across the area spanned by both vertical boxes, see Figure~\ref{sti_fig_extend}.

\begin{remark}
Lemma~\ref{sti_lem_base-form} and Lemma~\ref{sti_lem_extend} below remain valid for translations and rotations of the overall box configurations and when replacing the coefficient field for the set evolution with the restriction to an area including the configuration and with respect to which $S$ remains stable.
To simplify the notation, we state these two lemmas for a box configuration anchored at $x=0$ and with respect to the unrestricted coefficient field.
\end{remark}

\begin{lemma}[Cutting through two complementing boxes to form a base]\label{sti_lem_base-form}
Consider a configuration of boxes $Q_0=[0,R_1]\times[0,R_2]$, and $x_1+\Qr{r}$, $x_2+\Ql{r}$ as introduced in Section \ref{subs_aP}
where the main parts of the latter two cross $Q_0$, see Figure~\ref{sti_fig_base-form}.
That is, let 
\begin{align*}
			r	&>0, 
	&		R_1	&\in(0,r],
	&		R_2	&\geq w(r)+2,
	&		x_1	&\in	\Bp{0}\times[w(r)+1,R_2-1],\\
	&&&&&&	x_2	&\in	\Bp{R_1}\times[1,R_2-w(r)-1],
\end{align*}
where $w(r)$ is the width of the boxes $\Qr{r}$, $\Ql{r}$ as defined in \eqref{box_eq_Qup}.

Let $S\subset\Rt$ be a stable set with respect to some coefficient field $(A,F)\in\Omega$ which cuts vertically through $Q_0$.
If $(A,F)\in\Er[v]{r,x_1}\cap\El[v]{r,x_2}$ for some $v>0$, then the evolved set $\reach[(A,F)]{v^{-1}r}{S}$ cuts horizontally through the box 
	$\wt{Q}_0=[-(r-R_1),r]\times[0,R_2]$ 
spanning the main parts of $x_1+\Qr{r}$ and $x_2+\Ql{r}$.
\end{lemma}

\begin{proof}
We first consider the two horizontal layers $\wt{Q}_1,\wt{Q}_2\subset\wt{Q}_0$ corresponding to the thin boxes $x_1+\Qr{r}$ and $x_2+\Ql{r}$:
\begin{align*}
		\wt{Q}_1	&\coloneqq	x_1+[-(r-R_1),r]\times[-w(r)-1,1],\\
		\wt{Q}_2	&\coloneqq	x_2+[-r,r-R_1]\times[-1,w(r)+1].
\end{align*}
Let $U_1$ be the union of all connected components in $\pp{x_1+\Qr{r}+B_\frac{2}{3}}\setminus S$ which intersect the starting area $x_1+[-w(r),0]^2$.
Let $U_2$ be defined analogously for the box $x_2+\Ql{r}$.
Let $Q'_1$ be the extension of the main body of $x_1+\Qr{r}$ to cover the whole layer $\wt{Q}_1$ with $Q'_2$ chosen analogously for $\wt{Q}_2$, that is,
\begin{align*}
	Q'_1		&\coloneqq x_1+[-(r-R_1),0]\times[-w(r)-1,1],\\
	Q'_2		&\coloneqq x_2+[0,r-R_1]\times[-1,w(r)+1].
\end{align*}
We claim that 
\begin{align}\label{sti_eq_hor-cut-lay}
	&\reach[(A,F)]{v^{-1}r}{S}\cup U_i \cup Q_i'\text{ cuts horizontally through $\wt{Q}_i$ for $i=1,2$.}
\end{align}
We only show the claim for $i=1$ since it follows analogously for $i=2$.
First note that it is sufficient to show that $\reach[(A,F)]{v^{-1}r}{S}\cup U_1$ cuts horizontally through the main body $x_1+[0,r]\times[1,-w(r)-1]$ of $x_1+\Qr{r}$ since the extension $Q'_1$ fully covers the complement in $\wt{Q}_1$.
Let $S_1\subset\pp{x_1+[-w(r),0]^2}$ be the stable set propagating along $x_1+\Qr{r}$ whose existence is guaranteed with $(A,F)\in\Er[v]{r,x_1}$. 
By Lemma~\ref{sti_lem_charEup}, the evolved set $\reach[(A,F)_{x_1+\Qr{r}}]{v^{-1}r}{S_1}$ cuts horizontally through the main body of $x_1+\Qr{r}$. 
We hence only need to show that $\reach[(A,F)_{x_1+\Qr{r}}]{v^{-1}r}{S_1}\subset \reach[(A,F)]{v^{-1}r}{S}\cup U_1$.
This follows from an application of Lemma~\ref{restr_lem_stitch}, our stitching lemma, which yields that
\begin{align*}
	\reach[(A,F)_{x_1+\Qr{r}}]{v^{-1}r}{S_1}\cap\pp{\pp{x_1+\Qr{r}+B_{\frac{2}{3}}}\setminus U_1} \subset \reach[(A,F)]{v^{-1}r}{S},
\end{align*}
where we used that $S$ is stable, that 
	$\partial\pp{\pp{x_1+\Qr{r}+B_{\frac{2}{3}}}\setminus U_1}\subset S\cup \partial\pp{x_1+\Qr{r}+B_{\frac{2}{3}}}$, 
and that, as a subset of the starting area, $S_1\subset U_1\cup S$  by the definition of $U_1$.
After noting that $\reach[(A,F)_{x_1+\Qr{r}}]{v^{-1}r}{S_1}\subset x_1+\Qr{r}+B_{\frac{2}{3}}$ by Lemma~\ref{restr_lem_PropRestr} due to the restriction of the evolution, this establishes \eqref{sti_eq_hor-cut-lay}.

Now, suppose that $\reach[(A,F)]{v^{-1}r}{S}$ does not cut horizontally through $\wt{Q}_0$.
Then there exists a path $\gamma$ connecting the top and bottom of $\wt{Q}_0$ in $\wt{Q}_0\setminus \reach[(A,F)]{v^{-1}r}{S}$ and hence also in $\wt{Q}_0\setminus S$ due to the stability of $S$ and Lemma~\ref{stable_lem_preserv}.
In particular, there exist segments of $\gamma$ which connect the top and bottom of $\wt{Q}_1$ and $\wt{Q}_2$ inside these layers respectively. 
However, due to \eqref{sti_eq_hor-cut-lay}, this means that $\gamma$ has to pass through $U_i$ or $Q'_i$ for each $i=1,2$.
Since each connected component in $U_i$ is open, every point in $U_i$ is path-connected in $U_i\subset\wt{Q}_0\setminus S$ to the respective starting area of the associated box and hence $Q'_i$.
In summary, we obtain either a segment of $\gamma$ or an extension of such a segment along $U_1$ and/or $U_2$ which connects $Q'_1$ on the left side of $Q_0$ with $Q'_2$ on the right of $Q_0$ in $\wt{Q}_0\setminus S$.
This is a contradiction since we assumed that $S$ cuts vertically through $Q_0$.
Thus, we have shown that $\reach[(A,F)]{v^{-1}r}{S}$ cuts horizontally through $\wt{Q}_0$.
\end{proof}

\begin{lemma}[Extending a cut through a base]\label{sti_lem_extend}
Consider a configuration of boxes given by $Q_0=[0,R_1]\times[0,R_2]$, a smaller box $\wt{Q}_0=[0,R_1]\times[\wt{R}_2,R_2]$ covering its top, and a thin box $x+\Qup{r}$ as in \eqref{box_eq_Qup} with its starting area aligned below $\wt{Q}_0$, see Figure~\ref{sti_fig_extend}.
That is, let
\begin{align*}
			r		&>0, 
	&		R_1		&\geq w(r)+2,
	&		R_2		&>0,
	&		\wt{R}_2	&\in(R_2-r,R_2),\\
	&&&&&&	x		&\in	[1,R_1-w(r)-1]\times\lbrace\wt{R}_2\rbrace,
\end{align*}
where $w(r)$ is the width of the box $\Qup{r}$ as defined in \eqref{box_eq_Qup}.

Let $S\subset\Rt$ be a stable set with respect to some coefficient field $(A,F)\in\Omega$ which cuts vertically through $Q_0$ and horizontally through $\wt{Q}_0$.
If $(A,F)\in\Eup[v]{r,x}$ for some $v>0$, then the evolved set $\reach[(A,F)]{v^{-1}r}{S}$ cuts vertically through the box 
	$Q_1=[0,R_1]\times[0,\wt{R}_2+r]$, 
extending $Q_0$ by the height added via $x+\Qup{r}$.
\end{lemma}

\begin{proof}
Consider the main part $Q_2=x+[-1,w(r)+1]\times[0,r]$ of the box $x+\Qup{r}$.
Let $U$ be the union of all connected components in $\pp{x+\Qup{r}+B_\frac{2}{3}}\setminus S$ which intersect the starting area $x+[0,w(r)]\times[-w(r),0]$.
Analogously to \eqref{sti_eq_hor-cut-lay} in the Proof of Lemma~\ref{sti_lem_base-form}, we obtain that 
\begin{align}\label{sti_eq_ver-cut-lay}
	&\reach[(A,F)]{v^{-1}r}{S}\cup U \text{ cuts vertically through $Q_2$.}
\end{align}
Now, suppose that $\reach[(A,F)]{v^{-1}r}{S}$ does not cut vertically through $Q_1$.
Then there exists a path $\gamma$ connecting the left and right side of $Q_1$ in $Q_1\setminus \reach[(A,F)]{v^{-1}r}{S}$ and hence also in $Q_1\setminus S$ due to the stability of $S$ and Lemma~\ref{stable_lem_preserv}.
Since $S$ and hence by stability and Lemma~\ref{stable_lem_preserv} also $\reach[(A,F)]{v^{-1}r}{S}$ cut vertically through $Q_0$, at some point $\gamma$ has to visit the vertical extension $Q_1\setminus Q_0$.
This implies that $\gamma$ cannot visit the area $Q_0\setminus \wt{Q}_0$ below the horizontal layer $\wt{Q}_0$ since otherwise there would have to be a segment connecting the top and bottom of $\wt{Q}_0$, which would contradict the assumption that $S$ cuts horizontally through $\wt{Q}_0$.
With $\gamma$ thus staying inside the box $Q_1\setminus\pp{Q_0\setminus \wt{Q}_0}$ and connecting its left and right side, there has to be a segment of $\gamma$ which connects the left and right side of the vertical layer given by the main part $Q_2$ of the box $x+\Qup{r}$.
Due to \eqref{sti_eq_ver-cut-lay}, the path $\gamma$ hence has to visit $U$. 
Since each connected component in $U$ is open, every point in $U$ is path-connected in $U\subset Q_1\setminus S$ to the starting area of $x+\Qup{r}$.
With $\gamma$ visiting the vertical extension $Q_1\setminus Q_0$, we hence obtain a path in $Q_1\setminus S$ connecting $Q_1\setminus Q_0$ to the starting area of $x+\Qup{r}$ below the horizontal layer $\wt{Q}_0$.
This path thus connects the top to the bottom of $\wt{Q}_0$, which is a contradiction since we assumed that $S$ cuts horizontally through $\wt{Q}_0$.
Thus, we have shown that $\reach[(A,F)]{v^{-1}r}{S}$ cuts vertically through $Q_1$.
\end{proof}

\subsection{Measurability of events with respect to $\F(U)$}
We recall the measurability criterion from \cite{FI25} for the $\sigma$-algebras $\Bp{\F(U)}_U$ defined in \eqref{sett_eq_sigmaF}.

\begin{lemma}[{\cite[Lemma~5.12]{FI25}}]\label{meas_lem_closure}
Let $U\subset\Rt$ be a closed set.
Let $E\subset\Omega\subset C(\Rt\times\Sone,\Rtt\times\R)$ be an event which is closed in $\Omega$ with respect to the norm
\begin{equation*}
	\|(A,F)\|_{C_w(U\times\Sone,\Rtt\times\R)}\coloneqq \sup_{x\in U,\,e\in\Sone}\frac{1}{1+\abs{x}^2}\abs{(A(x,e),F(x,e))},
\end{equation*}
which corresponds to weighted uniform convergence and is well defined on $\Omega$ due to the uniform bounds from \eqref{sett_eq_assA} and \eqref{sett_eq_assF}.
Then $E$ is measurable with respect to $\F(U)$.
\end{lemma}

In order to show the well-posedness of the box criterion \ref{sett_aP_box}, we use Lemma~\ref{meas_lem_closure} to show the measurability of the event that there is a stable set in the starting area of a given box which reaches the opposite side before a prescribed time.

\begin{lemma}[Measurability of `propagation in a box']\label{meas_lem_Eup}
For any $r,v>0$ and $x\in\Rt$, the event $\Eup[v]{r,x}\subset\Omega$ defined in \eqref{box_E} is measurable with respect to $\F\pp{x+\Qup{r}+B_1}$, where $\Qup{r}= [0,w(r)]\times[-w(r),r]$ as in \eqref{box_eq_Qup}.
Analogous statements hold for $\Er[v]{r,x}$, $\Ed[v]{r,x}$, and $\El[v]{r,x}$.
\end{lemma}

\begin{proof}
By Lemma~\ref{meas_lem_closure}, it is enough to show that the event $\Eup[v]{r,x}\subset\Omega$ is closed with respect to $C_w(\pp{x+\Qup{r}+B_1}\times\Sone,\Rtt\times\R)$.
That is, assuming there are $(A_k,F_k)\in\Eup[v]{r,x}$ for $k\in\N$ with $(A_k,F_k)\rightarrow(A,F)$ uniformly on $\pp{x+\Qup{r}+B_1}\times\Sone$ as $k\rightarrow\infty$ for some $(A,F)\in\Omega$, we only need to show that $(A,F)\in\Eup[v]{r,x}$.
By the definition of restricted coefficient fields in \eqref{box_eq_defAQ} and \eqref{box_eq_defFQ}, we immediately have that $(A_k,F_k)_{x+\Qup{r}}\rightarrow(A,F)_{x+\Qup{r}}$ uniformly on $\Rt\times\Sone$ as $k\rightarrow\infty$.

Recalling the definition of $\Eup[v]{r,x}$, for $k\in\N$ let $S_k\subset[0,w(r)]\times[-w(r),0]$ satisfy \eqref{box_E_stable} and \eqref{box_E_top} with respect to $(A_k,F_k)$.
Let $u_k=\ur[(A_k,F_k)]{S^0_k}=1-\charfun[{\reach[(A_k,F_k)]{\cdot}{S^0_k}}]$ be the supersolutions of \eqref{intro_eq_uls} from Lemma~\ref{not_lem_ur}. 

The non-empty compact subsets of a compact set are a compact metric space when endowed with the Hausdorff metric. 
Hence, there exist (without loss of generality non-relabeled) subsequences and $S^0,\wt{S}^{v^{-1}r}\neq 0$ such that
\begin{align*}
	S_k					
	&\rightarrow S^0\subset x+\bp{0,w(r)}\times\bp{-w(r),0},\\
	\reach[{(A_k,F_k)_{x+\Qup{r}}}]{v^{-1}r}{S_k}\cap \pp{x+\bp{0,w(r)}\times\Bp{r}}
	&\rightarrow \wt{S}^{v^{-1}r}\subset \pp{x+\bp{0,w(r)}\times\Bp{r}}.
\end{align*}
We define $u\coloneqq \liminf\!{}_*u_k$, with $\liminf\!{}_*u_k$ denoting the lower relaxed limit as in \cite[(A.5)]{FI25}, which does not just take the infimum pointwise but also with respect to approaches from the surrounding neighborhood.
By Corollary~\ref{stable_cor_HDcont}, the evolutions of the starting sets are uniformly continuous with respect to the Hausdorff metric and hence 
	$u(\cdot, t) = \liminf\!{}_*\pp{1-\charfun[{\reach[(A_k,F_k)]{t}{S_k}}]}$.
Therefore, $u(\cdot,0)= 1-\charfun[S^0]$ 
and $\Bp{y\,:\, u(y,v^{-1}r)\leq 0}\supset \wt{S}^{v^{-1}r}$.
Due to the stability of viscosity solutions (see \cite[Theorem~A.7]{FI25}) and that $(A_k,F_k)_{x+\Qup{r}}\rightarrow(A,F)_{x+\Qup{r}}$ uniformly, we further know that $S_0$ is stable and that $u$ is a supersolution to \eqref{intro_eq_uls} with respect to $(A,F)_{x+\Qup{r}}$.
Applying the comparison principle \cite[Theorem~A.8]{FI25} to $u$ and $\uls{S^0}$ from Definition~\ref{not_def_reach}, we in particular obtain
\begin{equation*}
	\reach[{(A,F)_{x+\Qup{r}}}]{v^{-1}r}{S^0}
	=	\Bp{y\,:\,\uls{S^0}(y,v^{-1}r)\leq 0}
	\supset \Bp{y\,:\, u(y,v^{-1}r)\leq 0}	\supset \wt{S}^{v^{-1}r}.
\end{equation*}
Thus, $(A,F)\in\Eup[v]{r,x}$ with $S^0$ satisfying \eqref{box_E_stable} and \eqref{box_E_top}, showing that $\Eup[v]{r,x}$ is closed and therefore measurable.
\end{proof}

\section{Lifting the box criterion to larger scales}\label{s_vmin}
In this section, we aim to lift the probability bound \eqref{sett_eq_assBox} from the box criterion \ref{sett_aP_box} to arbitrary scales. 
That is, for arbitrary $r\geq r_0$ we aim to obtain probability bounds for the events that the interface propagates along the boxes $\Qup{r}$, $\Qr{r}$, $\Qd{r}$, and $\Ql{r}$ respectively, see Figure~\ref{intro_fig_box}.
We are looking to find speeds $v_i>0$ and probability bounds $p_i\in[0,1)$ such that for all $i\in\N$
\begin{align}\label{vmin_eq_assBoxScaled}
	\min\Bp{\PM{\Eup[v_i]{r}},\PM{\Er[v_i]{r}},\PM{\Ed[v_i]{r}},\PM{\El[v_i]{r}}}&\geq 1-p_i
	&&	\text{for all }r\in[r_i,r_{i+1}),
\end{align}
where the sequence $(r_i)_{i\in\N}$ is recursively defined by 
\begin{equation}\label{vmin_eq_ri}
	r_i\coloneqq w^{-1}\pp{r_{i-1}+2w\pp{r_{i-1}}}.
\end{equation}
This means that the width of the box $\Qup{r_i}$ from Section~\ref{subs_box} is given by $r_{i-1}+2w(r_{i-1})$.
Therefore, translations of $\Qr{r_{i-1}}$ and $\Ql{r_{i-1}}$ fit inside $\Qup{r_i}$ without their starting areas overlapping, see Figure~\ref{sett_fig_r1}.

This setup already suggests that we will obtain \eqref{vmin_eq_assBoxScaled} iteratively by arranging smaller boxes, for which we will already have a probability bound, inside a given box and exploiting the random fields' stationarity and finite range of dependence.

\begin{proposition}[The box criterion for large scales]\label{vmin_prop_assBoxScaled}
Assume that \ref{sett_aP_stat}, \ref{sett_aP_fin}, and \ref{sett_aP_box} hold. 
Let $(r_i)_{i\in\N}$ be recursively defined box lengths as in \eqref{vmin_eq_ri}.
Let $0<\delta<\frac{q-1}{q}$.
Then there are sequences $(p_i)_{i\in\N}$ and $(v_i)_{i\in\N}$ such that \eqref{vmin_eq_assBoxScaled} holds for all $i\in\N$ with the sequences satisfying
\begin{align}\label{vmin_eq_vpBounds}
	p_{i+1}
	&\leq	\min\Bp{p_0, p_i^{\pp{\frac{r_i}{r_0}}^{\frac{q-1}{q}-\delta}-1}},
	&
	v_{i+1} 	
	&\geq	\pp{\frac{1}{1	+ \Cv r_{i+1}^{-\min\Bp{\delta,~ \frac{q-1}{q}-\delta}}}}^{\ell} v_i,
\end{align}
where $\Cv=\Cv(q,\cw,r_0,\delta)$ and $\ell=\ell(q,\cw,r_0,\delta)$ depend only on $q,\cw,r_0,\delta$.
\end{proposition}


Once we obtain Proposition~\ref{vmin_prop_assBoxScaled}, relating the probability bounds to the scale of the boxes yields stretched exponential bounds for the probability that propagation along the boxes fails.

\begin{corollary}\label{vmin_cor_assBoxScaled}
Assume that \ref{sett_aP_stat}, \ref{sett_aP_fin}, and \ref{sett_aP_box} hold.
For $\ratebox \in \pp{0,\frac{q-1}{q^3}}$, there exist $\Cbox=\Cbox(v_0,q,\cw,r_0,p_0,\ratebox)>0$ and a velocity $\vmin=\vmin(q,\cw,r_0,v_0,\ratebox)$ such that 
\begin{align}\label{vmin_eq_assBox_vmin}
				\min\Bp{\PM{\Eup[\vmin]{r}},\PM{\Er[\vmin]{r}},\PM{\Ed[\vmin]{r}},\PM{\El[\vmin]{r}}}&\geq 1-\exp\pp{-\Cbox r^{\ratebox}}
\end{align}
for all $r\geq r_0$.
\end{corollary}

\begin{proof}[Proof of Corollary~\ref{vmin_cor_assBoxScaled}]
Let $\ratebox \in \pp{0,\frac{q-1}{q^3}}$. 
Choose $\delta\in\pp{0,\frac{q-1}{q}}$ such that $\ratebox = \frac{q-1}{q^3}-\frac{\delta}{q^2}$.
Given $\delta$, let $(p_i)_i$, $(v_i)_i$ be as in Proposition~\ref{vmin_prop_assBoxScaled}. 
Using \eqref{vmin_eq_vpBounds}, we will first show that $\vmin=\vmin(q,\cw,r_0,v_0,\delta)\coloneqq \liminf_{i\rightarrow\infty} v_i >0$ and then provide the probability bound in \eqref{vmin_eq_assBox_vmin} based on \eqref{vmin_eq_assBoxScaled} and the bounds for $(p_i)_i$ in \eqref{vmin_eq_vpBounds}.

Regarding $\vmin$, we set $\beta\coloneqq \min\Bp{\delta,~\frac{q-1}{q}-\delta}>0$.
From \eqref{vmin_eq_vpBounds}, for any $i\in\N$, we obtain
\begin{equation*}
	v_i	\geq		\pp{\frac{1}{1+\Cv r_{i}^{-\beta}}}^\ell v_{i-1}
		\geq		\pp{\prod_{j=1}^{i}\frac{1}{1+\Cv r_j^{-\beta}}}^\ell v_{0}
		\geq		\exp\pp{-\Cv\ell \sum_{j=1}^i r_j^{-\beta}}v_{0},
\end{equation*}
where we used that $1+x\leq \exp(x)$.
Note that, using the formula for geometric sums,
\begin{multline*}
	r_i	=		w^{-1}\pp{r_{i-1}+2w(r_{i-1})}
		\geq 		w^{-1}\pp{r_{i-1}}\\
		=			\cw r_{i-1}^q
		\geq		\cw^{\sum_{j=0}^{i-1}q^j}r_0^{q^i}
		=		\cw^{-\frac{1}{q-1}}\pp{\cw^{\frac{1}{q-1}}r_0}^{q^i}.
\end{multline*}
Plugging this into the bound for $v_i$ above yields
\begin{align*}
	v_i	\geq		\exp\pp{-\Cv\ell\cw^{\frac{\beta}{q-1}} \sum_{j=1}^i \pp{\cw^{\frac{1}{q-1}}r_0}^{-\beta q^j}}v_{0}.
\end{align*}
For $\vmin=\liminf_{i\rightarrow\infty} v_i$, we hence obtain
\begin{align*}
	\vmin	
		&\geq		\exp\pp{-\Cv\ell\cw^{\frac{\beta}{q-1}} \sum_{j=1}^\infty \pp{\cw^{\frac{1}{q-1}}r_0}^{-\beta q^j}}v_{0}
		>0,
\end{align*}
where we used that $q>1$, $\beta>0$ and $\cw^{\frac{1}{q-1}}r_0 = \pp{\cw^\frac{1}{q}r_0^{\frac{q-1}{q}}}^{\frac{q}{q-1}}\geq 24^{\frac{q}{q-1}}>1$ due to \eqref{sett_eq_assBox_r} from Assumption~\ref{sett_aP_box}.

Regarding the probability bound, let $r\in[r_{i+1},r_{i+2})$, $i\geq 1$.
Then, since $\vmin\leq v_{i}$ for all $i\in\Nz$ by definition, \eqref{vmin_eq_vpBounds} yields the bound
\begin{align*}
	1-\PM{\Eup[\vmin]{r}} 
	\leq		p_{i+1} 
	&\leq		p_0^{\pp{\frac{r_i}{r_0}}^{\frac{q-1}{q}-\delta}-1}\\
	&\leq		p_0^{\pp{1-\pp{\frac{r_0}{r_1}}^{\frac{q-1}{q}-\delta}}\frac{r_i}{r_0}^{\frac{q-1}{q}-\delta}}\\
	&\leq		\exp\pp{-\log\pp{p_0^{-1}}\pp{r_0^{\delta-\frac{q-1}{q}}-r_1^{\delta-\frac{q-1}{q}}}r_i^{\frac{q-1}{q}-\delta}}.
\end{align*}
where we used that $r_i\geq 2r_0$ for $i\geq 1$.
Since we have $r_j\geq \frac{1}{2}w\pp{r_{j+1}}$  by the definition of the sequence $(r_j)_j$ from \eqref{vmin_eq_ri}, we have
\begin{equation*}
	r_i	\geq 	\frac{1}{2}w\pp{\frac{1}{2}w\pp{r_{i+2}}}
		\geq		\frac{1}{2}w\pp{\frac{1}{2}w\pp{r}}
		=		\pp{\frac{1}{2\cw^{\frac{1}{q}}}}^{\frac{q+1}{q}}r^{\frac{1}{q^2}}
\end{equation*}
and we hence obtain a constant $C(q,\cw,r_0,p_0,\delta)>0$ such that
\begin{equation*}
	1-\PM{\Eup[\vmin]{r}} 
	\leq 
	\exp\pp{-C(q,\cw,r_0,p_0,\delta)r^{\frac{q-1}{q^3}-\frac{\delta}{q^2}}}
	\eqqcolon \exp\pp{-\Cbox r^{\ratebox}}.	
\end{equation*}
\end{proof}

\begin{proof}[Proof of Proposition~\ref{vmin_prop_assBoxScaled}]
We obtain the result by induction. 
The box criterion \ref{sett_aP_box} yields \eqref{vmin_eq_assBoxScaled} for $i=0$. 
Assuming that \eqref{vmin_eq_assBoxScaled} holds for some $i\in\Nz$ with $v_i>0$ and $p_i>0$ small enough, we will show that \eqref{vmin_eq_assBoxScaled} also holds for $i+1$ with some $p_{i+1}, v_{i+1}$ satisfying \eqref{vmin_eq_vpBounds}.
The requirements on $r_0$ from \eqref{sett_eq_assBox_r} and $p_0$ from \eqref{sett_eq_assBox_p} in \ref{sett_aP_box} are chosen to ensure that $p_0$ is small enough and via \eqref{vmin_eq_vpBounds} also $p_i$ for all $i\in\N$.

To be more precise, we will proceed in seven steps as follows:
\begin{enumerate}[label=\arabic*.)]
\item	Given $r\geq r_1$, we introduce a configuration of smaller boxes inside the box $\Qup{r}$ consisting of $N$ layers. 
		Each layer consists of one subgroup of $\kvert$ parallel vertical boxes of length $\rup=\rup(r,\kvert)$ and two subgroups of $\khor$ parallel horizontal boxes of length $\rbase=\rbase(r)$. 
\item	Assuming that for each subgroup there is at least one box along which the interface propagates with minimum average speed $v$, using the stitching lemma we show that $\Eup[\wt{v}]{r}$ holds for suitable $\wt{v}=\wt{v}(r,\kvert,\khor)$.
\item	From the previous step, with the stationarity from Assumption~\ref{sett_aP_stat} and the finite range of dependence from Assumption~\ref{sett_aP_fin}, we deduce a probability bound for $\PM{\Eup[\wt{v}]{r}}$ depending on $\PM{\Eup[v]{\rup}}$, $\PM{\Er[v]{\rbase}}$, and $\PM{\El[v]{\rbase}}$, as well as $\kvert$ and $\khor$.
\item	We discuss the requirements on the number of boxes per subgroup and fix a choice for $\kvert= \kvert(r)$ and $\khor=\khor(r)$.
\item	Having fixed $\kvert$ and $\khor$, we quantify the loss of speed from $v$ to $\wt{v}$ from Step 2.
\item	In the same vein, we ensure that the probability bound from Step 3 is an improvement compared to the smaller boxes, explaining the bounds for $r_0$ in \eqref{sett_eq_assBox_r} and $p_0$ in \eqref{sett_eq_assBox_p} from Assumption~\ref{sett_aP_box}. 
\item	Repeatedly applying this argument yields the result. Note that, assuming \eqref{vmin_eq_assBoxScaled} holds for some $i\in\N$, multiple intermediate steps are necessary in order to expand the probability bounds from $r_{i+1}$ to $r_{i+2}$: while $\rbase(r)\in[r_i,r_{i+1})$ for $r\in [r_{i+1},r_{i+2})$, for $r$ closer to $r_{i+2}$ we have $\rup(r)>r_{i+1}$.
\end{enumerate}
Repeating the procedure for $\Qr{r}$, $\Qd{r}$, and $\Ql{r}$ yields the result.

\underline{Step 1: Setting up smaller boxes in  $\Qup{r}$.}
Given $r\geq r_1$, we set up multiple layers of smaller boxes within the large vertical box $\Qup{r}$, see Figure~\ref{vmin_fig_box-config}.
Except for the first layer, each consists of three groups of boxes:
\begin{itemize}
\item	Two groups of horizontal boxes form the `base' of the layer. 
		The first group consists of $\khor\in\N$ translations of $\Qr{\rbase}$, where $\rbase>0$ is chosen such that $\rbase+2w(\rbase)=w(r)$. 
		These boxes are stacked on top of each other with distance $3$ between them, their starting areas aligning with the left side of $\Qup{r}$.
\item	The second group respectively consists of $\khor\in\N$ translations of $\Ql{\rbase}$.
		These are stacked overlapping with the first group, but such that their starting areas align with the right side of $\Qup{r}$.
		Hence, each box from the second group overlaps one from the first, except for their starting areas to the right/left respectively.
\item	The third group consists of $\kvert\in\N$	translations of the vertical box $\Qup{\rup}$, which are stacked next to each other with distance $3$ between them. 
		Here, $\rup$ is chosen such that the overall width of the group is less than $\rbase-2$.
		This group of vertical boxes is placed such that the top of their starting areas is $1$ below the lowest of the horizontal boxes.
\end{itemize}
For the first layer, we scrap the horizontal base, leaving only the vertical boxes.
This first layer of vertical boxes is placed such that the top of their starting areas is aligned with the top of the starting area $[0,w(r)]\times[-w(r),0]\subset\Qup{r}$.
Each consecutive layer is then placed such that the base is $1$ below the top of the vertical boxes from the previous layer.
For the last layer, if necessary we adjust the position of the vertical boxes downwards until their top aligns with the top of $\Qup{r}$.

\begin{figure}[h]
\centering
\includegraphics[width=0.95\textwidth]{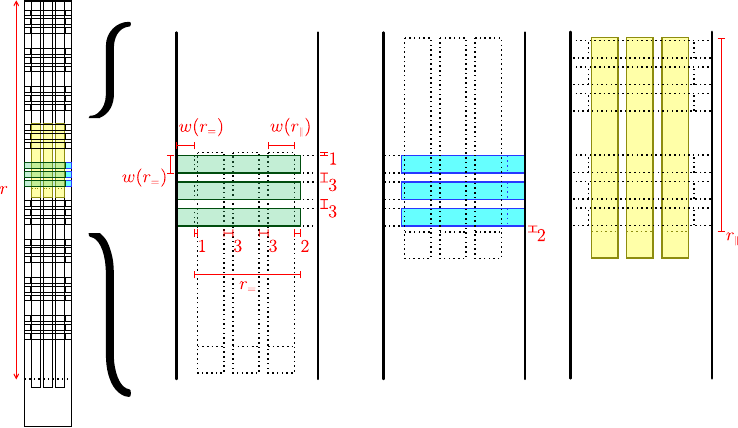}
\caption[An illustration of the box configuration introduced in Step 1 of the proof of Proposition \ref{vmin_prop_assBoxScaled}]{	An illustration of the box configuration introduced in Step 1 of the proof of Proposition~\ref{vmin_prop_assBoxScaled}.
			Left: The box $\Qup{r}$ with layers each consisting of three subgroups of smaller boxes stacked inside.
			Center left: For a single layer, the first subgroup consists of $\khor$ parallel boxes pointing to the right (depicted in green), cutting through the top of the vertical subgroup from the previous layer (depicted with dotted lines). 
			Center right: The second subgroup (depicted in blue) mirrors the first subgroup.
			Right:  The third subgroup consists of $\kvert$ parallel vertical boxes (depicted in yellow) spaced with distance $3$ between them and with their starting areas placed below the previous subgroups}
\label{vmin_fig_box-config}
\end{figure}

To be precise, given $r\geq r_1$ and $\khor,\kvert\in\N$, since $w\colon[0,\infty)\rightarrow[0,\infty)$ is onto and strictly increasing, there exist $\rbase\geq r_0$ and -- as along as $r$ is large enough compared to $\kvert$ -- further $\rup\geq r_0$ such that
\begin{align}\label{vmin_eq_rup_rvert}
	w(r)		&=\rbase+2w(\rbase),
	&
	\rbase 	&=\kvert w(\rup)+3\kvert.
\end{align}
The number of layers needed to reach the top of $\Qup{r}$ is 
\begin{equation}\label{vmin_eq_N}
	N\coloneqq \left\lceil\frac{r-\khor w(\rbase)-3\khor}{{\rup-\khor w(\rbase)-3\khor}} \right\rceil
\end{equation}
since for the height `gained' per layer (after the first one) we have to subtract the base of the layer, i.e. the horizontally stacked boxes.
The $n$-th layer, $n\in\Bp{0,\ldots,N-1}$, consists of the boxes 
\begin{align*}
	\xr{n}{j}+\Qr{\rbase},&\quad\text{for }j\in\Bp{0,\ldots,\khor-1},		\qquad\text{(only for $n\neq 0$)}
	\\
	\xl{n}{j}+\Ql{\rbase},&\quad\text{for }j\in\Bp{0,\ldots,\khor-1},		\qquad\text{(only for $n\neq 0$)}
	\\
	\xup{n}{j}+\Qup{\rup},&\quad\text{for }j\in\Bp{0,\ldots,\kvert-1},
\end{align*}
where, see Figure~\ref{vmin_fig_box-config},
\begin{align*}
	\xr{n}{j}	&\coloneqq		\begin{pmatrix}	w(\rbase)	\\	w(\rbase)+2					\end{pmatrix}
							+j	\begin{pmatrix}	0			\\	w(\rbase)+3					\end{pmatrix}
							+n	\begin{pmatrix}	0			\\	\rup-\khor\pp{w(\rbase)+3}	\end{pmatrix},\\
	\xl{n}{j}	&\coloneqq		\begin{pmatrix}	w(r)-w(\rbase)	\\	2							\end{pmatrix}
							+j	\begin{pmatrix}	0				\\	w(\rbase)+3					\end{pmatrix}
							+n	\begin{pmatrix}	0				\\	\rup-\khor\pp{w(\rbase)+3}	\end{pmatrix},\\
	\xup{n}{j}	&\coloneqq	\!\left\lbrace		\begin{aligned}
							&	\begin{pmatrix}	w(\rbase)+1	\\	0							\end{pmatrix}
							+j	\begin{pmatrix}	w(\rup)+3	\\	0							\end{pmatrix}
							+n	\begin{pmatrix}	0			\\	\rup-\khor\pp{w(\rbase)+3}	\end{pmatrix}
							&&\!\!\text{for }n<N\!-\!1,\\
							&	\begin{pmatrix}	w(\rbase)+1	\\	r-\rup	\end{pmatrix}
							+j	\begin{pmatrix}	w(\rup)+3	\\	0		\end{pmatrix}
							&&\!\!\text{for }n=N\!-\!1.
							\end{aligned}	\right.
\end{align*}

\underline{Step 2: $\Eup[\wt{v}]{r}$ holds for some $\wt{v}$, if enough of the smaller boxes allow propagation.}
Still assuming that $\khor,\kvert\in\N$ are already chosen, we require that
\begin{itemize}
\item	$r$ is large enough such that $\rbase,\rup$ are well defined,
\item	$0<N<\infty$, that is $\rup-\khor w(\rbase)-3\khor>0$,
\end{itemize}
In this step, we will show that if for each layer and for each subgroup there is at least one box along which the interface propagates with average speed $v>0$, then $\Eup[\wt{v}]{r}$ holds for some $\wt{v}=\wt{v}(v,\kvert,\khor,r)>0$:
\begin{align}\label{vmin_eq_EventImplication}
	\Eup[\wt{v}]{r}
	\supset
	\bigcap_{n=0}^{N-1}				\pp{\bigcup_{j=0}^{\kvert-1}	\Eup[v]{\rup,\xup{n}{j}}}
	\cap		\bigcap_{n=1}^{N-1}	\pp{\bigcup_{j=0}^{\khor-1}	\Er[v]{\rbase,\xr{n}{j}}}
	\cap		\bigcap_{n=1}^{N-1}	\pp{\bigcup_{j=0}^{\khor-1}	\El[v]{\rbase,\xl{n}{j}}}.
\end{align}
Hence, we will show that if for each $n\in\Bp{0,\ldots,N-1}$ there exist $\jup{n}\in \Bp{0,\ldots,\kvert-1}$, and $\jr{n}, \jl{n}\in\Bp{0,\ldots,\khor-1}$ (for $n>0$), such that $\Eup[v]{\rup,\xup{n}{\jup{n}}}$, $\Er[v]{\rbase,\xr{n}{\jr{n}}}$, and $\El[v]{\rbase,\xl{n}{\jl{n}}}$ hold, then $\Eup[v]{r}$ holds with $\wt{v}$ defined below in \eqref{vmin_eq_def_v}.

Since $\Eup[v]{\rup,\xup{0}{\jup{0}}}$ holds, there exists a closed set 
\begin{equation*}
	S_0\subset \pp{\xup{0}{\jup{0}}+\bp{0,w(\rup)}\times\bp{-w(\rup),0}}\subset \bp{0,w(r)}\times\bp{-w(r),0}
\end{equation*}
which is stable with respect to $(A,F)_{\xup{0}{\jup{0}}+\Qup{\rup}}$ and therefore by Lemma~\ref{restr_lem_PropRestr} also stable with respect to the less restricted coefficient field $(A,F)_{\Qup{r}}$.
Thus, $S_0$ satisfies \eqref{box_E} and \eqref{box_E_stable} for $\Eup[\wt{v}]{r}$.

In order to show that $S_0$ satisfies \eqref{box_E_top}, we use Lemma~\ref{sti_lem_base-form} and Lemma~\ref{sti_lem_charEup} to iteratively ``stitch'' together the guaranteed set evolutions along the the various small boxes, see Figure~\ref{vmin_fig_stitch}.
We will show that $\reach[(A,F)_{\Qup{r}}]{t}{S_0}$ passes the $n$-th layer, $n\in\Bp{0,\ldots,N-1}$, at the latest after time
\begin{equation*}
	t_n\coloneqq v^{-1}\rup + nv^{-1}\pp{\rbase+\rup}.
\end{equation*}
More precisely, we will show by induction that 
\begin{multline}\label{vmin_eq_Qn}
	\reach[(A,F)_{\Qup{r}}]{{t_n}}{S_0} \text{ cuts vertically through the box}\\ 
	\text{$Q_n\coloneqq[w(\rbase),w(r)-w(\rbase)]\times[0,\min\Bp{r,(n+1)\rup-n\khor\pp{w(\rbase)+3}}]$}
\end{multline}
in the sense of Definition \ref{sti_def_cut} for $n\in\Bp{0,\ldots,N-1}$.
For $n=N-1$, this implies that $\reach[(A,F)_{\Qup{r}}]{t_{N-1}}{S_0}$ must have reached the top of $\Qup{r}$.
Hence, $S_0$ satisfies \eqref{box_E_top} with 
\begin{equation}\label{vmin_eq_def_v}
	\wt{v}	\coloneqq 	\frac{r}{t_{N-1}}
			=			v\frac{r}{(N-1)\pp{\rup+\rbase}+\rup}.
\end{equation}
Thus, $\Eup[\wt{v}]{r}$ and therefore $\Eup[\hat{v}]{r}$ for any $\hat{v}\leq \wt{v}$ hold.

\begin{figure}[h]
\centering
\includegraphics[width=0.95\textwidth]{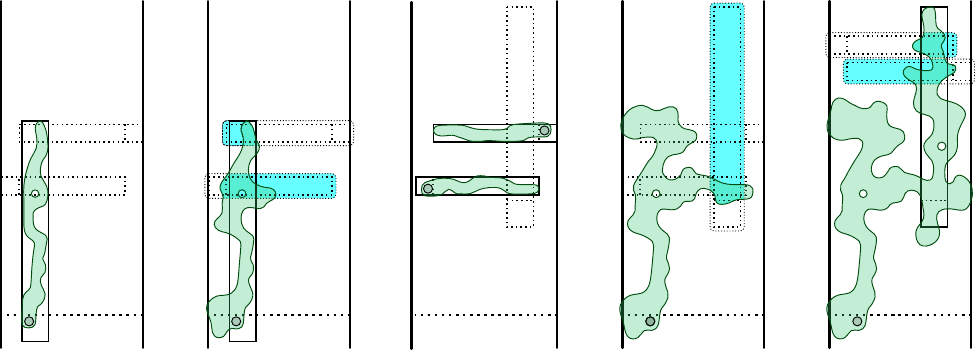}
\caption[An illustration of the stitching argument from Step 2 of the proof of Proposition \ref{vmin_prop_assBoxScaled}]{	An illustration of the stitching argument from Step 2 of the proof of Proposition~\ref{vmin_prop_assBoxScaled}.
			Left: After time $v^{-1}\rup$, the restricted evolution (depicted in green) of the set $S_0$ (depicted in grey) associated to one of the vertical boxes in the first layer cuts through the horizontal boxes from the second layer.
			Center left: This evolution is contained in the less restricted evolution $\reach[(A,F)_{\Qup{r}}]{v^{-1}\rup}{S_0}$ (depicted in green), which hence cuts vertically through the horizontal boxes separating their end (depicted in blue) from their starting area.
			Center: The stable sets associated to the horizontal boxes (depicted in grey)  after evolving with respect to the restriction for time $v^{-1}\rbase$ (depicted in green) cut horizontally through their boxes.
			Center right: Due to the stitching lemma, $\reach[(A,F)_{\Qup{r}}]{v^{-1}\rup+v^{-1}\rbase}{S_0}$ (depicted in green) cuts horizontally through the vertical box of the second layer along which the interface propagates, separating its top (depicted in blue) from its starting area. 
			Right: Another application of the stitching lemma yields that $\reach[(A,F)_{\Qup{r}}]{2v^{-1}\rup+v^{-1}\rbase}{S_0}$ (depicted in green) extends the vertical cut through the horizontal boxes of the third layer}
\label{vmin_fig_stitch}
\end{figure}

It remains to prove \eqref{vmin_eq_Qn}.
For $n=0$, due to the characterization of $\Eup[v]{\rup}$ from Lemma~\ref{sti_lem_charEup}, the set evolution of $S_0$ restricted to $\xup{0}{\jup{0}}+\Qup{\rup}$ cuts vertically through the main part of that box after time $t_0=v^{-1}\rup$.
By Lemma~\ref{restr_lem_PropRestr}, the same is hence true for $\reach[(A,F)_{\Qup{r}}]{{t_0}}{S_0}$ after evolving in a less restricted environment and for the wider box $Q_0$, for which the main part of $\xup{0}{\jup{0}}+\Qup{\rup}$ is a vertical layer. 

Assuming that \eqref{vmin_eq_Qn} holds for some $n\in\Bp{0,\ldots,N-2}$, we know that $\reach[(A,F)_{\Qup{r}}]{{t_n}}{S_0}$ cuts vertically through the box 
\begin{equation*}
	\wt{Q}_n= \begin{pmatrix}	w(\rbase)	\\	(n+1)\pp{\rup-\khor\pp{w(\rbase)+3}}	\end{pmatrix} + [0,w(r)-2w(\rbase)]\times[0,\khor\pp{w(\rbase)+3}]
\end{equation*}
which includes the main bodies of ${\xr{n+1}{\jr{(n+1)}}+\Qr{\rbase}}$ and ${\xl{n+1}{\jl{(n+1)}}+\Ql{\rbase}}$, the horizontal boxes of the $(n+1)$-th layer.
An application of Lemma~\ref{sti_lem_base-form} hence yields that $\reach[(A,F)_{\Qup{r}}]{{t_n+v^{-1}\rbase}}{S_0}$ cuts horizontally through $\wt{Q}_n$.
Since the bottom side of $\wt{Q}_n$ is aligned with the starting area of $\xup{n+1}{\jup{(n+1)}}+\Qup{\rup}$, the vertical box of the $(n+1)$-th layer, we can apply Lemma~\ref{sti_lem_extend}.
With $Q_{n+1}$ corresponding to the extension of $Q_n$ by the height of this vertical box, we obtain that $\reach[(A,F)_{\Qup{r}}]{{t_n+v^{-1}\rbase+v^{-1}\rup}}{S_0}$ cuts vertically through $Q_{n+1}$.
By induction, we have thus shown that \eqref{vmin_eq_Qn} holds for all $n\in\Bp{0,\ldots,N-1}$.

\underline{Step 3: A probability bound.}
In the previous step we showed \eqref{vmin_eq_EventImplication}, that is, that $\Eup[\wt{v}]{r}$ holds if for each of the $N$ groups of vertical boxes and $2(N-1)$ groups of horizontal boxes there is at least one box such that the corresponding event holds for speed $v$.
Using that $\PM{E_1\cap E_2}\geq 1 - \PM{\comp{(E_1)}}-\PM{\comp{(E_2)}}$ and \ref{sett_aP_stat}, the stationarity of the probability distribution, \eqref{vmin_eq_EventImplication} yields
\begin{multline*}
	\PM{\Eup[\wt{v}]{r}}	\geq
	1	-	N\,\PM{\bigcap_{j=0}^{\kvert-1}	\comp{\pp{\Eup[v]{\rup,\xup{0}{j}}}}}\\
	-	(N-1) \PM{\bigcap_{j=0}^{\khor-1}\comp{\pp{\Er[v]{\rbase,\xr{1}{j}}}}}
	-	(N-1) \PM{\bigcap_{j=0}^{\khor-1}\comp{\pp{\El[v]{\rbase,\xl{1}{j}}}}}.
\end{multline*}
Since the distance between any two boxes in one group is larger than $3$, and since by Lemma~\ref{meas_lem_Eup} whether the interface propagates along the box only depends on the coefficient field in a $1$-neighborhood, the corresponding events are independent due to the finite range of dependence from Assumption~\ref{sett_aP_fin}.
Hence, again invoking the stationarity of the probability distribution, we obtain
\begin{align}
	\PM{\Eup[\wt{v}]{r}}	
	&\geq	1	-	N \PM{\comp{\pp{\Eup[v]{\rup}}}}^{\kvert}
			-	N \PM{\comp{\pp{\Er[v]{\rbase}}}}^{\khor}
			-	N \PM{\comp{\pp{\El[v]{\rbase}}}}^{\khor}		\notag\\
	&= 	1	-	N \pp{1-\PM{\Eup[v]{\rup}}}^{\kvert} -N\pp{1-\PM{\Er[v]{\rbase}}}^{\khor}-N\pp{1-\PM{\El[v]{\rbase}}}^{\khor}.
	\label{vmin_eq_ProbBound}
\end{align}

\underline{Step 4: Choosing the parameters $\kvert$ and $\khor$.}
Given $r\in [r_i,r_{i+1})$ for some $i\in\N$, in particular we have $r\geq r_1$. 
Hence, for the length $\rbase=\rbase(r)$ of the horizontal boxes as defined in \eqref{vmin_eq_rup_rvert} we automatically obtain 
\begin{enumerate}[label=(C\arabic*)]
\item	\label{vmin_lab_rbase}
		$\rbase=\rbase(r)\in [r_{i-1},r_{i})$ and $\rbase=\rbase(r)\geq r_0$ by the definition of $(r_i)_i$ in \eqref{vmin_eq_ri}.
\end{enumerate}
Considering \eqref{vmin_eq_ProbBound}, we set $\khor=\kvert=k$ with the intent to choose $k=k(r)\geq 2$ as large as possible. 
For $\rbase$ and the length $\rup=\rup(r)$ of the vertical boxes as defined in \eqref{vmin_eq_rup_rvert} we still want to guarantee though that 
\begin{enumerate}[resume*]
\item	\label{vmin_lab_rup}
		$\rup=\rup(r)$ exists with $\rbase\leq \rup \leq k^{-q}r$ (the last inequality holds by definition),
\item	\label{vmin_lab_three}
		$w(\rbase)\geq 3$,
\item	\label{vmin_lab_N}
		$N\in \bp{k^q,2k^q}$ (by definition $N\geq k^q$, the constant $2$ is chosen arbitrarily for an initial bound on the loss of speed\footnote{
				In order to obtain a positive limit speed with this iteration, we need $\frac{\wt{v}}{v}\rightarrow 1$ as $r\rightarrow\infty$. 
				This requires $\frac{N\rbase}{r}$, the additional time corresponding to traversing the horizontal boxes relative to $r$, to disappear. 
				Since $N\geq k^q$, this implies $\frac{k^q}{\rbase^{q-1}}\rightarrow 0$ and thus plugging this into \eqref{vmin_eq_N} also $\frac{N}{k^q}\rightarrow 1$ as $r\rightarrow\infty$.
				Requiring that $\frac{k^q}{\rbase^{q-1}}\rightarrow 0$ also explains why below we choose $k\sim \rbase^{\frac{q-1}{q}-\delta}$ and not $k\sim \rbase^{\frac{q-1}{q}}$.
			}).
\end{enumerate}
First, \ref{vmin_lab_three} is equivalent to 
\begin{equation}
 	\rbase \geq 3^q\cw.\label{vmin_eq_C3}
 \end{equation} 
With \eqref{vmin_eq_rup_rvert}, the definition of $\rup$, existence is guaranteed for $\rbase>3k$.
Using \ref{vmin_lab_three}, we obtain that the inequality $\rbase\leq \rup$ in \ref{vmin_lab_rup} is implied via
\begin{align}
	&&	\rup			&	=			w^{-1}\pp{\frac{\rbase}{k}-3} \geq \rbase
	&				&	\Leftarrow
	&	\rbase		&	\geq 		k\pp{w(\rbase)+3}\notag\\
					&	\Leftarrow
	&	\rbase		&	\geq 		2k w(\rbase)
	&				&	\Leftarrow
	&	\rbase^{q-1}	&	\geq 		2^q k^q\cw^{-1}.\label{vmin_eq_C2}
\end{align}
The second inequality, $\rup \leq k^{-q}r$, follows from \eqref{vmin_eq_rup_rvert}, the definition of $\rbase$, via 
\begin{equation*}
	\rup		=		w^{-1}\pp{\frac{\rbase}{k}-3} \leq w^{-1}\pp{\frac{w(r)}{k}} = \frac{r}{k^{q}}.
\end{equation*}
With respect to \ref{vmin_lab_N}, using that $s\mapsto \frac{C_1-s}{C_2-s}$ increases for $C_1>C_2>s$ we have 
\begin{equation*}
	N 	= 		\left\lceil\frac{r-k w(\rbase)-3k}{{\rup-k w(\rbase)-3k}} \right\rceil
		\geq		\frac{r}{\rup}
		\geq 	\frac{w^{-1}(\rbase)}{w^{-1}\pp{\frac{\rbase}{k}}}
		=		k^q,
\end{equation*}
where for the last inequality we used that $r\geq w^{-1}(\rbase)$ and that $\rup\leq w^{-1}\pp{\frac{\rbase}{k}}$ due to the definition of $\rbase, \rup$ in \eqref{vmin_eq_rup_rvert} and the monotonicity of $w(\cdot)$.
On the other hand, using that $r=w^{-1}\pp{\rbase+2w(\rbase)}$ and $\rup=w^{-1}\pp{\frac{\rbase}{k}-3}$, as well as that $w(\rbase)\geq 3$  from \ref{vmin_lab_three}, we obtain 
\begin{align*}
	N	&\leq 	\frac{r}{{\rup-k w(\rbase)-3k}}	+1\\
		&=		\pp{\frac{\cw\pp{\frac{\rbase}{k}-3}^q}{\cw\pp{\rbase+2w(\rbase)}^q}-\frac{k\pp{w(\rbase)+3}}{\cw\pp{\rbase+2w(\rbase)}^q}}^{-1}+1\\
		&=		k^q	\pp{\pp{1-\frac{3k+2w(\rbase)}{\rbase+2w(\rbase)}}^q-\frac{k^{q+1}\pp{w(\rbase)+3}}{\cw\pp{\rbase+2w(\rbase)}^q}}^{-1}+1\\
		&\leq		k^q	\pp{\pp{1-\frac{(k+2)w(\rbase)}{\rbase}}^q-\frac{2k^{q+1}w(\rbase)}{\cw\rbase^q}}^{-1}+1\\
		&\leq		k^q	\pp{1-q\frac{(k+2)w(\rbase)}{\rbase}-\frac{2k^{q+1}w(\rbase)}{\cw\rbase^q}}^{-1}+\frac{k^q}{2}
\end{align*}
where for the last inequality we used that $k\geq 2$, $q>1$.
In particular, $N\leq 2k^q$ holds if 
\begin{align*}
	\frac{(k+2)w(\rbase)}{\rbase}			&\leq 	\frac{1}{6q},	&&\text{	which is implied by}	&	\rbase^{\frac{q-1}{q}}\geq 6q(k+2)\cw^{-\frac{1}{q}},\\
	\frac{2k^{q+1}w(\rbase)}{\cw\rbase^q}	&\leq 	\frac{1}{6},	&&\text{	which is implied by}	&	\rbase^{\frac{q^2-1}{q}}\geq 12k^{q+1}\cw^{-\frac{q+1}{q}}.
\end{align*}
Combining this with \eqref{vmin_eq_C3} and \eqref{vmin_eq_C2}, we obtain that \ref{vmin_lab_rup}, \ref{vmin_lab_three}, and \ref{vmin_lab_N} hold if
\begin{equation}\label{vmin_eq_req_rbase}
	\rbase
	\geq		\max\Bp{3^q\cw, \pp{12 q k\cw^{-\frac{1}{q}}}^{\frac{q}{q-1}}}.
\end{equation}
This prompts us to fix $0<\delta<\frac{q-1}{q}$, $k_0\in\N$ and choose 
\begin{align}\label{vmin_eq_k}
		k=k(r)	&\coloneqq 	\left\lfloor k_0 \pp{\frac{\rbase(r)}{r_0}}^{\frac{q-1}{q}-\delta} \right\rfloor,
	&	k_0 		&\coloneqq 	\left\lfloor\frac{1}{12q}\cw^{\frac{1}{q}}r_0^{\frac{q-1}{q}}\right\rfloor,
\end{align}
where $k_0$ is the largest $k\in\N$ such that \eqref{vmin_eq_req_rbase} holds for $\rbase\pp{r_1}=r_0$ and hence $k(r_1)=k_0$.
The bound \eqref{sett_eq_assBox_r} for $r_0$ from Assumption~\ref{sett_aP_box} guarantees $k_0\geq 2$ and hence $k(r)\geq 2$ for all $r\geq r_1$.
Note that \eqref{sett_eq_assBox_r} also guarantees $\rbase(r)\geq r_0\geq 3^q\cw$ for $r\geq r_1$.

Hence, choosing $k(r)$ and $k_0$ as in \eqref{vmin_eq_k} yields \eqref{vmin_eq_req_rbase} and thus \ref{vmin_lab_rup}, \ref{vmin_lab_three}, and \ref{vmin_lab_N} for all $r\geq r_1$.
In addition, using that $k(r)\geq 2$ for all $r\geq r_1$, we obtain the bound
\begin{align}
	\label{vmin_eq_kboundlow1}
	k(r)		&\geq 	k_0\pp{\frac{\rbase(r)}{r_0}}^{\frac{q-1}{q}-\delta}-1\\
	\label{vmin_eq_kboundlow2}
			&\geq	\frac{k_0}{2}\pp{\frac{\rbase(r)}{\rbase(r_1)}}^{\frac{q-1}{q}-\delta}
			\geq		\frac{k_0}{2}\pp{\frac{w(r)}{w(r_1)}}^{\frac{q-1}{q}-\delta}
			=		\frac{k_0}{2}\pp{\frac{r}{r_1}}^{\frac{q-1}{q^2}-\frac{\delta}{q}},
\end{align}
where we used that, with the definition of $\rbase(r)$ from \eqref{vmin_eq_rup_rvert} and due to $\rbase\mapsto \frac{\rbase}{\rbase+2w(\rbase)}$ increasing, we have
\begin{equation*}
	\frac{\rbase(r)}{\rbase(r_1)}
	=		\frac{w(r)}{w(r_1)}\frac{\frac{\rbase(r)}{\rbase(r)+2w(\rbase(r))}}{\frac{\rbase(r_1)}{\rbase(r_1)+2w(\rbase(r_1))}}
	\geq		\frac{w(r)}{w(r_1)}.
\end{equation*}

\underline{Step 5: Quantifying the loss of speed per step.}
In Step 2, given $r\geq r_1$ and $v>0$ for the speed associated with the smaller boxes inside $\Qup{r}$, for the minimum average speed in $\Qup{r}$ we obtained
\begin{equation}
	\wt{v}	=			\wt{v}(v,r)
			\coloneqq	v\frac{r}{(N-1)\pp{\rup+\rbase}+\rup}.\tag{\ref{vmin_eq_def_v}}
\end{equation}
We bound the loss factor by bounding the excess distance from following the path along the smaller boxes. 
With \eqref{vmin_eq_N}, the definition of $N$, we have
\begin{align}
	(N-1)	&	\pp{\rup+\rbase}+\rup \notag\\
	&=
	(N-1)\pp{\rup-k\pp{w(\rbase)+3} + k\pp{w(\rbase)+3}} + (N-1)\rbase + \rup\notag\\
	&\leq 
	r+(N-1)k\pp{w(\rbase)+3} + (N-1)\rbase + \rup.\label{vmin_eq_extra-distances}
\end{align}
The three additional terms correspond to the overlap of the layers, the horizontally traversed distance, and the potential overshoot in the last layer.
For the overlap, using $\rbase\leq\rup$ from \ref{vmin_lab_rup}, the monotonicity of $w(\cdot)$, and the definition of $\rup$ from \eqref{vmin_eq_rup_rvert}, we obtain a bound corresponding to the horizontal distance,
\begin{align}\label{vmin_eq_extra-overlap}
	(N-1)k\pp{w(\rbase)+3} 	\leq		(N-1)k\pp{w(\rup)+3}		
							=	 	(N-1)\rbase.
\end{align}
For the horizontal distance, with \ref{vmin_lab_N}, the choice of $k$, and $\rbase\leq w(r)$, we have
\begin{align}
	(N-1)\rbase
	\leq		2k^q\rbase
	\leq		2k_0^q	\pp{\frac{\rbase}{r_0}}^{q-1-\delta q}\rbase
	\leq		r	~	2k_0^q	r_0^{\delta q-q+1} \cw^{-(1-\delta)}r^{-\delta}.\label{vmin_eq_extra-horizontal}
\end{align}
For the overshoot, first applying \ref{vmin_lab_rup} and then the bound on $k(r)$ from\eqref{vmin_eq_kboundlow2}, we obtain
\begin{align}\label{vmin_eq_overshoot}
	\rup 
	\leq 	r k^{-q}
	\leq 	r	~	 \pp{\frac{r_1}{r}}^{\frac{q-1}{q}-\delta}.
\end{align}
In summary, plugging \eqref{vmin_eq_extra-distances} with \eqref{vmin_eq_extra-overlap}, \eqref{vmin_eq_extra-horizontal}, and \eqref{vmin_eq_overshoot} into \eqref{vmin_eq_def_v} and then using that $k_0\leq \frac{1}{12q}\cw^{\frac{1}{q}}r_0^{\frac{q-1}{q}}$ yields
\begin{align}
	\wt{v}	
	&\geq		v\frac{1}{1 	+ 4 k_0^{q}r_0^{\delta q-(q-1)}\cw^{\delta-1}r^{-\delta}
							+ \pp{\frac{r_1}{r}}^{\frac{q-1}{q}-\delta} }
	\geq		v\frac{1}{1	+ \Cv r^{-\min\Bp{\delta, \frac{q-1}{q}-\delta}}}		\label{vmin_eq_vBound}	
\end{align}
with 
	$\Cv=\Cv(q,\cw,r_0,\delta)
	\coloneqq 
	\frac{4}{(12q)^q} r_0^{\delta q}\cw^{\delta} 
	+r_1^{\frac{q-1}{q}-\delta} $
where we used that $k_0\leq \frac{1}{12q}\cw^{\frac{1}{q}}r_0^{\frac{q-1}{q}}$.

\underline{Step 6: Ensuring improved probability per step.}
Referring to the probability bound \eqref{vmin_eq_ProbBound} from Step 3, for $0<p<1$ and $r\geq r_1$ with $\rbase(r)$, $N(r)$, and $k(r)$ as above, we set
\begin{equation}\label{vmin_eq_pbardef}
	\ov{p}	=			\ov{p}(p,r)	
			\coloneqq 	3 N(r) p^{k(r)}.
\end{equation}
In this step, we will provide bounds for $\ov{p}$.
In particular, we are interested in sufficient conditions to guarantee that $\ov{p}< p$.
Note that due to the bounds on $N$ from \ref{vmin_lab_N} it holds that
\begin{align}\label{vmin_eq_pbar}
	\ov{p}	\leq		6k^q p^k
			=		p \exp\pp{\log(6)+q\log(k)-(k-1)\log(p^{-1})}.
\end{align}
We obtain the desired bound if the exponent is negative.
For $p=p_0$ and $r=r_1$ with $k(r_1)=k_0$ as chosen in \eqref{vmin_eq_k}, we hence have $\ov{p}(p_0,r_1)<p_0$ if
\begin{align*}
	&&	0			&>	\log(6)+q\log(k_0)-(k_0-1)\log(p_0^{-1})\\
	\Leftarrow	&&	
	\log(p_0^{-1})	&> \frac{\log(6)+q\log(k_0)}{k_0-1},
\end{align*}
which is guaranteed by \eqref{sett_eq_assBox_p} from Assumption~\ref{sett_aP_box}.

Since $r\mapsto k(r)$ is increasing and we only want to apply the iteration for $r\geq r_1$ and $p\leq p_0$, we now only need to check that
\begin{equation*}
	g\colon	(p,k)\mapsto	\log(6)+q\log(k)-(k-1)\log(p^{-1}) 
\end{equation*}
is decreasing in $k$ for $k\geq k_0=k(r_1)$ and increasing in $p$ for $p\leq p_0$.
We have 
\begin{align*}
	\frac{\partial g}{\partial p}(p,k)	&=	\frac{k-1}{p}>0,\\
	\frac{\partial g}{\partial k}(p,k)	&=	\frac{q}{k}+\log(p)	\leq		\frac{1}{2}\log(p)		
	&&\Leftarrow
	&	\log(p^{-1})	&\geq	\frac{2q}{k},
\end{align*}
which holds for $p\leq p_0$ and $k\geq k_0$ due to \eqref{sett_eq_assBox_p} from Assumption~\ref{sett_aP_box}.
In particular, we have $g(p,k)<0$ for $p\leq p_0$ and $k\geq k_0$ and, due to the mean value theorem, we obtain
\begin{equation*}
	g(p,k)\leq g(p,k_0)+(k-k_0)\frac{\log(p)}{2}.
\end{equation*}
In summary, plugging these bounds into \eqref{vmin_eq_pbar} yields
\begin{align}
	\label{vmin_eq_pbarsmall}
	\ov{p}(p,r)	&\leq p \exp\pp{g(p,k(r))}< p,\\	
	\ov{p}(p,r)	&\leq p\exp\pp{g(p,k_0)+(k(r)-k_0)\frac{\log(p)}{2}}
	\label{vmin_eq_p_growth}
					\leq	p^{\frac{1}{2}\pp{k(r)-k_0}+1}
\end{align}
for all $r\geq r_1$ and $0< p\leq p_0$.

\underline{Step 7: The recursive formulas for $v_i$ and $p_i$.}
Assume that \eqref{vmin_eq_assBoxScaled} holds for some $i\in\Nz$, that is, assume that there exist $0\leq p_i\leq p_0$ and $v_i>0$ such that
\begin{equation*}
				\min\Bp{\PM{\Eup[v_i]{r}},\PM{\Er[v_i]{r}},\PM{\Ed[v_i]{r}},\PM{\El[v_i]{r}}}\geq 1-p_i,
				\quad\text{for all }	r\in[r_i,r_{i+1}).
\end{equation*}
Let $r\in [r_{i+1},r_{i+2})$. 
We want to apply the iterative scheme to obtain a probability bound for $\PM{\Eup[v_{i+1}]{r}}$ with some $v_{i+1}$ via \eqref{vmin_eq_ProbBound}.
By the definition of $(r_i)_i$, we have $\rbase(r)\in [r_i,r_{i+1})$. 
However, while $\rup(r)\geq r_i$ by \ref{vmin_lab_rup}, $\rup(r)>r_{i+1}$ for $r$ close enough to $r_{i+2}$.
Hence, we will first obtain bounds for $r\in [r_{i+1},r_{i+1}^{(1)})$, with $r_{i+1}^{(1)}$ chosen such that $\rup(r)<r_{i+1}^{(0)}\coloneqq r_{i+1}$ for all such $r$.
Having extended the range of the probability bounds, we repeat this step until we have obtained a bound for all $r< r_{i+2}$, see Figure~\ref{vmin_fig_steps}.
To be more precise, we set 
\begin{align*}
	r_{i+1}^{(0)}&\coloneqq r_{i+1},
	&
	r_{i+1}^{(j)}&\coloneqq	\inf\Bp{r\in[r_{i+1},\infty)	\,:\,	\rup(r)\geq r_{i+1}^{(j-1)}}
	\quad\text{for }j\in\N.
\end{align*}
The amount of such steps necessary to reach $r_{i+2}$ we denote with
\begin{equation}\label{vmin_eq_ell-def}
	\ell_{i+1}\coloneqq \min\Bp{j\in\N	\,:\,	r_{i+1}^{(j)}\geq r_{i+2}}.
\end{equation}
Now, since $\rbase(r),\rup(r)\in [r_{i},r_{i+1}^{(j-1)})$ for $r\in [r_{i+1}^{(j-1)},r_{i+1}^{(j)})$ by \ref{vmin_lab_rbase}, \ref{vmin_lab_rup} and the definition of $(r_{i+1}^{(j)})_j$, iteratively applying \eqref{vmin_eq_ProbBound} yields
\begin{align}
	&\PM{\Eup[v_{i}^{(j)}]{r}}	\geq 1-p_{i+1},	\quad\text{for all }	r\in[r_i,r_{i+1}^{(j)}\wedge r_{i+2}),~j\in\Bp{1,\ldots, \ell_{i+1}}, \notag\\
	&\begin{aligned}	\text{with }&
						&	v_{i}^{(j)}&\coloneqq	\inf\Bp{\wt{v}\pp{v_{i}^{(j-1)},r}	\,:\,	r\in \big[r_{i+1}^{(j-1)},r_{i+1}^{(j)}\big)} , 	
						&	v_{i}^{(0)}&\coloneqq	v_{i},\\
						&&	p_{i+1}		&\coloneqq	\sup\Bp{\ov{p}(p_i,r)	\,:\,	r\in [r_{i+1}, r_{i+2})},
	\end{aligned}\label{vmin_eq_vij-pi-bounds}
\end{align}
where $\wt{v}$ and $\ov{p}$ were defined in \eqref{vmin_eq_def_v} and \eqref{vmin_eq_pbardef}.
We also used that $\wt{v}(v,r)\leq v$ and $\ov{p}(p,r)<p$ for $r\geq r_1$ and $p\leq p_0$ by \eqref{vmin_eq_pbarsmall}.
We obtain the same for $\PM{\Er[v_{i}^{(j)}]{r}}$, $\PM{\Ed[v_{i}^{(j)}]{r}}$, and $\PM{\El[v_{i}^{(j)}]{r}}$ by respectively adjusting the procedure in the previous steps of the proof.
In summary, this yields \eqref{vmin_eq_assBoxScaled} with $(p_i)_i$ defined recursively as above and $(v_i)_i$ defined recursively via $v_{i+1}\coloneqq v_{i}^{(\ell_{i+1})}$.

\begin{figure}[h]
\centering
\includegraphics[width=0.7\textwidth]{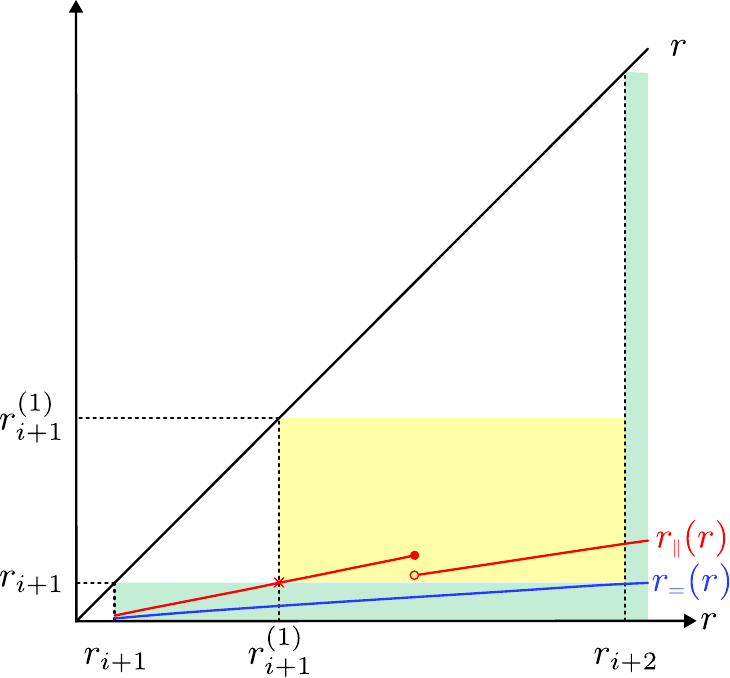}
\caption[An illustration of the functions {$r\mapsto (\rbase(r),\rup(r),r)$} with the thresholds from Step 7 of the proof of Proposition \ref{vmin_prop_assBoxScaled}]{	An illustration of the functions $r\mapsto (\rbase(r),\rup(r),r)$ with the thresholds from Step 7 of the proof of Proposition~\ref{vmin_prop_assBoxScaled}.
		Assuming we have bounds for boxes in the range $[r_{i},r_{i+1})$ (depicted in green), we can use these to obtain bounds for boxes of length $r\geq r_{i+1}$ (depicted in black) until $\rup(r)$ (depicted in red) leaves this range at $r=r_{i+1}^{(1)}$ (marked with a dotted line). 
		In the next step, these bounds can be used until $\rup(r)$ leaves the range $[r_{i+1},r_{i+1}^{(1)})$ (depicted in yellow) or until $\rbase(r)$ (depicted in blue) leaves the original range at $r=r_{i+2}$.
		This process is repeated until the second case occurs, which is when we can update the propagation probability for all smaller boxes and restart the loop.
		The discontinuous jumps in $\rup(r)$ correspond to increases of $k(r)\in\N$}
\label{vmin_fig_steps}
\end{figure}

For the probabilities, by \eqref{vmin_eq_ProbBound}, \eqref{vmin_eq_pbardef}, and \eqref{vmin_eq_p_growth}, we have
\begin{equation}\label{vmin_eq_pi-final-bound}
	p_{i+1}	\leq		\sup_{r\in [r_{i+1}, r_{i+2})} 3N(r)p_i^{k(r)}
			\leq		p_i^{\frac{1}{2}\pp{k(r_{i+1})-k_0}+1}
			\leq		p_i^{\frac{k_0}{2}\pp{\pp{\frac{r_i}{r_0}}^{\frac{q-1}{q}-\delta}-1}},
\end{equation}
using the lower bound for $k(r)$ from \eqref{vmin_eq_kboundlow1}.

It remains to prove the recursive bounds for $(v_{i})_i$. 
To achieve this, we first determine a uniform bound on the amount of steps $i\mapsto \ell_i$ for the respective intervals as defined in \eqref{vmin_eq_ell-def}.
To obtain such a bound we need to show that $j\mapsto r_{i+1}^{(j)}$ grows fast enough.
Note that, since $r\mapsto \rup(r)$ is right-continuous by the definition of $r\mapsto k(r)$, we obtain
	$\rup\pp{r_{i+1}^{(j)}} \geq r_{i+1}^{(j-1)}$
with the definition of $r_{i+1}^{(j)}$.
By the definition of $\rbase(r),\rup(r)$ from \eqref{vmin_eq_rup_rvert}, this implies
\begin{align*}
	w\pp{r_{i+1}^{(j)}}	&=	k\pp{r_{i+1}^{(j)}} \pp{w\pp{\rup\pp{r_{i+1}^{(j)}}} +3} + 2 w\pp{\rbase\pp{r_{i+1}^{(j)}}}
						\geq		k\pp{r_{i+1}^{(j)}}w\pp{r_{i+1}^{(j-1)}}.
\end{align*}
Combining this inequality with the lower bound for $k(r)\geq 2$ from \eqref{vmin_eq_kboundlow2}, we obtain 
\begin{equation*}
	r_{i+1}^{(j)}	\geq 	k^q\pp{r_{i+1}^{(j)}} r_{i+1}^{(j-1)}
					\geq		\max\Bp{2^q,\pp{\frac{r_{i+1}^{(j)}}{r_1}}^{\frac{q-1}{q}-\delta}} r_{i+1}^{(j-1)},
\end{equation*}
implying that $\ell_{i+1}$ is finite.
From this inequality we deduce that
\begin{align*}
	r_{i+1}^{(j)}	
	\geq		r_1^{1-\frac{q}{1+\delta q}} \pp{r_{i+1}^{(j-1)}}^{\frac{q}{1+\delta q}}
	&=		\pp{r_{i+1}^{(j-1)}}^{\frac{q}{1+\delta q}-\pp{\frac{q}{1+\delta q}-1}\frac{\log\pp{r_1}}{\log\pp{r_{i+1}^{(j-1)}}}}\\
	&\geq		\pp{r_{i+1}^{(j-1)}}^{\frac{q}{1+\delta q}-\pp{\frac{q}{1+\delta q}-1}\frac{\log\pp{r_1}}{\log\pp{r_{i+1}}}},
\end{align*}
where we used that $r_{i+1}^{(j)}\geq r_{i+1}$ for all $j\in\Nz$.
Iteratively applying this inequality yields
\begin{equation}\label{vmin_eq_rij-lower-bound}
	r_{i+1}^{(j)}	
	\geq		r_{i+1}^{\pp{\frac{q}{1+\delta q}-\pp{\frac{q}{1+\delta q}-1}\frac{\log\pp{r_1}}{\log\pp{r_{i+1}}}}^j},
\end{equation}
which provides a good lower bound for the growth of $j\mapsto r_{i+1}^{(j)}$ in terms of $j$.
In order to bound the amount of steps necessary to reach $r_{i+2}$, note that with $2w(\rbase)\leq 2w(\rup)\leq \rbase$ due to \ref{vmin_lab_rup} we have 
\begin{equation}\label{vmin_eq_ri-upper-bound}
	r_{i+2}	=		w^{-1}\pp{r_{i+1}+2w(r_{i+1})}
			\leq		w^{-1}\pp{2r_{i+1}}
			=		2^q\cw r_{i+1}^q
			=		r_{i+1}^{q+\frac{\log(2^q\cw)}{\log(r_{i+1})}}.
\end{equation}
Hence, by the definition of $\ell_{i+1}$ from \eqref{vmin_eq_ell-def} and by comparing the exponents in \eqref{vmin_eq_rij-lower-bound} and \eqref{vmin_eq_ri-upper-bound}, we know that 
\begin{equation*}
	\pp{\frac{q}{1+\delta q}-\pp{\frac{q}{1+\delta q}-1}\frac{\log\pp{r_1}}{\log\pp{r_{i+1}}}}^{\ell_{i+1}-1}
	< q+\frac{\log(2^q\cw)}{\log(r_{i+1})}
\end{equation*}
and can conclude that for any $i\geq 1$ we obtain the bound
\begin{align*}
	\ell_{i+1} 
	&\leq 
			\frac{\log\pp{q+\frac{\log(2^q\cw)}{\log(r_{i+1})}}}
			{\log\pp{\frac{q}{1+\delta q}-\pp{\frac{q}{1+\delta q}-1}\frac{\log\pp{r_1}}{\log\pp{r_{i+1}}}}}+1\\
	&\leq 
			\frac{\log\pp{q+\max\Bp{0,\frac{\log(2^q\cw)}{\log(r_{2})}}}}
			{\log\pp{\frac{q}{1+\delta q}-\pp{\frac{q}{1+\delta q}-1}\frac{\log\pp{r_1}}{\log\pp{r_{2}}}}}+1\\
	&\eqqcolon
	\ov{\ell} = \ov{\ell}(q,\cw,r_0,\delta).
\end{align*} 
That is, for all $i\in\N$ we have $\ell_i\leq \ell = \ell(q,\cw,r_0,\delta) \coloneqq \max\Bp{\ell_1,\ov{\ell}(q,\cw,r_0,\delta)}$.

Concluding, we obtain \eqref{vmin_eq_assBoxScaled} from \eqref{vmin_eq_vij-pi-bounds} with the iterative bounds on $(p_i)_i$ due to \eqref{vmin_eq_pi-final-bound} and the iterative bounds on $(v_i)_i$ due to $\ell$ applications of \eqref{vmin_eq_vBound}, given by
\begin{align*}
	v_{i+1} 	
	&\geq	\pp{\frac{1}{1	+ \Cv(q,\cw,r_0,\delta)r_{i+1}^{-\min\Bp{\delta,~ \frac{q-1}{q}-\delta}}}}^{\ell(q,\cw,r_0,\delta)} v_i,\\
	p_{i+1}
	&\leq	p_i^{\pp{\frac{r_i}{r_0}}^{\frac{q-1}{q}-\delta}-1}.
\end{align*}
\end{proof}

\section{Establishing an almost guaranteed effective minimum speed}\label{s_veff}
We will now show that, via Proposition~\ref{vmin_cor_assBoxScaled}, the box criterion \ref{sett_aP_box} implies that there is an almost guaranteed effective minimum speed in the sense of \ref{veff_aP_star}, see Section~\ref{subs_sett_veff}. 
That is, we need to find a family of events $(\Espeed{h}(M))_{h,M}$ for areas $M\subset\Rt$ and scales $h\gg 1$ which are \ref{veff_Pmeas} measurable with respect to a neighborhood of $M$, \ref{veff_Porder} ordered with respect to set inclusions, and, most importantly,  \ref{veff_Pspeed} for given $M\subset\Rt$ and scales $h\gg 1$  guarantees that all sufficiently $h$-fat sets with boundary in $M$ have a stable approximation which propagates with some effective minimum speed.
These notions have been introduced in Definition~\ref{veff_def_fatstab}.
In addition, we need to provide strong probability bounds in the form of \eqref{veff_eq_probBound}.

\begin{figure}[h]
\centering
\includegraphics[width=0.85\textwidth]{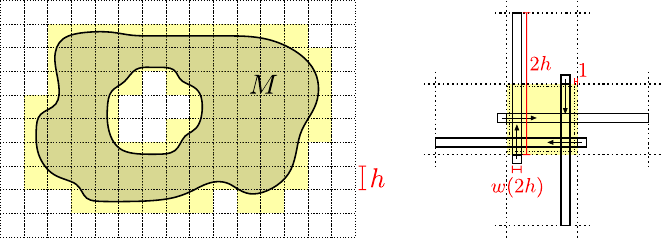}
\caption[An illustration of the event {${\Egrid[v]{h}(M)}$} introduced in Section \ref{s_veff}]{	An illustration of the event $\Egrid[v]{h}(M)$ introduced in Section~\ref{s_veff}.
			Left: For a set $M\subset\Rt$ (depicted in grey), we consider all cells of the grid $h\Zt$ which intersect $M$ (depicted in yellow).
			Right: If $\Egrid[v]{h}(M)$ holds, then for each of these cells and for each direction there exists a box along which the interface propagates cutting through that cell as well as through the neighboring cell in the respective direction
		}
\label{veff_fig_Egrid}
\end{figure}

Our candidate for $\Espeed{h}(M)$ is the event that there is a net of boxes with mesh width $\sim h$ covering $M$ along which the interface is guaranteed to propagates with some minimum average speed.
That is, given $h\geq \frac{r_0}{2}$, $M\subset\Rt$ and $v>0$, we define the event
\begin{align}
	\Egrid[v]{h}(M)	&\coloneqq \bigcap_{x\in h\Zt\,:\,\pp{x+[0,h]^2}\cap M\neq\emptyset} \Esquare[v]{h,x},\label{veff_eq_Egrid}
\intertext{where, with the measurability to be established in Lemma~\ref{veff_lem_meas},}	
	\Esquare[v]{h,x}	&\coloneqq			\pp{		\bigcup_{y\in x+\bp{1,h-w(2h)-1}\times\Bp{0}}		\Eup[v]{2h, y}	}
								\cap		\pp{		\bigcup_{y\in x+\Bp{0}\times\bp{1+w(2h),h-1}}		\Er[v]{2h, y}	}\notag\\
					&\quad		\cap		\pp{		\bigcup_{y\in x+\bp{1+w(2h),h-1}\times\Bp{h}}		\Ed[v]{2h, y}	}
								\cap		\pp{		\bigcup_{y\in x+\Bp{h}\times\bp{1,h-w(2h)-1}}		\El[v]{2h, y}	}.\notag
\end{align}
That is, if $\Egrid[v]{h}(M)$ holds, then for each cell of the grid $h\Zt$ touching $M$ and for each direction there exists a box along which the interface propagates crossing both this cell and then the neighboring cell in the corresponding direction, see Figure~\ref{veff_fig_Egrid}.

In this section, we will show that this event indeed guarantees 
\begin{itemize}
\item	an effective minimum speed of propagation on the scale of $h$, and
\item	approximate stability on the scale of $h$ for $2h$-fat sets.
\end{itemize} 
Due to  Corollary~\ref{vmin_cor_assBoxScaled} we will obtain stretched exponential bounds for the complement of $\Egrid[\vmin]{h}(M)$ as required in \eqref{veff_eq_probBound} for \ref{veff_aP_star}.
Finally, we will show that choosing $\Espeed{h}(M)\coloneqq \Egrid[v]{\frac{h}{2}}\pp{M+2\bp{-h,h}^2}$ satisfies the requirements for \ref{veff_aP_star} and hence we will establish Theorem~\ref{main_thm}.

\subsection{The scaled up version of the box criterion implies \ref{veff_aP_star}}

Regarding the measurability for \ref{veff_Pmeas} in \ref{veff_aP_star}, we obtain the following lemma analogously to Lemma~\ref{meas_lem_Eup} on the measurability of the event that the interface percolates along a single box.

\begin{lemma}[Measurability of ${\Egrid[v]{h}}(M)$]\label{veff_lem_meas}
Let $M\subset \Rt$, $h\geq \frac{r_0}{2}$, $v>0$. 
Then $\Egrid[v]{h}(M)$ is measurable with respect to $\F\pp{M+3[-h,h]^2}$.
\end{lemma}

\begin{proof}
We only need to show that  $\Esquare[v]{h,x}\in \F\pp{x+[-2h,3h]^2}$ for any $x\in h\Zt$ because if $\pp{x+[0,h]^2}\cap M\neq \emptyset$ then $\F\pp{x+[-2h,3h]^2}\subset \F\pp{M+3[-h,h]^2}$ .
As in the proof of Lemma~\ref{meas_lem_Eup} for the measurability of an event such as $\Eup[v]{2h}$, we apply Lemma~\ref{meas_lem_closure}.
That is, we have to show that the event is closed with respect to weighted uniform convergence of the coefficient field.
A small modification to the proof for Lemma~\ref{meas_lem_Eup} suffices, accounting for the fact that the box for which the respective event holds might vary along a given sequence of coefficient fields.

To be precise, let there be a sequence $(A_k,F_k)\rightarrow(A,F)$ converging uniformly on $\pp{x+[-2h,3h]^2}\times\Sone$ with 
\begin{equation*}
	(A_k,F_k)\in \Eup[v]{2h,y_k},\quad
	y_k\in x+\bp{1,h-w(2h)-1}\times\Bp{0}.
\end{equation*}
There is a (without loss of generality non-relabeled) subsequence with $y_k\rightarrow y$ for some $y\in x+\bp{1,h-w(2h)-1}\times\Bp{0}$.
Thus, $(A_k,F_k)_{y_k+\Qup{2h}}\rightarrow (A,F)_{y+\Qup{2h}}$ uniformly.
Analogously to the proof for Lemma~\ref{meas_lem_Eup}, we obtain $(A,F)\in\Eup[v]{2h,y}$.
An application of Lemma~\ref{meas_lem_closure} hence yields
\begin{equation*}
	\bigcup_{y\in x+\bp{1,h-w(2h)-1}\times\Bp{0}}		\Eup[v]{2h, y} \subset \F\pp{x+[-2h,3h]^2}.
\end{equation*}
The other directions in the definition of $\Esquare[v]{h,x}$ are treated analogously.
\end{proof}

Next, in view of \ref{veff_Pspeed} in \ref{veff_aP_star}, we show that ${\Egrid[v]{h}}(M)$ guarantees the existence of stable approximations as introduced in Definition~\ref{veff_def_fatstab} for sufficiently fat sets with boundary suitably covered by $M$.

\begin{lemma}[Approximate stability guaranteed by ${\Egrid[v]{h}}(M)$]\label{veff_lem_h-stability}
Let $h\geq \frac{r_0}{2}$ and $S\subset\Rt$ be a $\pp{\sqrt{2}+\frac{1}{2}}h$-fat set.
Assume that $(A,F)\in \Egrid[v]{h}\pp{\partial S + \pp{\sqrt{2}+\frac{1}{2}}B_h}$ for some velocity $0<v\leq\vmax$.
Then $S$ is 
	$\pp{h,\,12v^{-1}}$-approximately stable 
with respect to $(A,F)$.
\end{lemma}

\begin{proof}
We will construct a stable $\pp{h,\,12v^{-1}}$-approximation $S_h$ for $S$, which is a stable set with $S_h \supset \Rt\setminus \pp{\comp{S}+12\frac{\vmax}{v} B_h}$ and $S\hsubs \reach{12v^{-1}h}{S_h}$.

In order to relate the set to the grid cells used to define the events $\pp{\Egrid[v]{h}(M)}_M$, we select those inner cells which together with a small neighborhood are fully included in $S$, that is, 
\begin{equation}\label{veff_eq_def-int}
	\inner[h]{S}\coloneqq \bigcup	\Bp{		y+[0,h]^2	\,:\,	y\in h\Zt \text{ and }y+[0,h]^2+\frac{1}{2}B_h \subset S}.
\end{equation}
Note that in each ball of radius $\pp{\sqrt{2}+\frac{1}{2}}h$, the grid cell $y+[0,h]^2$ with $y\in h\Zt$ which includes the center of the ball is itself included in the ball with $\frac{h}{2}$ to spare and thus an inner cell of the ball, see Figure~\ref{veff_fig_int-ball}. 
\begin{figure}[h]
\centering
\includegraphics[width=0.35\textwidth]{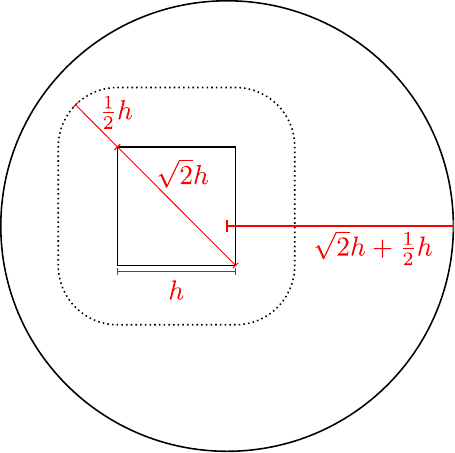}
\caption[An illustration of the fact that a sufficiently large ball contains an `inner cell' as defined in \eqref{veff_eq_def-int}]{	An illustration of the fact that a sufficiently large ball contains an `inner cell' as defined in \eqref{veff_eq_def-int}.
			The $\frac{1}{2}h$-neighborhood of a cell given by $x+[0,h]^2+\frac{1}{2}B_h$ (depicted with dotted lines) has diameter $\sqrt{2}h+h$.
			Hence, any ball with radius at least $\sqrt{2}h+\frac{1}{2}h$ and center in the cell fully includes this neighborhood. 
		}
\label{veff_fig_int-ball}
\end{figure}
In particular, since $S$ is $\pp{\sqrt{2}+\frac{1}{2}}h$-fat, we obtain
\begin{equation}\label{veff_eq_intEnlarge}
	S\subset	\inner[h]{S}+\pp{\sqrt{2}+\frac{1}{2}}\ov{B}_h.
\end{equation}
The grid cells making up $\inner[h]{S}$ can be split into two groups, see Figure~\ref{veff_fig_int-groups}:
The boundary cells are those, where at least one of the eight surrounding cells is not in $\inner[h]{S}$.
Accordingly, the interior cells are those with all surrounding cells also in $\inner[h]{S}$.
\begin{figure}[h]
\centering
\includegraphics[width=0.95\textwidth]{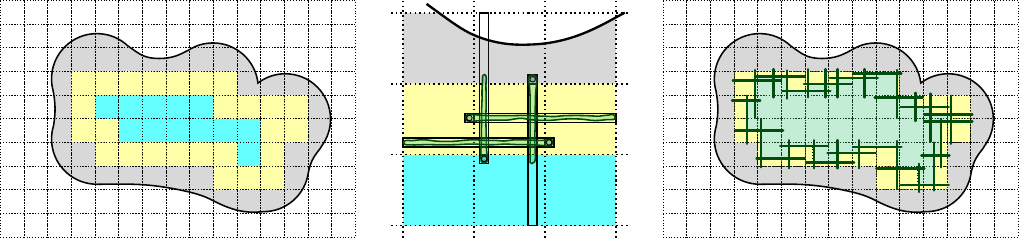}
\caption[An illustration of the construction of the stable approximation in the proof of Lemma \ref{veff_lem_h-stability}]{	An illustration of the construction of the stable approximation in the proof of Lemma~\ref{veff_lem_h-stability}.
			Left: For a fat set $S$, the boundary cells of $\inner[h]{S}$ (depicted in yellow) approximate $\partial S$ (depicted in black) from the inside, separating the interior cells of $\inner[h]{S}$ (depicted in blue) from the outside.
			Center: For each boundary cell and each direction, we take the stable `seed' guaranteed by $\Egrid[v]{h}\pp{\partial S + [-h,h]^2+\frac{1}{2}B_h}$.
			Its evolution (depicted in green) with respect to the coefficient field restricted to the respective box remains stable and cuts through the cell and a potential neighboring boundary cell.
			If the respective neighbor is not a boundary cell, we instead consider the set evolution at the moment it enters that cell.
			Right: The stable approximation is the union of all these sets with the surrounded interior (depicted in green)
			}
\label{veff_fig_int-groups}
\end{figure}
Note that this means that for $x\in h\Zt$ we have
\begin{align}\label{veff_eq_char-int}
	x+[0,h]^2 &\text{ is an interior cell of }\inner[h]{S}
	&&\Leftrightarrow
	&	x+[0,h]^2+[-h,h]^2+\frac{1}{2}B_h\subset S.
\end{align}
This means that 
\begin{align*}
	S \subset\bigcup\Bp{x+[0,h]^2\,:\, x+[0,h]^2 \text{ is an interior cell}}\cup \pp{\partial S+ 2[-h,h]^2+\frac{1}{2}B_h}
\end{align*}
and therefore
\begin{align}\label{veff_eq_intCells}
	\bigcup\Bp{x+[0,h]^2\,:\, x+[0,h]^2 \text{ is an interior cell}}\supset \Rt\setminus \pp{\comp{S}+ \pp{2\sqrt{2}+\frac{1}{2}}B_h}.
\end{align}
In fact, since our candidate $S_h$ for the stable approximation will include the union of the interior cells, this will ensure $S_h\supset \Rt\setminus \pp{\comp{S}+cv^{-1}\vmax  B_h}$ for $c\geq \pp{2\sqrt{2}+\frac{1}{2}}\vmax^{-1}$. 

To satisfy the other properties of a stable approximation, we now aim to find a stable interface surrounding the union of the interior cells inside $S$, which will evolve to envelope $S$.
Contrasting \eqref{veff_eq_char-int}, for each boundary cell $x+[0,h]^2\subset \inner[h]{S}$ we have
\begin{equation*}
	\pp{x+[0,h]^2}\cap\pp{\partial S + [-h,h]^2+\frac{1}{2}B_h}\neq\emptyset.
\end{equation*}
Thus, if $\Egrid[v]{h}\pp{\partial S + [-h,h]^2+\frac{1}{2}B_h}$ holds, then $\Esquare[v]{h,x}$ holds for the boundary cell $x+[0,h]^2$.
This means that for each direction there is a box crossing the cell and its neighbor along which the interface propagates.
To be more precise, there is a stable set (`seed') which within the respective translation and rotation of $\Qup{2h}+B_1$ would evolve to cross $x+[0,h]^2$ and then the corresponding neighboring cell.
By Lemma~\ref{stable_lem_preserv} and \ref{restr_lem_PropRestr}, these sets remain stable with respect to $(A,F)$ throughout this evolution.
For each of these sets, we stop the evolution when it leaves $\bigcup\Bp{y+[0,h]^2\,:\, y+[0,h]^2 \text{ is a boundary cell}}+B_{w(2h)+1}$, using that by  Corollary~\ref{stable_cor_HDcont} the evolution is continuous. 
By a straightforward adaption of Lemma~\ref{sti_lem_charEup} and in the sense of Definition \ref{sti_def_cut}, the resulting sets cut through $x+[0,h]^2$ in the respective direction, as well as through the next neighboring cell if that one is also a boundary cell, see Figure~\ref{veff_fig_int-groups}.
We take the union of all these evolved sets, denoted with $\wt{S}_h$.
By Lemma~\ref{stable_lem_union}, $\wt{S}_h$ is again stable as a locally finite union of stable sets.
Furthermore, using arguments as in the proofs of Lemma~\ref{sti_lem_base-form} and Lemma~\ref{sti_lem_extend}, we obtain that $\wt{S}_h$ separates the interior cells from the complement of $\inner[h]{S}$ in the sense that there exists no path in $\Rt\setminus\wt{S}_h$ connecting the two.

We now define $S_h$ as the union of $\wt{S}_h$ with all the path-connected components of $\comp{\pp{\wt{S}_h}}$ that intersect interior cells.
Note that $S_h$ is also stable since $\partial S_h \subset \partial \wt{S}_h$ with $S_h\supset \wt{S}_h$ and hence $\subJ \pp{1-\charfun[S_h]}(\cdot) \subset  \subJ \pp{1-\charfun[\wt{S}_h]}(\cdot)$ for the subjets in $\Rt$, see \cite[Definition~A.1]{FI25}. 
Moreover, due to \eqref{veff_eq_intCells}, we have $S_h \supset \Rt\setminus \pp{\comp{S}+ \pp{2\sqrt{2}+\frac{1}{2}}B_h}$.

It remains to show that $S\hsubs \reach{4\ceil{c+1}v^{-1}h}{S_h}$ for suitable $c>0$.
We obtain this result by consecutively applying the stitching lemma in the form of Lemma~\ref{sti_lem_base-form} to $S_h$  with time steps of size $2hv^{-1}$.
We know that for each of these time steps, $S_h$ evolving with respect to $(A,F)$ will propagate along the boxes associated to the surrounding cells which have been cut in the sense of Definition \ref{sti_def_cut} in the previous step, see Figure~\ref{veff_fig_int-expand} and \ref{veff_fig_int-expand-small}.

\begin{figure}[h]
\centering
\includegraphics[width=0.45\textwidth]{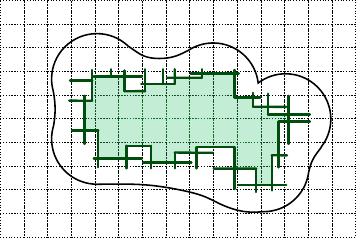}\hspace{10pt}
\includegraphics[width=0.45\textwidth]{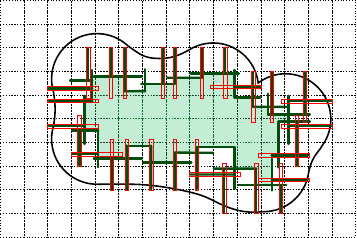}\\\vspace{15pt}
\includegraphics[width=0.45\textwidth]{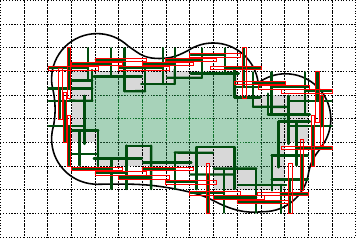}
\caption[An illustration of the first steps of the guaranteed expansion of the stable approximation in the proof of Lemma \ref{veff_lem_h-stability}]{	An illustration of the first steps of the guaranteed expansion of the stable approximation in the proof of Lemma~\ref{veff_lem_h-stability}.
			Top left: The constructed stable approximation $S_h$ (depicted in green) of the set $S$ (depicted with black boundary).
			Top right: After time steps $\Delta t=2v^{-1}h$, due to the stitching lemma, the stable approximation $S_h$ expands at least according to the evolution of the stable seeds along previously crossed boxes associated to the surrounding cells (depicted in red). 
			Bottom: After two such steps, all inner cells (depicted in grey) are $h$-enveloped.
			After every two additional cells, all horizontally, diagonally and vertically neighboring cells are $h$-enveloped. 
			In particular, the original set $S$  is $h$-enveloped after six steps }
\label{veff_fig_int-expand}
\end{figure}

\begin{figure}[h]
\centering
\includegraphics[width=0.23\textwidth]{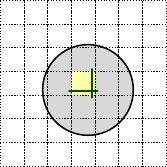}
\hspace{0pt}
\includegraphics[width=0.23\textwidth]{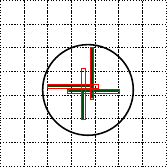}
\hspace{0pt}
\includegraphics[width=0.23\textwidth]{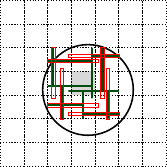}
\hspace{0pt}
\includegraphics[width=0.23\textwidth]{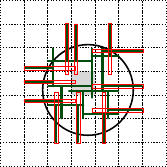}\\\vspace{10pt}
\includegraphics[width=0.23\textwidth]{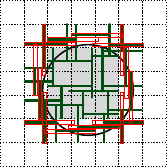}
\hspace{0pt}
\includegraphics[width=0.23\textwidth]{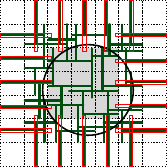}
\hspace{0pt}
\includegraphics[width=0.23\textwidth]{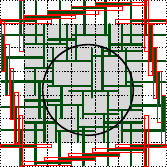}
\caption[An illustration of the guaranteed expansion of the stable approximation in the proof of Lemma \ref{veff_lem_h-stability} for a single inner cell]{	An illustration of the guaranteed expansion of the stable approximation in the proof of Lemma~\ref{veff_lem_h-stability} for a single inner cell.
			Top left: The stable approximation $S_h$ (depicted in green) constructed for a $\pp{\sqrt{2}+\frac{1}{2}}h$-fat set $S$ (depicted in grey) with $\inner[h]{S}$ consisting of a single inner cell (depicted in yellow).
			Center: After every time step $\Delta t=2v^{-1}h$, due to the stitching lemma, the stable approximation $S_h$ expands at least according to the evolution of the stable seeds along boxes associated to the surrounding cells which have been crossed in the previous step (depicted in red). 
			After two steps, $\inner[h]{S}$ is $h$-enveloped.
			Every two further steps, the set of cells which are guaranteed to be $h$-enveloped (depicted in grey) expands by all surrounding cells.
			Bottom right: After six steps, the evolution of $S_h$ is guaranteed to $h$-envelop $S$ (depicted with black boundary)
			}
\label{veff_fig_int-expand-small}
\end{figure}

More precisely, let $K_{h,0}=\inner[h]{S}$ be the union of all interior and boundary cells. 
For $i\in\N$, we define $K_{h,i+1}$ by adding to $K_{h,i}$ all cells $x+[0,h]^2$, which are neighboring $K_{h,i}$ horizontally, vertically, or diagonally, and for whose eight surrounding cells and itself $\Esquare[v]{h}$ holds.
Using Lemma~\ref{sti_lem_base-form}, $\reach{2v^{-1}h}{S_h}$ will horizontally (or vertically) cut through all cells that horizontally (or vertically) border $K_{h,0}$. 
Next, due to another application of Lemma~\ref{sti_lem_base-form}, $\reach{4v^{-1}h}{S_h}$ will orthogonally cut across those cells and their orthogonal neighbors, which includes all cells diagonally bordering $K_{h,0}$.
In summary, 
\begin{multline}\label{veff_eq_Kh0-cuts}
	\text{$\reach{4v^{-1}h}{S_h}$ cuts horizontally and vertically}\\
	\text{through any two neighboring cells in $K_{h,0}+[-h,h]^2$}.
\end{multline}
We claim that this implies $K_{h,0}\hsubs\reach{4v^{-1}h}{S_h}$.
Consider a cell in $K_{h,0}$ and its eight surrounding cells.
The nine cells are included in $\reach{4v^{-1}h}{S_h}+\ov{B}_h$ since otherwise there would be a straight segment in $\Rt\setminus \reach{4v^{-1}h}{S_h}$ crossing one of the boxes, contradicting \eqref{veff_eq_Kh0-cuts}.
Hence, for a point in the center cell to not be $h$-enveloped, it would have to be connected to the exterior of the nine cells in $\Rt\setminus \reach{4v^{-1}h}{S_h}$.
Assuming that this was the case, since $\Rt\setminus \reach{4v^{-1}h}{S_h}$ is open, there exists a path in $\Rt\setminus \reach{4v^{-1}h}{S_h}$ connecting the center cell to the boundary of the nine cells.
Without loss of generality, assume that the path ends in the bottom side of the nine cells, meaning that there exists a path segment connecting the top and bottom side of the bottom row of cells, see Figure~\ref{veff_fig_reach-paths}. 
Now there are two cases: 
First, if all three of the bottom cells are visited, there has to be a segment connecting the left and right side within the center cell, contradicting \eqref{veff_eq_Kh0-cuts}.
Second, the path segment connects the top and bottom side of two adjacent cells, contradicting \eqref{veff_eq_Kh0-cuts}.
Thus, the center cell has to be $h$-enveloped, which yields $K_{h,0}\hsubs\reach{4v^{-1}h}{S_h}$.
Repeating this argument, we obtain for $i\in\N$ that $K_{h,i-1} \hsubs \reach{4iv^{-1}h}{S_h}$.

\begin{figure}[h]
\centering
\includegraphics[width=0.66\textwidth]{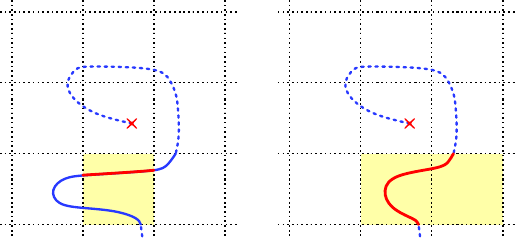}
\caption[A path connecting a grid cell to the exterior across the eight surrounding cells at some point has to cut either across a single cell or across two neighboring cells]{ A path connecting a grid cell to the exterior across the eight surrounding cells at some point has to cut either across a single cell or across two neighboring cells.
		 Left: If the last segment of the path crosses all three cells on the side of its exit, a sub-segment has to cut the center cell (depicted in yellow).
		 Right: Otherwise, the last segment cuts through the union of two neighboring cells (depicted in yellow)}
\label{veff_fig_reach-paths}
\end{figure}

To finish the proof, we show that $S\subset K_{h,2}$ and hence $S\hsubs \reach{12v^{-1}h}{S_h}$.
Recall that in \eqref{veff_eq_intEnlarge} we had $S\subset	\inner[h]{S}+\pp{\sqrt{2}+\frac{1}{2}}\ov{B}_h$ and thus in particular 
\begin{equation*}
	y\in \partial S + \pp{\sqrt{2}+\frac{1}{2}}\ov{B}_h	\qquad\text{ for all }y\in S\setminus \inner[h]{S}.
\end{equation*}
Since $\Egrid[v]{h}\pp{\partial S + \pp{\sqrt{2}+\frac{1}{2}}\ov{B}_h}$ holds, $\Esquare[v]{h}$ holds for all of the grid cells that intersect $S\setminus \inner[h]{S} = S\setminus K_{h,0}$ as well as their respective eight surrounding cells.
Since there are at most two cells necessary to connect $y\in S\setminus \inner[h]{S}$ to $K_{h,0}=\inner[h]{S}$, we obtain $S\subset K_{h,2}$.
\end{proof}

Continuing the argument from the proof above, we would obtain that the evolution of the constructed stable approximation $S_h$ essentially stays ahead of the original set expanding at some fixed speed. 
However, in order to obtain the effective minimum speed required for \ref{veff_Pspeed} in \ref{veff_aP_star}, ${\Egrid[v]{h}}(M)$ has to guarantee not just that $S_h$ expands accordingly in $M$, but also the set evolution $\reach{t}{S_h}$ at any time $t\geq 0$.

\begin{lemma}[An effective minimum speed of propagation guaranteed by ${\Egrid[v]{h}}(M)$]\label{veff_lem_h-speed}
Let $M\subset\Rt$, $h\geq \frac{r_0}{2}$, $v>0$.
Let $(A,F)\in{\Egrid[v]{h}}\pp{M+2[-h,h]^2}$ and $S\subset\Rt$ be a stable set with respect to $(A,F)$.
If $x\in \reach[(A,F)]{t}{S}\cap M$ for some $t>0$ with $\dist{x,S}\geq 3h$, then $x+\ov{B}_h\hsubs\reach[(A,F)]{t+14v^{-1}h}{S}$.
\end{lemma}

\begin{figure}[h]
\centering
\includegraphics[width=0.2\textwidth]{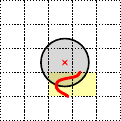}
\hspace{10pt}
\includegraphics[width=0.2\textwidth]{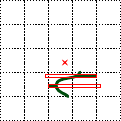}
\hspace{10pt}
\includegraphics[width=0.2\textwidth]{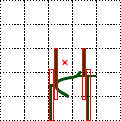}\\\vspace{15pt}
\includegraphics[width=0.2\textwidth]{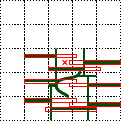}
\hspace{10pt}
\includegraphics[width=0.2\textwidth]{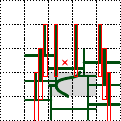}
\caption[An illustration of the guaranteed further expansion around a point $x$ reached by the set evolution for Case 1 in the proof of Lemma \ref{veff_lem_h-speed}]{ 	An illustration of the guaranteed further expansion around a point $x$ reached by the set evolution for Case 1 in the proof of Lemma~\ref{veff_lem_h-speed}.
		Top left: The stable set $\reach{t}{S}$ cuts vertically through the union of two neighboring cells below $x$, providing a base for the stitching lemma (depicted in red). 
		Center: Consecutive applications of the stitching lemma to the stable seeds in the boxes associated to the surrounding cells (depicted in red) guarantee the further expansion of the set evolution within time steps $\Delta t=2v^{-1}h$.
		Bottom right: After four time steps, the cells with the path segment (depicted in grey) are guaranteed to be $h$-enveloped.
		With every further time step, this set of guaranteed $h$-enveloped cells alternately expands by all vertically neighboring cells and all horizontally neighboring cells.
		Hence, the ball $\ov{B}_h(y)$ (depicted in black) is guaranteed to be $h$-enveloped after seven such time steps
		}
\label{veff_fig_reach-expand1}
\end{figure}

\begin{figure}[h]
\centering
\includegraphics[width=0.2\textwidth]{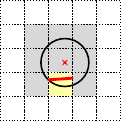}
\hspace{10pt}
\includegraphics[width=0.2\textwidth]{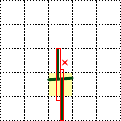}
\hspace{10pt}
\includegraphics[width=0.2\textwidth]{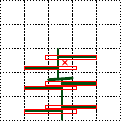}
\hspace{10pt}
\includegraphics[width=0.2\textwidth]{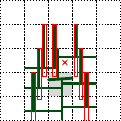}
\caption[An illustration of the guaranteed further expansion around a point $x$ reached by the set evolution for Case 2 in the proof of Lemma \ref{veff_lem_h-speed}]{	An illustration of the guaranteed further expansion around a point $x$ reached by the set evolution for Case 2 in the proof of Lemma~\ref{veff_lem_h-speed}.
		Left: The stable set $\reach{t}{S}$ cuts horizontally through the cell below $x$, providing a base for the stitching lemma (depicted in red). 
		Center: Consecutive applications of the stitching lemma to the stable seeds from ${\Egrid[v]{h}}\pp{M+2[-h,h]^2}$ in the boxes associated to the surrounding cells (depicted in red) guarantee the further expansion of the set evolution within time steps $\Delta t=2v^{-1}h$.
		Right: After three time steps, the cell with the path segment (depicted in grey) is guaranteed to be $h$-enveloped.
		With every further time step, this set of guaranteed $h$-enveloped cells alternately expands by all vertically neighboring cells and all horizontally neighboring cells.
		Hence, the ball $\ov{B}_h(y)$ (depicted in black) is guaranteed to be $h$-enveloped after six such time steps
			}
\label{veff_fig_reach-expand2}
\end{figure}

\begin{proof}
Let $x\in \reach[(A,F)]{t}{S}\cap M$ be in the grid cell $y+[0,h]^2$ with $y\in h\Zt$.
First, since $\dist{x,S}\geq 3h$, the initial set has no intersection with the $3\times 3$ group of cells $y+[-h,2h]^2$ surrounding $y+[0,h]^2$.
We will show that this implies that $\reach[(A,F)]{t}{S}$ cuts through the eight cells surrounding $y+[0,h]^2$ in the sense of Definition \ref{sti_def_cut}.
Second, since $(A,F)\in{\Egrid[v]{h}}\pp{M+2[-h,h]^2}$ holds, the event $\Esquare[v]{h}$ holds for all cells in the surrounding $5\times 5$ group $y+[-2h,3h]^2$. 
Starting with the cut from the first step, we will obtain the $h$-envelopment due to the guaranteed propagation along the boxes associated to ${\Egrid[v]{h}}\pp{M+2[-h,h]^2}$ in those cells.

Regarding the first step, we claim that $\reach[(A,F)]{t}{S}$ cuts vertically through the top or bottom row of $1\times 3$ cells in $y+[-h,2h]^2$ or horizontally through the left or right column of $3\times 1$ cells.
If this was not the case, then by definition there exist paths in $\Rt\setminus\reach[(A,F)]{t}{S}$ connecting the opposing narrow sides in all four triples. 
Combining these paths yields a loop which does not intersect $\reach[(A,F)]{t}{S}$ and hence neither $\reach[(A,F)]{s}{S}$ for any $s\in[0,t]$ due to the stability of $S$ and Lemma~\ref{stable_lem_preserv}.
The interior of the loop contains the center cell and is itself included in the surrounding $3\times 3$ group.
Therefore, the set evolution of $S$ does not intersect the parabolic boundary of the loop's interior until time $t$.
By the comparison principle -- as in the proof of Lemma \ref{sti_lem_charEup} -- this implies that $\reach[(A,F)]{t}{S}\cap \pp{y+[0,h]^2} = \emptyset$, contradicting the existence of $x\in \reach[(A,F)]{t}{S}\cap \pp{y+[0,h]^2}$.

Hence, without loss of generality, let $\reach[(A,F)]{t}{S}$ cut vertically through the bottom row of $1\times 3$ cells.
Now there are two cases, similar to Figure~\ref{veff_fig_reach-paths}:
In the first case, $\reach[(A,F)]{t}{S}$ cuts vertically through two adjacent cells, say without loss of generality through $x+[0,2h]\times[-h,0]$.
In the second case, $\reach[(A,F)]{t}{S}$ has to horizontally cut through the bottom center cell $x+[0,h]\times[-h,0]$.
If this was not the case, then there would exist paths in $\pp{x+[0,2h]\times[-h,0]}\setminus \reach[(A,F)]{t}{S}$ connecting left and right side of the two bottom left cells, connecting the top and bottom of the bottom center cell, and connecting the left and right side of the two bottom right cells.
Combining these three paths we could obtain a path connecting the left and right side of the full bottom row, contradicting the assumed vertical cut.

Regarding the second step, for both of these cases, we now consider the expansion of $\reach[(A,F)]{t}{S}$ along the boxes associated to ${\Egrid[v]{h}}\pp{M+2[-h,h]^2}$.
As in the proof of Lemma~\ref{veff_lem_h-stability}, we consecutively apply the stitching lemma in the form of Lemma~\ref{sti_lem_base-form} using the respective translations and rotations of $\Qup{2h}$ associated to ${\Egrid[v]{h}}\pp{M+2[-h,h]^2}$, which as observed above includes boxes for all cells in the surrounding $5\times 5$ grid $y+[-2h,3h]^2$.

In the first case, we obtain $y+[-h,2h]^2\hsubs \reach{7\cdot 2hv^{-1}}{S}$, see Figure~\ref{veff_fig_reach-expand1}.
In the second case, we obtain $y+[-h,2h]^2\hsubs \reach{6\cdot 2hv^{-1}}{S}$, see Figure~\ref{veff_fig_reach-expand2}.
In summary, using the stability of $S$, we obtain $x+\ov{B}_h\hsubs \reach{14v^{-1}h}{S}$.
\end{proof}

We obtain probability bounds for $\Egrid[\vmin]{h}(M)$ corresponding to \eqref{veff_eq_probBound} for \ref{veff_aP_star} based on the results from Section~\ref{s_vmin}, in particular Corollary~\ref{vmin_cor_assBoxScaled}.

\begin{corollary}\label{veff_cor_Egrid}
Assume that Assumption~\ref{sett_aP_stat}, Assumption~\ref{sett_aP_fin}, and the box criterion \ref{sett_aP_box} hold. 
Let $h\geq \frac{1}{2} r_0$.
For any $\ratebox \in \pp{0,\frac{q-1}{q^3}}$ and $\vmin=\vmin(q,\cw,r_0,v_0,\ratebox)$ as in Corollary~\ref{vmin_cor_assBoxScaled}, there exists $\Cgrid=\Cgrid(q,\cw,p_0,r_0,\ratebox)>0$ such that for any $M\subset\Rt$
\begin{equation*}
	\PM{\Egrid[\vmin]{h}(M)}
	\geq 
	1-16\pp{1+\frac{\diam(M)}{h}}^2\exp\pp{-\Cgrid h^{\ratebox+\frac{q-1}{q}}}.
\end{equation*}
\end{corollary}

\begin{figure}[h]
\centering
\includegraphics[width=0.35\textwidth]{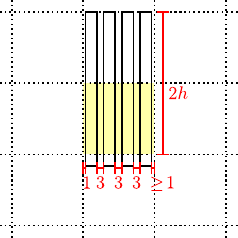}
\caption[An illustration of the box positions considered to bound the probability of the event {${\Egrid[v]{h}(M)}$}]{	 An illustration of the box positions considered to bound the probability of the event $\Egrid[v]{h}(M)$.
		 The distance between the boxes guarantees that the events of propagation along the respective boxes are independent.
		}
\label{veff_fig_Esquare-boxes}
\end{figure}

\begin{proof}
Note that  
\begin{equation*}
	\#\Bp{x\in h\Zt\,:\,\pp{x+[0,h]^2}\cap M\neq\emptyset}\leq	 \pp{2+\frac{\diam(M)}{h}}^2
\end{equation*}
and therefore by Assumption~\ref{sett_aP_stat} (Stationarity), it holds that
\begin{equation}\label{veff_eq_probEgrid-Esquare}
	\PM{\Egrid[\vmin]{h}(M)}\geq 1 - 4\pp{1+\frac{\diam(M)}{h}}^2 \pp{1- \PM{\Esquare[\vmin]{h,0}}}.
\end{equation}
With respect to bounding the probability of $\Esquare[v]{h,0}$, for each direction we use only finitely many box positions rather than along the whole interval.
However, if the distance between them is sufficiently large we can apply the finite range of dependence from Assumption~\ref{sett_aP_fin}.
That is, for boxes as in Figure~\ref{veff_fig_Esquare-boxes} we consider the union of the associated events,
\begin{equation*}
	\bigcup_{y\in \bp{1,h-w(2h)-1}\times\Bp{0}}	\Eup[v]{2h, y}
	\supset
	\bigcup_{y\in \Bp{1, 1+w(2h)+3, \ldots, 1+(j-1)(w(2h)+3)}\times\Bp{0}}	\Eup[v]{2h, y}
\end{equation*}
with
\begin{equation*}
	j	\coloneqq	\left\lfloor\frac{h+1}{w(2h)+3}\right\rfloor
		\geq			\frac{h}{4w(2h)}
		=			 \frac{\cw^{\frac{1}{q}}}{2^{\frac{1}{q}}4}h^{1-\frac{1}{q}}
\end{equation*}
where we used that $h+1\geq h \geq 4w(2h) \geq 2(w(2h)+3)$ since $2h\geq r_0$ and $w(r_0)\geq 3$ by \eqref{sett_eq_assBox_r} from Assumption~\ref{sett_aP_box}.
Hence, with the measurability from Lemma~\ref{meas_lem_Eup}, the finite range of dependence from Assumption~\ref{sett_aP_fin} and the stationarity from Assumption~\ref{sett_aP_stat}, using the bounds from Corollary~\ref{vmin_cor_assBoxScaled} we obtain
\begin{align*}
	\PM{\Esquare[\vmin]{h,0}}	&\geq	1 - \pp{1-\PM{\Eup[\vmin]{2h}}}^{j} - \pp{1-\PM{\Er[\vmin]{2h}}}^{j}\\
							&\hspace{30ex} - \pp{1-\PM{\Ed[\vmin]{2h}}}^{j} - \pp{1-\PM{\El[\vmin]{2h}}}^{j} \\
								&\geq	1 - 4\exp\pp{-\Cbox 2^\ratebox \cw^{\frac{1}{q}} \frac{1}{2^{\frac{1}{q}}4} 2^\ratebox h^{\ratebox+\frac{q-1}{q}}}.
\end{align*}
Plugging this bound into \eqref{veff_eq_probEgrid-Esquare} yields the result.
\end{proof}

We now obtain our main result by combining the previous considerations of the section.

\begin{theorem}[Explicit version of Theorem~\ref{main_thm}]\label{veff_thm_Ps}
Assume that Assumption~\ref{sett_aP_stat}, Assumption~\ref{sett_aP_fin}, and the box criterion \ref{sett_aP_box} hold. 
Let $\ratebox \in \pp{0,\frac{q-1}{q^3}}$.
Then \ref{veff_aP_star} holds with $h_0\coloneqq r_0$ and
\begin{align*}
		\ratespeed	&\coloneqq	\frac{q-1}{q}+\ratebox,
	&	\cspeed		&\coloneqq	\frac{\Cgrid(q,\cw,r_0,p_0,\ratebox)}{2^{\frac{q-1}{q}+\ratebox}},
	&	\Cspeed		&\coloneqq	2500,
	\\	\Cmeas		&\coloneqq	\frac{7}{\sqrt{2}},
	&	\vmineff		&\coloneqq	\frac{1}{16}\vmin\pp{v_0,q,\cw,r_0,p_0,\ratebox},
	&	\cstable		&\coloneqq	6\vmin^{-1},
\end{align*}
where $\Cgrid(q,\cw,p_0, r_0,\ratebox)$ is as in Corollary~\ref{veff_cor_Egrid} and $\vmin\pp{v_0,q,\cw,r_0,p_0,\ratebox}$ is  as in Corollary~\ref{vmin_cor_assBoxScaled}.
\end{theorem}

\begin{proof}
For the family of events $\big(\Espeed{\tilde{h}}(M)\big)_{\tilde{h},M}$ in \ref{veff_aP_star} which will guarantee the effective minimum speed on the respective scale $\tilde{h}\geq h_0$ in $M\subset\Rt$, we choose a slightly finer scale $h\coloneqq \frac{\tilde{h}}{2}$ and set 
\begin{equation*}
	\Espeed{\tilde{h}}(M)=\Egrid[\vmin]{h}\pp{M+4[-h,h]^2}
	\qquad\qquad\text{for $\tilde{h}\geq h_0\coloneqq r_0$.}
\end{equation*}
Let $\ratebox\in \pp{0,\frac{q-1}{q^3}}$.
Then the probability bound \eqref{veff_eq_probBound} is satisfied for all $\tilde{h}\geq h_0$ due to Corollary~\ref{veff_cor_Egrid} with
\begin{align*}
		\Cspeed		&\geq	16 (8\sqrt{2}+1)^2,
	&	\cspeed		&\coloneqq	\frac{\Cgrid(q,\cw,r_0,p_0,\ratebox)}{2^{\frac{q-1}{q}+\ratebox}},
	&	\ratespeed	&\coloneqq	\frac{q-1}{q}+\ratebox.
\end{align*}
Due to Lemma~\ref{veff_lem_meas}, using that $\F(M_1)\subset \F(M_2)$ for $M_1\subset M_2$ by the definition of $\F(\cdot)$ in \eqref{sett_eq_sigmaF}, we obtain 
\begin{equation*}
	\Espeed{h}(M)	\in 		\F\pp{M+7[-h,h]^2}
					\subset	\F\pp{M+7\sqrt{2}B_h}
					=		\F\pp{M+\frac{7\sqrt{2}}{2}B_{\tilde{h}}}
\end{equation*}
and hence the measurability in \ref{veff_Pmeas} is satisfied with $\Cmeas\coloneqq \frac{7}{\sqrt{2}}$.
In addition, the monotonicity of the events in \ref{veff_Porder} is satisfied since $\Egrid[\vmin]{h}(M_1)\supset \Egrid[\vmin]{h}(M_2)$ for any sets $M_1\subset M_2$ by the definition in \eqref{veff_eq_Egrid}.

Regarding \ref{veff_Pspeed}, let $S\subset\Rt$ be $\tilde{h}$-fat and thus $\pp{\sqrt{2}+\frac{1}{2}}h$-fat with $\partial S \subset M$ and let $(A,F)\in\Espeed{\tilde{h}}(M)$.
Then, due to Lemma~\ref{veff_lem_h-stability} we know that $S$ is $\pp{h,12\vmin^{-1}}$-approximately stable and thus $\pp{\tilde{h},6\vmin^{-1}}$-approximately stable.
Let $S_h$ be the stable approximation of $S$ as constructed in the proof of Lemma~\ref{veff_lem_h-stability}.

It remains to show that $(A,F)$ admits an effective minimum speed on the scale $\tilde{h}$ for $S_h$ in $M$.
Let $x\in\reach{t}{S_h}\cap M$ for some $t\geq 0$.
Now there are two cases, considering the $3\times 3$ block of cells from $h\Zt$ with $x$ in the central cell.
 
For the first case, assume that the intersection of this block with $S_h$ is empty. 
From this point on, following the arguments in the proof of Lemma~\ref{veff_lem_h-speed} we obtain
\begin{equation*}
	x+\ov{B}_h	\hsubs	\reach{t+14\vmin^{-1}h}{S_h}.
\end{equation*}
For the second case, assume that the intersection of the block with $S_h$ is non-empty.
Due to the construction of $S_h$, this means that either one of these nine cells or one of the neighboring cells is in $\inner[h]{S}$ as defined in the proof of Lemma~\ref{veff_lem_h-stability}.
With $\pp{K_{h,i}}_{i\in\Nz}$ denoting the same family of sets as in the proof of Lemma~\ref{veff_lem_h-stability}, we hence obtain $x\in K_{h,2}$ and
\begin{equation*}
	x+\ov{B}_h	\subset	K_{h,3} \hsubs	\reach{4\cdot 4\vmin^{-1}h}{S_h}	\subset	\reach{t+ 16\vmin^{-1}h}{S_h}.
\end{equation*}
Repeating the argument once more, we obtain
\begin{equation*}
	x+\ov{B}_{\tilde{h}}	=	x+\ov{B}_h+\ov{B}_h	\hsubs	\reach{t+16\vmin^{-1}h+16\vmin^{-1}h}{S_h}	
						=	\reach{t+16\vmin^{-1}\tilde{h}}{S_h}.	
\end{equation*}
In summary, \ref{veff_Pspeed} is satisfied with
\begin{align*}
		\cstable		&\coloneqq	6\vmin^{-1},
	&	\vmineff		&\coloneqq	\frac{1}{16}\vmin.
\end{align*}
\end{proof}

\section{Sparse Poisson obstacles admit an effective minimum speed}\label{s_poisson}
In this section, we will prove Theorem~\ref{intro_thm_Poisson}, that is, we will show that the box criterion \ref{sett_aP_box} is satisfied if the intensity $\rho$ of the Poisson point process $\Bp{X_i}_i$ is small enough compared to the uniform forcing $F_{uni}$, irrespective of the obstacle shape $\vphi\in C^{0,1}(\R)$ with $\supp \vphi \subset \ov{B}_1$. 

It is straightforward to see that, for all directions $e\in\Sone$, the characteristic function  $1-\charfun[S_t]$ of the set evolution 
\begin{align*}
	t\mapsto S_t\coloneqq \bigcup_{s\in[0,t]} \ov{B}_{\frac{2}{F_{uni}}}\pp{\frac{F_{uni}}{2}se}
\end{align*}
is a supersolution to the level-set formulation corresponding to the interface motion
\begin{equation*}
	V = H+F_{uni}.
\end{equation*}
The curvature $H$ is always at least $-\frac{F_{uni}}{2}$ and the normal velocity $V$ is at most $\frac{F_{uni}}{2}$.
In particular, $S_0$ is a stable set.

This means that if there is no obstacle in the $1$-neighborhood of a box with width at least $\frac{4}{F_{uni}}$, then a translation of $S_0$ is a stable set in the starting area which propagates to the opposite end of the box with average speed at least $v_0\coloneqq\frac{F_{uni}}{2}$.
That is, for all $r>0$ with $w(r)\geq \frac{4}{F_{uni}}$, using the definition of a Poisson point process we obtain
\begin{align}\label{poisson_eq_prob}
	\PM{\Eup[v_0]{r}} \geq \PM{\#\Bp{i\in\N\,:\, X_i\in\Qup{r}+\ov{B}_1}=0} = \exp\pp{-\rho \abs{\Qup{r}+\ov{B}_1}}
\end{align}
with $\#$ denoting the number of elements in a set and $\abs{\cdot}$ the area of a set.
Analogous bounds hold for $\PM{\Er[v_0]{r}}$, $\PM{\Ed[v_0]{r}}$,  and $\PM{\El[v_0]{r}}$.

For a given uniform forcing $F_{uni}>0$ and arbitrary width-scaling $q>1$, we choose the width function $w(r)=\pp{\frac{r}{\cw}}^{\frac{1}{q}}$ such that $w(r_0)= \frac{4}{F_{uni}}$ and such that the ratio of the sides of $\Qup{r_0}$ is the same for all $F_{uni}$, that is,
\begin{align}
	\cw	&\coloneqq 24q \pp{\frac{F_{uni}}{4}}^{q-1},	&
	r_0 &\coloneqq \frac{96q}{F_{uni}}.
\end{align}
The remaining degree of freedom has been used to further guarantee that
\begin{equation*}
	\pp{24q}^{\frac{q}{q-1}}\cw^{-\frac{1}{q-1}} 
	= \pp{24q}^{\frac{q-1}{q-1}}\frac{4}{F_{uni}}
	= r_0
\end{equation*}
and hence \eqref{sett_eq_assBox_r} from \ref{sett_aP_box} is automatically satisfied if in addition $w(r_0)\geq 3$, which corresponds to $F_{uni}\leq \frac{4}{3}$.

Due to the scaling invariance of the shape of $\Qup{r_0}$ with respect to $F_{uni}$, the required probability bounds \eqref{sett_eq_assBox_p} in \ref{sett_aP_box} for the first range of boxes are independent of $F_{uni}$, that is,
\begin{equation*}
	k_0 = \left\lfloor \frac{1}{12q}\cw^{\frac{1}{q}}r_0^{\frac{q-1}{q}} \right\rfloor
		= 2.
\end{equation*}
It remains to show the probability bound \eqref{sett_eq_assBox} for some $p_0<\min\Bp{\frac{2^{-q}}{6},e^{-q}}$.
Note that, for the sides of the the largest box $\Qup{r_1}$ in the relevant range, we have
\begin{align*}
	w(r_1)+2 	&= r_0+2w(r_0)+ 2 \leq 2r_0 = \frac{192q}{F_{uni}},\\
	r_1+2		&\leq w^{-1}(2r_0) = \frac{4(24q)^{q+1}2^q}{F_{uni}}.
\end{align*}
For $r\in[r_0,r_1)$, the bound \eqref{poisson_eq_prob} hence yields
\begin{multline*}
	\PM{\Eup[v_0]{r}}	\geq  \exp\pp{-\rho \abs{\Qup{r_1}+\ov{B}_1}}
						\geq	\exp\pp{-\rho\pp{w(r_1)+2}\pp{r_1+2}}\\
						\geq	\exp\pp{- 2^q 24^{q+1}q^{q+2}768\frac{\rho}{F_{uni}^2} },
\end{multline*}
where $\rho$ is the intensity of the Poisson point process.
With analogous bounds for $\PM{\Er[v_0]{r}}$, $\PM{\Ed[v_0]{r}}$,  and $\PM{\El[v_0]{r}}$, we thus obtain the box criterion if 
\begin{equation*}
	\rho	< 768^{-1}2^{-q} 24^{-q-1}q^{-q-2}\min\Bp{\frac{2^{-q}}{6},e^{-q}}F_{uni}^2 \eqqcolon c_0(q) F_{uni}^2
\end{equation*} 
and $F_{uni}\leq \frac{4}{3}$.
By the comparison principle, the box criterion is of course also satisfied for any stronger uniform forcings.
We have therefore established Theorem~\ref{intro_thm_Poisson} with the constants $C_0=2c_0^{-\frac{1}{2}}(q)$ and $\rho_0= \frac{16}{9}c_0(q)$ for any $q>1$.

This scaling relationship independent of the obstacle shape is obviously not preserved for arbitrarily large intensities -- it breaks at the very least when the critical threshold is reached where the support of the obstacles forms an infinite percolation cluster.

%

\sloppy
\printbibliography[title={References}]
\fussy
\end{document}